\documentclass[twoside, 11pt]{article}

\usepackage{amssymb}
\usepackage{amsfonts}
\usepackage{amsthm}
\usepackage{amsmath}
\usepackage{amscd}

\begin{document}

\pagestyle{myheadings}
\markboth{BASARAB}
{Residue Structures\/}
\title{{\huge Arithmetic-Arboreal Residue Structures
induced by Pr\" ufer Extensions : 
An axiomatic approach\/}}

\author{\c SERBAN A. BASARAB\\
Institute of Mathematics "Simion Stoilow" of the
Romanian Academy \\
P.O. Box 1--764\\
RO -- 70700 Bucharest 1, ROMANIA\\
\texttt{{\itshape e-mail}: Serban.Basarab@imar.ro}}
\date{} 

\maketitle
\vspace{8mm}

\newtheorem{te}{Theorem}[section]
\newtheorem{ex}[te]{Example}
\newtheorem{pr}[te]{Proposition}
\newtheorem{rem}[te]{Remark}
\newtheorem{rems}[te]{Remarks}
\newtheorem{co}[te]{Corollary}
\newtheorem{lem}[te]{Lemma}
\newtheorem{prob}[te]{Problem}
\newtheorem{exs}[te]{Examples}
\newtheorem{defs}[te]{Definitions}
\newtheorem{de}[te]{Definition}
\newtheorem{que}[te]{Question}
\newtheorem{aus}[te]{}


\newenvironment{aussage}%
\renewcommand{\theequation}{\alph{equation}}\setcounter{equation}{0}%

\def\R{\mathbb{R}}
\def\bbbm{{\rm I\!M}}
\def\N{\mathbb{N}}
\def\F{\mathbb{F}}
\def\H{\mathbb{H}}
\def\I{\mathbb{I}}
\def\K{\mathbb{K}}
\def\P{\mathbb{P}}
\def\D{\mathbb{D}}
\def\Q{\mathbb{Q}}
\def\Z{\mathbb{Z}}
\def\C{\mathbb{C}}
\def\T{\mathbb{T}}
\def\L{\mathbb{L}}

\def\cA{\mathcal{A}}
\def\cB{\mathcal{B}}
\def\cC{\mathcal{C}}
\def\cD{\mathcal{D}}
\def\cE{\mathcal{E}}
\def\cF{\mathcal{F}}
\def\cG{\mathcal{G}}
\def\cH{\mathcal{H}}
\def\cI{\mathcal{I}}
\def\cJ{\mathcal{J}}
\def\cK{\mathcal{K}}
\def\cL{\mathcal{L}}
\def\cM{\mathcal{M}}
\def\cN{\mathcal{N}}
\def\cO{\mathcal{O}}
\def\cP{\mathcal{P}}
\def\cQ{\mathcal{Q}}
\def\cR{\mathcal{R}}
\def\cS{\mathcal{S}}
\def\cT{\mathcal{T}}
\def\cU{\mathcal{U}}
\def\cV{\mathcal{V}}
\def\cW{\mathcal{W}}
\def\cX{\mathcal{X}}
\def\cY{\mathcal{Y}}
\def\cZ{\mathcal{Z}}
\def\epl{e^{+}}
\def\eb{e^{\bullet}}
\def\Ep{E^{+}}
\def\Eb{E^{\bullet}}


\newcommand{\notdiv}{{\not{|}\,}}
\newcommand{\ctg}{{\rm ctg\,}}
\newcommand{\sh}{{\rm sh\,}}
\newcommand{\ch}{{\rm ch\,}}
\newcommand{\Spek}{{\rm Spek\,}}
\newcommand{\Ker}{{\rm Ker\,}}
\newcommand{\di}{\mbox{\rm d\,}}
\newcommand{\cont}{\mbox{\rm cont\,}}
\newcommand{\rank}{\mbox{\rm rank\,}}
\newcommand{\codim}{\mbox{\rm codim\,}}
\newcommand{\ssqrt}[2]{\sqrt[\scriptstyle{#1}]{#2}}
\newcommand{\dbigcup}{\displaystyle\bigcup}
\newcommand{\dprod}{\displaystyle\prod}
\newcommand{\dbigcap}{\displaystyle\bigcap}
\newcommand{\doplus}{\displaystyle\bigoplus}
\newcommand{\dinf}{\displaystyle\inf}
\newcommand{\dint}{\displaystyle\int}
\newcommand{\dsum}{\displaystyle\sum}
\newcommand{\comb}[2]{\left(\begin{array}{c}#1\\#2\end{array}\right)}
\newcommand{\dlim}{\displaystyle\lim_{\stackrel{\longleftarrow}{}}}
\def\be{\begin{equation}}
\def\ee{\end{equation}}
\newcommand{\rs}{respectiv }
\newcommand{\G}{Z^1(\Gam,A)}

\def\bp{\begin{proof}}
\def\ep{\end{proof}}
\def\hf{\hfill $\square$}
\def\ben{\begin{enumerate}}
\def\een{\end{enumerate}}
\def\ba{\begin{eqnarray*}}
\def\ea{\end{eqnarray*}}

\def\ve{\varepsilon}
\def\ls{\leqslant}
\def\gs{\geqslant}
\def\lla{\longleftarrow}
\def\lra{\longrightarrow}
\def\Llra{\Longleftrightarrow}
\def\Lra{\Longrightarrow}
\def\p{\perp}
\def\wp{\,\widehat{1/p}}
\def\w4{\,\widehat{1/4}}
\newcommand{\wh}{\widehat}
\newcommand{\la}{\langle}
\newcommand{\ra}{\rangle}
\newcommand{\wt}{\widetilde}
\newcommand{\sm}{\setminus}
\newcommand{\sse}{\subseteq}
\newcommand{\es}{\varnothing} 

\newcommand{\f}{\frac}
\newcommand{\q}{\quad}
\newcommand{\n}{\vartriangleleft}


\newcommand{\al}{h }
\newcommand{\De}{\Delta}
\newcommand{\del}{\delta}
\newcommand{\eps}{\varepsilon}
\newcommand{\gam}{\gamma}
\newcommand{\Gam}{\Gamma}
\newcommand{\Lam}{\Lambda}
\newcommand{\lam}{\lambda}
\newcommand{\om}{\omega}
\newcommand{\Om}{\Omega}
\newcommand{\ovG}{\overline G}
\newcommand{\si}{\sigma}


\def\mA{\mathbb{A}}
\def\mC{\mathbb{C}}
\def\mN{\mathbb{N}}
\def\mQ{\mathbb{Q}}
\def\mR{\mathbb{R}}
\def\mZ{\mathbb{Z}}
\def\mF{\mathbb{F}}
\def\cA{\mathcal{A}}
\def\cC{\mathcal{C}}
\def\cD{\mathcal{D}}
\def\cE{\mathcal{E}}
\def\cF{\mathcal{F}}
\def\cG{\mathcal{G}}
\def\cM{\mathcal{M}}
\def\cP{\mathcal{P}}
\def\cS{\mathcal{S}}
\def\cO{\mathcal{O}}
\def\cL{\mathcal{L}}
\def\cQ{\mathcal{Q}}
\def\cK{\mathcal{K}}
\def\cH{\mathcal{H}}
\def\cX{\mathcal{X}}
\def\i{\rm i\/}

\begin{abstract}
We present an axiomatic framework for the residue structures
induced by Pr\" ufer extensions with a stress upon the intimate
connection between their arithmetic and arboreal theoretic properties.
The main result of the paper provides an adjunction relationship
between two naturally defined functors relating Pr\" ufer extensions
and superrigid directed commutative regular quasi-semirings. 
\smallskip


\smallskip

\end{abstract}

\bigskip


\section*{Introduction}

$\q$ Let $K$ be a valued field with Krull valuation $v$, valuation ring
$A$, maximal ideal $\bf{m}$, residue field $k := A/{\bf m}$,
and value group $\Lam := v K \cong K^\ast/A^\ast$. For any
$\lam \in \Lam$, we consider the $A$-submodules ${\bf m}_\lam :=
\{x \in K\,|\,v(x) > \lam\} \sse A_\lam := \{x \in K\,|\,v(x) \geq \lam\}$
of $K$, whence ${\bf m}_0 = {\bf m}, A_0 = A$. We denote by
$\cR_\lam := K/{\bf m}_\lam, R_\lam := K/A_\lam$ the associated
quotient $A$-modules. The disjoint union $\cR := \bigsqcup_{\lam \in \Lam}
\cR_\lam$ is identified with the set of all {\em open balls}
$\cB(a, \lam) := \{x \in K\,|\,v(x - a) > \lam\}$ for $a \in K,
\lam \in \Lam$, while the disjoint union $R := \bigsqcup_{\lam \in \Lam}
R_\lam$ is identified with the set of all {\em closed balls}
$B(a, \lam) := \{x \in K\,|\,v(x - a) \geq \lam\}$ for 
$a \in K, \lam \in \Lam$. Notice that $\cR = R$ if and only if
the maximal ideal ${\bf m} = \pi A$ is principal, so $v(\pi)$ is
the smallest positive element of the totally ordered Abelian
group $\Lam$.

The arithmetic and geometric structure of the valued field $K$
induces a suitable structure on the residue set $\cR \cup R$
which can be interpreted as a {\em deformation} of the original
structure of $K$. Fragments of this residue structure involving
suitable families of balls and relations between them play a
basic role in various interconnected algebraic-geometric and model theoretic 
contexts (elimination of quantifiers, cell decomposition, 
$p$-adic and motivic integration) [3, 8, 10-15, 17, 18, 20, 23]. 

Such residue structures can be also considered for more general
objects than valued fields. Thus we investigate in \cite{Res}
the residue structure $R$ associated to a Pr\" ufer domain
$A$ with field of fractions $K$, where the totally ordered 
value group $\Lam$ above is replaced by the Abelian $l$-group of 
the non-zero fractional ideals of finite type of $K$. In
particular, we introduce on $R$ an arboreal structure which
generalizes Tits' construction of the tree associated to
a valued field \cite{Trees}. The relation of this arboreal
structure with the arithmetic operations induced from $K$
as well as the action of $GL_2(K)$ are described.

In the present work we extend our paper \cite{Res} in two
directions : on the one hand, we consider the larger category
of {\em Pr\" ufer extensions} \cite{Grif,Kne}, and on the other hand we
develop an axiomatic framework for the {\em arithmetic-arboreal
residue structures} induced by such ring extensions.
Roughly speaking these residue structures are algebras
$(R, +, \bullet, -, {}^{- 1}, \eps)$ of signature (2, 2, 1, 1, 0)
satisfying a finite set of equational axioms which are mainly
suitable deformations of ring axioms. Thus these algebras,
called {\em directed commutative regular quasi-semirings} form a
variety. The main result of the paper (Theorem 7.2.) provides
a relationship of adjunction between two naturally defined
functors relating the categories (with suitably defined morphisms)
of Pr\" ufer extensions and of directed commutative regular
quasi-semirings satisfying an additional $\forall \exists$-axiom
called {\em superrigidity} \footnote{$^)$Not to be confounded with 
Margulis superrigidity.}$^)$.

Let us briefly review the content of this long and technical paper. In 
Section 1 some basic notions and constructions to be used later
in the paper are considered. Among them we mention the {\em commutative
$l$-monoid extension} $\wh{\Lam}$ of an Abelian $l$-group $\Lam$ - 
a nontrivial generalization of the totally ordered monoid $\Lam \cup
\{\infty\}$ associated to a totally ordered Abelian group $\Lam$, necessary 
to define the notion of {\em $l$-valuation} on a commutative ring
$w : B \to \wh{\Lam}$ - a generalization of the wellknown
notion of {\em valuation}. In particular, the basic notion of 
{\em Manis valuation} is generalized to the so called {\em Manis
$l$-valuation}. In the last subsection 1.4. of Section 1, some
basic facts on Pr\" ufer extensions are recalled, and the
{\em Pr\" ufer-Manis $l$-valuation} associated to a Pr\" ufer
extension is defined. 

The next four sections are devoted to the axiomatic framework
for the arithmetic-arboreal residue structures and the investigation 
of the main properties of the structures $(R, +, \bullet, -, {}^{- 1},
\eps)$ of signature $(2, 2, 1, 1, 0)$ introduced step by step
by suitable axioms. Thus, in Section 2, we introduce the class of 
{\em commutative regular semirings} which contains as proper
subclasses the Abelian $l$-groups and the {\em commutative regular rings} 
(in the sense of von Neumann), in particular, the fields and the 
Boolean rings. By relaxing suitably the distributive law, we define
the variety of {\em commutative regular quasi-semirings} which contains
the class of commutative regular semirings as a proper subvariety.  
Notice that in any commutative regular quasi-semiring $R$, the subset
$E^+$ of the idempotents of the commutative regular semigroup $(R, +)$
has a natural structure of Abelian $l$-group $\Lam$ 
with the group operation induced by $\bullet$ and the meet-semilattice
operation defined by $+$. Another interesting subvariety, 
whose members are the so called {\em directed commutative regular 
quasi-semirings} is introduced and studied in Section 3. Notice 
that the intersection of this subvariety with the subvariety of 
commutative regular semirings is the class of all Abelian $l$-groups.
Section 4 is devoted to the $\Lam$-metric structure of the directed
commutative regular quasi-semirings as well as to other related
questions (congruences, rigidity and superrigidity, subdirectly
irreducible structures). Using the $\Lam$-metric introduced in Section 4,
we consider in Section 5 two classes of structures with nice
arboreal theoretic properties : the class of {\em median directed commutative
regular quasi-semirings} and its proper subclass of {\em locally
faithfully full directed commutative regular quasi-semirings}.
The underlying arboreal structure on the latter objects is
a nontrivial generalization to Abelian $l$-groups $\Lam$ of 
the wellknown $\Lambda$-trees \cite{AB, MS}, where 
$\Lam$ is a totally ordered Abelian group. The classes of
structures above are also varieties in suitably extended signatures.

The last two sections of the paper are devoted to the study of
the relationship between some ring theoretic structures ($l$-valued
commutative rings, Pr\" ufer extensions) and their deformations
(the induced arithmetic-arboreal residue structures, axiomatized and
investigated in the previous sections). In Section 6 we associate a 
$l$-valued commutative ring $(B, w)$ to a superrigid directed 
commutative regular quasi-semiring $R$, and we characterize those
$R$ for which the associated ring $B$ is a Pr\" ufer extension of its
$l$-valuation subring (Corollary 6.7.). Thus we obtain a covariant
functor $\cB : \frak{R} \lra \frak{P}$ from an adequate category of directed
commutative regular quasi-semirings $\frak{R}$ to the category $\frak{P}$ of
Pr\" ufer extensions, with suitably defined morphisms. In Section 7
we define the {\em deformation} functor $\cR : \frak{P} \lra \frak{R}$
which is a left adjoint of the functor $\cB : \frak{R} \lra \frak{P}$
such that the endofunctor $\cB \circ \cR : \frak{P} \lra \frak{P}$
sends a Pr\" ufer extension $A \sse B$ to its completion $\wh{A} \sse \wh{B}$ 
(Theorem 7.2.). Consequently, the full subcategory $\frak{CP}$ of $\frak{P}$
consisting of all {\em complete Pr\" ufer extensions} is equivalent
with a suitable full subcategory $\frak{CR}$ of $\frak{R}$, whose objects
are explicitely described. 

\bigskip


\section{Notation and Preliminaries}

\subsection{The commutative $l$-monoid extension of an Abelian
$l$-group}

$\q$ Let $\Lam$ be an Abelian $l$-group with a multiplicative
group operation, the neutral element $\eps$, the partial
order $\leq$, and the (distributive) lattice operations 
$\wedge$ and $\vee$. For any $\alpha \in \Lam$, put 
$\alpha_+ := \alpha \vee \eps, \alpha_- := (\alpha^{- 1})_+, 
| \alpha | := \alpha \vee \alpha^{- 1} = \alpha_+ \alpha_-$. 
Let $\Lam_+ := \{\alpha\,|\,\alpha \geq \eps\}$ denote the
commutative $l$-monoid of all nonnegative elements of $\Lam$.
The Abelian $l$-group $\Lam$ has a canonical subdirect representation
into the product $\prod_{P \in \cP(\Lam)} \Lam/P$ of
its maximal totally ordered factors, where $P$ ranges over
the set $\cP(\Lam)$ of the minimal prime convex $l$-subgroups 
of $\Lam$, in bijection with the set of the minimal prime convex 
submonoids of $\Lam_+$ as well as with the set of the 
ultrafilters of the distributive lattice with a least element
$\Lam_+$.

In the present paper we shall use frequently the following
construction providing a natural embedding of an Abelian
$l$-group $\Lam$ into a commutative $l$-monoid $\wh{\Lam}$.

Setting $\alpha\downarrow := \{\gam \in \Lam\,|\,\gam \leq \alpha\}$
for any $\alpha \in \Lam$, we consider the inverse system
consisting of the distributive lattices (with a last element)
$\alpha \downarrow$ for $\alpha \in \Lam_+$, with the 
natural connecting epimorhisms $\beta \downarrow \lra
\alpha \downarrow, \gam \mapsto \gam \wedge \alpha$ for
$\alpha \leq \beta$. The inverse limit $\wh{\Lam} :=
\dlim \alpha\downarrow$ is a distributive lattice 
identified with the set of all maps 
$\varphi : \Lam_+ \lra \Lam$ satisfying $\varphi(\alpha)
= \varphi(\beta) \wedge \alpha$ for $\alpha \leq \beta$,
in particular, $\varphi(\alpha) \leq \alpha$,
with the induced pointwise partial order and lattice 
operations, and with $\om : \Lam_+ \lra \Lam,\alpha \mapsto \alpha$
as the last element. Notice that for any $\varphi \in \wh{\Lam}$, 
$\varphi(\alpha)_- = \varphi(\eps)^{- 1}$ for all $\alpha \in \Lam_+$,
and hence $\varphi(\alpha) = \varphi(\eps) \leq \eps$ whenever 
$\eps \leq \alpha \leq \varphi(\eps)^{- 1}$. Set
$\varphi_- := \varphi(\eps)^{- 1} \in \Lam_+$ for 
any $\varphi \in \wh{\Lam}$. 

The underlying unbounded distributive lattice of the
Abelian $l$-group $\Lam$ is identified with an ideal, 
in particular, with a sublattice of $\wh{\Lam}$, trough the 
embedding $\Lam \lra \wh{\Lam}, \gam \mapsto \iota_\gam$, 
where $\iota_\gam(\alpha) = \gam \wedge \alpha$
for all $\alpha \in \Lam_+$, so $\iota_\gam(\alpha) = \gam$ for
all $\alpha \geq \gam_+$. Notice that $\varphi \wedge \iota_\gam =
\iota_{\gam \wedge \varphi(\gam_+)}$ for all $\varphi \in \wh{\Lam},
\gam \in \Lam$, in particular, $\varphi \wedge \iota_\eps =
\iota_{\varphi(\eps)}$ for all $\varphi \in \wh{\Lam}$. In
the following we shall write $\gam$ instead of $\iota_\gam$ for
any $\gam \in \Lam$. $\Lam_+$ is identified with an ideal of 
the bounded distributive lattice 
$\wh{\Lam}_+ := \{\varphi \in \wh{\Lam}\,|\,\varphi \geq \eps\}
= \{\varphi \in \wh{\Lam}\,|\,\varphi(\eps) = \eps\}$ 
with the least element $\eps$ and the last element $\om$,
the inverse limit of the inverse system consisting of the
bounded distributive lattices $\alpha\downarrow_+ := [\eps, \alpha] =
\{\gam \in \Lam_+\,|\,\gam \leq \alpha\}$ for $\alpha \in \Lam_+$,
with the induced connecting epimorphisms. 

Let $\partial \wh{\Lam}_+$ denote the boolean subalgebra 
of the bounded distributive lattice $\wh{\Lam}_+$ consisting 
of those $\theta \in \wh{\Lam}_+$ which admit a (unique) 
complement $\neg \theta$, i.e. $\theta \wedge \neg \theta =
\eps$ and $\theta \vee \neg \theta = \om$.
Thus $\partial \wh{\Lam}_+$ is the inverse limit of the 
inverse system consisting of the boolean algebras 
$\partial [\eps, \alpha] := \{\gam \in [\eps, \alpha]\,|\,
\gam \wedge \frac{\alpha}{\gam} = \eps\}$ for $\alpha \in \Lam_+$,
with the natural connecting morphisms which, in general, are not 
necessarily surjective. Notice that $\partial \wh{\Lam}_+ \cap
\Lam_+ = \{\eps\}$, each $\theta \in \partial \wh{\Lam}_+$
extends uniquely to an endomorphism $\wt{\theta}$ of 
the $l$-group $\Lam$, and
$\wt{\theta} \circ \wt{\zeta} = \wt{\zeta} \circ \wt{\theta} =
\wt{\theta \wedge \zeta}$ for all 
$\theta, \zeta \in \partial \wh{\Lam}_+$, in particular,
$\wt{\theta} \circ \wt{\theta} = \wt{\theta}$ for all 
$\theta \in \partial \wh{\Lam}_+$. It follows that
$\wt{\theta}(\Lam) \sse \wt{\zeta}(\Lam) \Llra \theta \leq 
\zeta$, $\wt{\theta}(\Lam) \cap \wt{\zeta}(\Lam) = 
\wt{\theta \wedge \zeta}(\Lam)$, and 
$\wt{\theta \vee \zeta}(\Lam)$ is the convex $l$-subgroup
of $\Lam$ generated by $\wt{\theta}(\Lam) \cup \wt{\zeta}(\Lam)$
for all $\theta, \zeta \in \partial \wh{\Lam}_+$. In particular,
the Abelian $l$-group $\Lam = \wt{\om}(\Lam)$ is the direct 
sum of its convex $l$-subgroups $\wt{\theta}(\Lam)$ and
Ker$(\wt{\theta}) = \wt{\neg \theta}(\Lam)$ for any
$\theta \in \partial \wh{\Lam}_+$, where $\wt{\neg \theta}
(\gam) = \dfrac{\gam}{\wt{\theta}(\gam)}$ for all $\gam \in \Lam$.

The distributive lattice $\wh{\Lam}$ becomes a commutative 
$l$-monoid extending the Abelian $l$-group $\Lam$, with the  
multiplication defined by
$$(\varphi \psi)(\alpha) := (\varphi(\alpha \psi_-)
\psi(\alpha \varphi_-)) \wedge \alpha
= \displaystyle\lim_{\stackrel{\gam \rightarrow \infty}{}} 
(\varphi(\gam) \psi(\gam) \wedge \alpha)$$
for $\varphi, \psi \in \wh{\Lam}, \alpha \in \Lam_+$.
In particular, the identity $\varphi \psi = (\varphi \wedge \psi)
(\varphi \vee \psi)$ holds for all $\varphi, \psi \in \wh{\Lam}$.
The last element $\om$ of the distributive lattice
$\wh{\Lam}$ is a zero element of the monoid $\wh{\Lam}$,
i.e. $\om \varphi = \om$ for all $\varphi \in \wh{\Lam}$.
$\Lam$ is identified with the subgroup $\wh{\Lam}^\ast$ 
of all invertible elements of the monoid $\wh{\Lam}$, 
and $\wh{\Lam} = \Lam_+^{- 1} \wh{\Lam}_+$
is the monoid of fractions of the monoid $\wh{\Lam}_+$ 
relative to its submonoid $\Lam_+$. More precisely, any
element $\varphi \in \wh{\Lam}$ is uniquely represented
in the form $\varphi = \frac{\varphi_+}{\varphi_-}$
with $\varphi_+ \in \wh{\Lam}_+, \varphi_- \in \Lam_+$
subject to $\varphi_+ \wedge \varphi_- = \eps$ : 
$\varphi_+(\alpha) = \varphi(\alpha)_+$ for $\alpha \in \Lam_+$,
while $\varphi_- = \varphi(\eps)^{- 1}$ as defined above.
The idempotent elements of the monoid $\wh{\Lam}$ are 
exactly the elements of the boolean algebra $\partial \wh{\Lam}_+$, 
and $\varphi \psi = \varphi \vee \psi$ for all $\varphi, \psi \in 
\partial \wh{\Lam}_+$. Notice also that $\varphi \cdot \psi =
\om \Llra \varphi \vee \psi = \om$ for 
$\varphi, \psi \in \wh{\Lam}$; in particular, 
$\varphi^n = \om \Llra \varphi = \om$ for $\varphi \in
\wh{\Lam}, n \geq 1$.

If the Abelian $l$-group $\Lam$ is totally
ordered then $\wh{\Lam} = \Lam \cup \{\om\}$ and
$\partial \wh{\Lam}_+ = \{\eps, \om\}$, with 
$\eps = \om \Llra \Lam = \{\eps\}$.

\begin{ex} \em 
As a suggestive example, let $\Lam := \Q_{> 0}$ be the 
multiplicative Abelian group of the positive rationals,
freely generated by the subset $\P$ of all prime natural numbers.
$\Lam$ is an $l$-group with the lattice operations
$$x \wedge y := {\rm g.c.d.}(x, y) = \prod_{p \in \P} p^{{\rm min}
(v_p(x),v_p(y))},$$
and
$$x \vee y := {\rm l.c.m.}(x, y) = \prod_{p \in \P} p^{{\rm max}(v_p(x),
v_p(y))},$$
where $v_p$ denotes the $p$-adic valuation for $p \in \P$. One
checks easily that $\wh{\Lam}$ consists of all formal products
$\prod_{p \in \P} p^{n_p}$ with $n_p \in \Z \cup \{\infty\}$ such 
that the set $\{p \in \P\,|\,n_p < 0\}$ is finite (with the
corresponding map $\Lam_+ = \Z_{> 0} \lra \Lam,
m \mapsto \prod_{p \in \P} p^{{\rm min}(v_p(m), n_p)}$), while 
$\wh{\Lam}_+$ consists of those formal products with 
$n_p \in \N \cup \{\infty\}$ for all $p \in \P$, the so called 
{\em supernatural numbers} \cite{GCoh}. The boolean algebra
$\partial \wh{\Lam}_+$, identified with the power set 
$\mathbf{2}^\P$, consists of those formal products with
$n_p \in \{0, \infty\}$ for all $p \in \P$. 
\end{ex}

The next statements provide some basic properties of 
the correspondence $\Lam \mapsto \wh{\Lam}$.

\begin{lem}
Let $f : \Lam \lra \Gam$ be a morphism of Abelian
$l$-groups. Then $f$ extends canonically to a
morphism $\wh{f} : \wh{\Lam} \lra \wh{\Gam}$ of commutative
$l$-monoids provided one of the following conditions 
is satisfied.
\ben
\item[\rm (1)] $\Gam = {\rm ch}_\Gam(f(\Lam))$, the
{\em convex hull} of the image $f(\Lam)$ in $\Gam$; in
particular, if $f$ is onto.
\item[\rm (2)] $\Gam$ is totally ordered; in particular,
$\wh{f}$ is onto whenever $f$ is onto.
\een
\end{lem}

\bp
(1) Assuming that $\Gam = {\rm ch}(f(\Lam))$, let
$\varphi \in \wh{\Lam}, \gam \in \Gam_+$. By assumption
there is $\alpha \in \Lam_+$ such that $\gam \leq f(\alpha)$.
For any two such elements $\alpha, \beta$, we obtain 
$$f(\varphi(\alpha)) \wedge \gam = f(\varphi(\alpha
\vee \beta) \wedge \alpha) = f(\varphi(\alpha \vee \beta))
\wedge f(\alpha) \wedge \gam = f(\varphi(\alpha \vee \beta)) 
\wedge \gam = f(\varphi(\beta)) \wedge \gam,$$
therefore the element $\wh{f}(\varphi)(\gam) :=
f(\varphi(\alpha)) \wedge \gam \in \Gam$ for some
(for all) $\alpha$ as above is well defined. One
checks easily that the map $\wh{f}(\varphi) : \Gam_+ \lra \Gam$
so defined belongs to $\wh{\Gam}$, and moreover the
map $\wh{f} : \wh{\Lam} \lra \wh{\Gam}$ is a morphism
of commutative $l$-monoids extending $f$.

(2) Assume that $\Gam$ is totally ordered. If $f$ is
the null morphism then $\wh{f}(\varphi) = \eps$
for all $\varphi \in \wh{\Lam}$. If $f$ is nontrivial,
in particular, $\Gam \neq \{\eps\}$ and 
$\wh{\Gam} = \Gam \sqcup \{\om\}$, then for any $\varphi
\in \wh{\Lam}$ we distinguish the following two cases :

(i) There is $\alpha \in \Lam_+$ such that $f(\varphi(\alpha))
< f(\alpha)$. In this case, the element $\wh{f}(\varphi) :=
f(\varphi(\alpha)) \in f(\Lam) \sse \Gam$ does not depend on the
choice of $\alpha$ with the property above.

(ii) For all $\alpha \in \Lam_+$, $f(\varphi(\alpha)) = f(\alpha)$,
in particular, $\varphi \not \in \Lam$.
In this case we define $\wh{f}(\varphi) = \om$.

The map $\wh{f} : \wh{\Lam} \lra \wh{\Gam}$ is a morphism
of $l$-monoids extending $f$, with image $\wh{f(\Lam)} = f(\Lam)
\sqcup \{\om\}$.
\ep

\begin{lem}
Given a family $(\Lam_i)_{i \in I}$ of nontrivial
Abelian $l$-groups, let $\Lam := \prod_{i \in I} \Lam_i$.
Then $\wh{\Lam} \cong \prod_{i \in I} \wh{\Lam_i}$.
\end{lem}

\bp
Let us denote by $\pi_i : \Lam \lra \Lam_i, i \in I$, the 
natural projections. By Lemma 1.2, $\pi_i$ extends to
a morphism $\wh{\pi_i} : \wh{\Lam} \lra \wh{\Lam_i}$ of
$l$-monoids for all $i \in I$, and hence we obtain a
canonical morphism $\xi : \wh{\Lam} \lra \prod_{i \in I}
\wh{\Lam_i}, \varphi \mapsto (\wh{\pi_i}(\varphi))_{i \in I}$ 
that is the identity on $\Lam$. It follows that $\xi$ is
an isomorphism of $l$-monoids and 
$\xi^{- 1}((\psi_i)_{i \in I})(\alpha) = 
(\psi_i(\pi_i(\alpha)))_{i \in I}$ for $\psi_i \in \wh{\Lam_i},
i \in I$, and $\alpha \in \Lam_+$.
\ep

As a consequence of Lemmas 1.2, 1.3, we obtain

\begin{co}
Let $(\Lam_i)_{i \in I}$ be a family of nontrivial Abelian 
$l$-groups, and let \\ $f : \Gam \lra \Lam := \prod_{i \in I} \Lam_i$
be a morphism of Abelian $l$-groups such that 
${\rm ch}_{\Lam_i}(\pi_i(f(\Gam))) = \Lam_i$ for all $i \in I$.
Then $f$ extends canonically to a morphism $\wh{f} : \wh{\Gam} \lra
\wh{\Lam} \cong \prod_{i \in I} \wh{\Lam_i}$, and 
$\wh{f}(\wh{\Gam}) \sse \{\varphi \in \wh{\Lam}\,|\,
\forall \gam \in \Gam,\,\varphi \wedge f(\gam) \in f(\Gam)\}
\cong \wh{f(\Gam)}$.
\end{co}

\begin{rem} \em 
Let $\Lam$ be an Abelian $l$-group. For each $l$-subgroup
$\Gam \sse \Lam$, the set $\Om(\Lam\,|\,\Gam) := 
\{\varphi \in \wh{\Lam}\,|\,\forall \gam \in \Gam,\,
\varphi \wedge \gam \in \Gam\} = \{\varphi \in \wh{\Lam}\,|
\,\forall \gam \in \Gam_+,\,\varphi(\gam) \in \Gam\}$ is
a sub-$l$-monoid of $\wh{\Lam}$ with $\Om(\Lam\,|\,\Gam)^\ast = \Gam$, 
and the restriction map $r : \Om(\Lam\,|\,\Gam) \lra \wh{\Gam}, 
\varphi \mapsto \varphi|_{\Gam_+}$ is a well defined morphism 
of commutative $l$-monoids which is the identity on $\Gam$.
As in Corollary 1.4, if $\Lam = \prod_{i \in I} \Lam_i$ is a
product of nontrivial Abelian $l$-groups, and ${\rm ch}(\pi_i(\Gam))
= \Lam_i$ for all $i \in I$ then $r$ is an isomorphism. 
The simplest example with $r$ onto but not injective is obtained
by taking $\Lam$ a nontrivial Abelian $l$-group and $\Gam = \{\eps\}$,
in which case, $\Om(\Lam\,|\,\Gam) = \wh{\Lam}_+$ and $r$ is
the constant map $\eps$.
\end{rem}

\begin{co}
For each nontrivial Abelian $l$-group $\Lam$, its canonical 
subdirect representation into the product 
$\Gam := \prod_{P \in \cP(\Lam)} \Lam/P$ of its maximal 
totally ordered factors extends to a subdirect representation
of the associated commutative $l$-monoid $\wh{\Lam}$ into
the product $\prod_{P \in \cP(\Lam)} \wh{\Lam/P}$ of its 
maximal totally ordered factors, isomorphic to the commutative 
$l$-monoid $\wh{\Gam}$ associated to the Abelian $l$-group
$\Gam$. 

In particular, the boolean algebra 
$\partial \wh{\Lam}_+$ is canonically embedded
into the power set $\mathbf{2}^{\cP(\Lam)}$, and the
injective map $\cP(\Lam) \lra {\rm Max}(\partial \wh{\Lam}_+)$,
$$P \mapsto \{\theta \in \partial \wh{\Lam}_+\,|\,
\forall \alpha \in \Lam_+,\,\theta(\alpha) \in P\} =
\{\theta \in \partial \wh{\Lam}_+\,|\,\exists \alpha
\in \Lam_+, \frac{\alpha}{\theta(\alpha)} \not \in P\}$$ 
identifies $\cP(\Lam)$ with a dense subset of the profinite
space ${\rm Max}(\partial \wh{\Lam}_+)$, the Stone dual of
the boolean algebra $\partial \wh{\Lam}_+$, consisting of 
all maximal ideals of $\partial \wh{\Lam}_+$.
\end{co}

\begin{rem} \em
Let $M$ be a commutative {\em semilattice-ordered monoid} 
(for short, {\em sl-monoid}), i.e., a commutative monoid, 
with neutral element $\eps$, together
with a meet-semilattice operation $\wedge$ defining the partial
order $x \leq y \Llra x = x \wedge y$, which satisfies the 
compatibility condition $x (y \wedge z) = x y \wedge x z$
for all $x, y, z \in M$. Denote by $M^\ast$ 
the subgroup of all invertible elements of the monoid $M$,
and consider the submonoid 
$\wt{M} := \{x \in M\,|\,\exists y \in M, x y \leq \eps\}$.
Note that $M^\ast \cup M_- \sse \wt{M}$,
where $M_- := \{x \in M\,|\,x \leq \eps\}$, and $x \wedge y
\in \wt{M}$ for all $x \in \wt{M}, y \in M$, i.e.
$\wt{M}$ is a lower subset of $M$. It follows that 
$$\wt{M} = M^\ast \Llra \forall x \in M^\ast, \forall y \in M,
x \wedge y \in M^\ast,$$
so, in this case, $M^\ast$ is an Abelian $l$-group and the embedding
$M^\ast \lra \wh{M^\ast}$ extends uniquely to the morphism
of sl-monoids $\wh{w} : M \lra \wh{M^\ast}$, defined by
$\wh{w}(x)(\alpha) = x \wedge \alpha$ for
$x \in M, \alpha \in M^\ast_+$. For any such $sl$-monoid
$M$, the meet-semilattice operation $\wedge$ induces on 
the subset $\Eb(M) := \{x \in M\,|\,x^2 = x\}$ of idempotent 
elements a natural structure of distributive lattice 
with the least element $\eps$ and the join $x \vee y := x \cdot y$.
The distributive lattice $\Eb(M)$ is bounded, with the last
element $\om$, provided the element $\om$
is a zero as well as a last element of $M$. The morphism
$\wh{w} : M \lra \wh{M^\ast}$ induces by restriction 
a lattice morphism from the distributive lattice $\Eb(M)$ 
to the boolean algebra $\Eb(\wh{M^\ast}) = \partial \wh{M^\ast}_+$.
\end{rem}

The simplest way to extend the correspondence $\Lam \mapsto
\wh{\Lam}$ to a covariant functor is to define it on
the subcategory $\cC$ of the category of Abelian $l$-groups 
whose morphisms $f : \Gam \lra \Lam$ satisfy the condition 
${\rm ch}_\Lam(f(\Gam)) = \Lam$. On the other hand, let
$\cD$ be the category having as objects the $sl$-monoids
$M$ which satisfy the condition from Remark 1.7 : 
$x \in M, y \in M^\ast \Lra x \wedge y \in M^\ast$, with 
morphisms $F : M \lra N$ of $sl$-monoids for which 
${\rm ch}_{N^\ast}(F(M^\ast)) = N^\ast$.
By Lemma 1.2.(1), the correspondence $\Lam \mapsto \wh{\Lam}$
extends to a covariant functor $\,\wh{}\, : \cC \lra \cD$.
The next statement is immediate.

\begin{pr}
\ben
\item[\rm (1)] The functor $\,\wh{}\, : \cC \lra \cD$ is 
a right adjoint of the covariant functor 
$^\ast : \cD \lra \cC, M \mapsto M^\ast$.
\item[\rm (2)] The {\em counit} $\,^\ast \circ\, \wh{} \lra 1_\cC$
of the adjunction is a natural isomorphism, while the
{\em unit} $1_\cD \lra \wh{}\, \circ\, ^\ast$ is
the natural transformation sending each $M$ in $\cD$
to the canonical morphism $\wh{w} : M \lra \wh{M^\ast}$ as
defined in {\em Remark 1.7.} In particular, the category 
$\cC$ is equivalent with the full subcategory 
of $\cD$ consisting of those $M$ for which 
the canonical morphism $\wh{w} : M \lra \wh{M^\ast}$ 
is an isomorphism. 
\een
\end{pr}

\begin{rem} \em 
As a right adjoint, the functor 
$\,\wh{} : \cC \lra \cD$ is {\em continuous}, in particular, 
$\wh{\prod_{i \in I} \Lam_i} \cong \prod_{i \in I}
\wh{\Lam_i}$ for each finite family $(\Lam_i)_{i \in I}$ of 
Abelian $l$-groups. Though, by Lemma 1.3, the isomorphism above 
holds for arbitrary index sets $I$, we cannot deduce this 
fact from the continuity of the functor $\,\,\wh{}\,\,$ since
the category $\cC$ does not have arbitrary products. Indeed,
given a family $\Lam_i, i \in I$, of Abelian $l$-groups,
and a family $f_i : \Gam \lra \Lam_i$ of morphisms in $\cC$,
the canonical morphism $f : \Gam \lra \prod_{i \in I} \Lam_i$ of
Abelian $l$-groups is not necessarily a morphism in $\cC$.
However $f$ is a morphism in $\cC$ provided the index set $I$ is finite.
\end{rem}

\begin{lem}
 For any $\theta \in \partial \wh{\Lam}_+$, put $\wh{\Lam} \theta :=
\{\varphi \cdot \theta\,|\,\varphi \in \wh{\Lam}\} = \{\varphi \in
\wh{\Lam}\,|\,\varphi \cdot \theta = \varphi\}$, and let
$\wt{\theta}$ be the endomorphism of the Abelian $l$-group $\Lam$ 
induced by the endomorphism $\theta$ of the commutative $l$-monoid
$\Lam_+$, with $\Ker(\wt{\theta}) = \wt{\neg \theta}(\Lam) =
\{\frac{\alpha}{\wt{\theta}(\alpha}\,|\,\alpha \in \Lam\}$. Then :
\begin{enumerate}
\item[\rm (1)] For all $\varphi \in \wh{\Lam}$, the following 
assertions are equivalent.
\begin{enumerate}
 \item[\rm (i)] $\varphi \in \wh{\Lam} \theta$.
\item[\rm (ii)] $\varphi \cdot \neg \theta = \om$.
\item[\rm (iii)] $\wt{\theta} \circ \varphi = \theta$.
\item[\rm (iv)] $\theta \leq \varphi_+$.
\item[\rm (v)] $\varphi(\alpha) = \varphi(\frac{\alpha}{\theta(\alpha)})
\vee (\theta(\alpha) \cdot \varphi(\eps))$ for all $\alpha \in \Lam_+$.
\item[\rm (vi)] $\varphi(\alpha) = \theta(\alpha) \cdot 
\varphi(\frac{\alpha}{\theta(\alpha)})$ for all 
$\alpha \in \Lam_+$.
\end{enumerate}
Thus $\wh{\Lam} \theta = \{\varphi \in \wh{\Lam}\,|\,
\varphi_+ \geq \theta\}$ is a filter of $\wh{\Lam}$
and also a commutative $l$-monoid with neutral element 
$\theta$, and $(\wh{\Lam} \theta)_+ = \wh{\Lam}_+ \theta =$
$\{\varphi \in \wh{\Lam}\,|\,\varphi \geq \theta\}$. 
\item[\rm (2)] The restriction map 
$\varphi \mapsto \varphi\,|_{\Ker(\wt{\theta})_+}$ 
is an isomorphism of commutative
$l$-monoids $\wh{\Lam} \theta \lra \wh{\Ker(\wt{\theta})}$,
whose inverse sends $\psi \in \wh{\Ker(\wt{\theta})}$ to
$\varphi \in \wh{\Lam} \theta$ defined by $\varphi(\alpha) = 
\theta(\alpha) \cdot \psi(\frac{\alpha}{\theta(\alpha)})$ for
$\alpha \in \Lam_+$. In particular, $(\wh{\Lam} \theta)^\ast =
\Lam \theta \cong \Ker(\wt{\theta})$, and $\partial (\wh{\Lam}
\theta)_+ = (\partial \wh{\Lam}_+) \theta = \{\zeta \in \partial
\wh{\Lam}_+\,|\,\zeta \geq \theta\} \cong \partial 
\wh{\Ker(\wt{\theta})}_+$.
\item[\rm (3)] For $\theta, \zeta \in \partial \wh{\Lam}_+$,
$\wh{\Lam} \theta \sse \wh{\Lam} \zeta \Llra \theta \geq \zeta 
\Llra \theta = \theta \cdot \zeta$, $\wh{\Lam} \theta \cap 
\wh{\Lam} \zeta = \wh{\Lam} (\theta \cdot \zeta) = \wh{\Lam}
(\theta \vee \zeta)$, and $\wh{\Lam} \theta \vee \wh{\Lam} \zeta :=
\{\varphi \wedge \psi\,|\,\varphi \in \wh{\Lam} \theta,
\psi \in \wh{\Lam} \zeta\} = \wh{\Lam} (\theta \wedge \zeta)$.
In particular, $\wh{\Lam} \theta \cap \wh{\Lam} \neg \theta =
\wh{\Lam} \om = \{\om\}$, $\wh{\Lam} \theta \vee \wh{\Lam} \neg 
\theta = \wh{\Lam} \eps = \wh{\Lam}$, and the map
$\wh{\Lam} \lra \wh{\Lam} \theta \times \wh{\Lam} \neg \theta,
\varphi \mapsto (\varphi \cdot \theta, \varphi \cdot \neg \theta)$
is an isomorphism of commutative $l$-monoids, with the
inverse $(\varphi, \psi) \mapsto \varphi \wedge \psi$. 
\end{enumerate}
\end{lem}

\bp
(1). (i) $\Lra$ (ii) : $\varphi \cdot \neg \theta = \varphi \cdot \theta
\cdot \neg \theta = \varphi \cdot \om  = \om$.

(ii) $\Lra$ (iii) : Let $\alpha \in \Lam_+$. By assumption,
$$\alpha = \om(\alpha) = (\varphi \cdot \neg \theta)(\alpha)
\leq \varphi(\alpha)\, \neg \theta(\alpha \varphi(\eps)^{- 1}) =
\frac{\alpha \varphi(\alpha) \wt{\theta}(\varphi(\eps))}
{\theta(\alpha) \varphi(\eps)},$$
therefore $\theta(\alpha) \varphi(\eps) \leq \varphi(\alpha)
\wt{\theta}(\varphi(\eps))$. In particular, setting $\alpha = \eps$,
we obtain $\wt{\theta}(\varphi(\eps)) = \eps$, and hence
$\theta(\alpha) \varphi(\eps) \leq \varphi(\alpha) \leq \alpha$.
Applying the endomorphism $\wt{\theta}$ to the inequalities 
above, we obtain the required identity 
$\wt{\theta}(\varphi(\alpha)) = \theta(\alpha)$ for all
$\alpha \in \Lam_+$.

(iii) $\Lra$ (iv) : Let $\alpha \in \Lam_+$. We obtain
$$\theta(\alpha) = \wt{\theta}(\varphi(\alpha)) \leq 
\theta((\varphi(\alpha))_+) \leq \varphi_+(\alpha),$$
therefore $\theta \leq \varphi_+$ as desired. 

(iv) $\Lra$ (v) and (vi) : Since $\theta \leq \varphi_+$
by assumption, it follows that 
$$\varphi_+ = \theta \vee (\varphi_+ \wedge \neg \theta) =
\theta \cdot (\varphi_+ \wedge \neg \theta),$$
and hence 
$$\varphi_+(\alpha) = \theta(\alpha) \vee \varphi_+
(\frac{\alpha}{\theta(\alpha)}) = \theta(\alpha) \cdot 
\varphi_+(\frac{\alpha}{\theta(\alpha)})$$
for all $\alpha \in \Lam_+$. The assertions (v) and (vi)
are now immediate by multiplication with $\varphi(\eps) =
\varphi(\alpha)_-^{- 1} = \varphi(\frac{\alpha}{\theta(\alpha)})_-^{- 1}$. 

The implication (v) $\Lra$ (iv) is obvious.

(vi) $\Lra$ (i) : As $\varphi(\alpha) = \varphi(\eps)$
provided $\eps \leq \alpha \leq \varphi(\eps)^{- 1}$, it
follows by assumption that $\theta(\varphi(\eps)^{- 1}) = \eps$,
therefore 
$$(\varphi \cdot \theta)(\alpha) = \varphi(\alpha) \theta(\alpha
\varphi(\eps)^{- 1}) \wedge \alpha = \theta(\alpha)(\varphi(\alpha)
\wedge \frac{\alpha}{\theta(\alpha)}) = \theta(\alpha) \varphi
(\frac{\alpha}{\theta(\alpha)}) = \varphi(\alpha)$$
for all $\alpha \in \Lam_+$, i.e. $\varphi \cdot \theta = \varphi$
as desired.

(2) By (1), (i) $\Lra$ (iii), (i) $\Lra$ (vi), it follows that
the restriction map \\ $\wh{\Lam} \theta \lra \wh{\Ker(\wt{\theta})},
\varphi \mapsto \varphi|_{\Ker(\wt{\theta})_+}$ is an injective
morphism of commutative $l$-monoids. To prove that it 
is an isomorphism, it suffices by (1), (vi) $\Lra$ (i), to show
that for any $\psi \in \wh{\Ker(\wt{\theta})}$, the map 
$$\varphi : \Lam_+ \lra \Lam, \alpha \mapsto \theta(\alpha)
\psi(\frac{\alpha}{\theta(\alpha)})$$
belongs to $\wh{\Lam}$, i.e. $\varphi(\alpha) = \varphi(\beta)
\wedge \alpha$ for all $\alpha, \beta \in \Lam_+$ such that 
$\alpha \leq \beta$. Since $\psi \in \wh{\Ker(\wt{\theta})}$, 
it follows that $\wt{\theta}(\varphi(\alpha)) = \theta(\alpha),
\wt{\neg \theta}(\varphi(\alpha)) = \psi(\frac{\alpha}{\theta(\alpha)})$
for all $\alpha \in \Lam_+$. Consequently, assuming that 
$\eps \leq \alpha \leq \beta$, we obtain 
$$\wt{\theta}(\frac{\varphi(\beta) \wedge \alpha}{\varphi(\alpha)}) =
\wt{\neg \theta}(\frac{\varphi(\beta) \wedge \alpha}{\varphi(\alpha)}) =
\eps,$$
and hence the required identity.

The proof of the statement (3) is straightforward.
\ep

\begin{co}
$\overline{\Lam} := \{\gam \cdot \theta\,|\,\gam \in \Lam, \theta \in 
\partial \wh{\Lam}_+\}$ is the smallest
$l$-submonoid of $\wh{\Lam}$ containing the union
$\Lam \cup \partial \wh{\Lam}_+$, and 
$\overline{\Lam}_+ = \{\gam \cdot \theta\,|\,\gam
\in \Lam_+, \theta \in \partial \wh{\Lam}_+\}$.
In addition, the commutative monoid
$\overline{\Lam}$ is regular \footnote{$^)$A commutative
(multiplicative) semigroup $S$ is {\em regular} if for all
$x \in S$ there is $y \in S$ such that $x^2 y = x$. A regular 
commutative semigroup $S$ is an {\em inverse} semigroup, i.e.
for all $x \in S$ there is a unique element $x^{- 1} \in S$,
called the {\em quasi-inverse} of $x$, satisfying the identities
$x^2 x^{- 1} = x, (x^{- 1})^2 x = x^{- 1}$.}$^)$, 
with the quasi-inverse naturally defined by 
$(\gam \cdot \theta)^{- 1} := \gam^{- 1} \cdot \theta$ 
for $\gam \in \Lam, \theta \in \partial \wh{\Lam}_+$.
\end{co}

\bp
By Lemma 1.10.(2), the surjective map $\Lam \times \partial \wh{\Lam}_+
\lra \overline{\Lam}, (\gam, \theta) \mapsto \gam \cdot \theta$
identifies $\overline{\Lam}$ with the disjoint union
$\bigsqcup_{\theta \in \partial \wh{\Lam}} \Ker(\wt{\theta})$
of convex $l$-subgroups of $\Lam$. For 
$\theta, \zeta \in \partial \wh{\Lam}_+, \gam \in 
\Ker(\wt{\theta}), \rho \in \Ker(\wt{\zeta})$, we obtain
$$(\gam \theta) \cdot (\rho \zeta) = (\gam \rho) \cdot (\theta
\zeta) = (\frac{\gam \rho}{\wt{\zeta}(\gam) \wt{\theta}(\rho)}) 
\cdot (\theta \zeta),$$
$$\gam \theta \leq \rho \zeta \Llra \theta \leq \zeta, \gam \leq 
\rho \zeta(\gam_+) \Llra \theta \leq \zeta, (\neg \zeta)(\gam_+)
\leq \rho_+, \rho_- \leq \gam_-,$$
$$(\gam \theta) \wedge (\rho \zeta) = (\gam \theta(\rho_+) 
\wedge \rho \zeta(\gam_+)) \cdot (\theta \wedge \zeta),$$
and
$$(\gam \theta) \vee (\rho \zeta) = (\frac{\gam}{\zeta(\gam_+)} 
\vee \frac{\rho}{\theta(\rho_+)}) \cdot (\theta \zeta).$$
\ep

For instance, in the case $\Lam = \Q_{>0}$ as in Example 1.1,
$\overline{\Lam}$ consists of those formal products 
$ x:= \prod_{p \in \P} p^{n_p}$, with $n_p \in \Z \cup \{\infty\}$,
for which the set $\{p \in \P\,|\,n_p \in \Z \sm \{0\}\}$ is
finite, and $x^{- 1} = \prod_{p \in \P} p^{m_p}$, where
$m_p = - n_p$ if $n_p \in \Z$, and $m_p = \infty$ if $n_p = \infty$.


\subsection{The $sl$-monoid associated to a commutative ring extension}

$\q$ Let $A \sse B$ be a commutative unital ring extension. We denote by
$M := M(A, B)$ the set of all finitely generated (in short f.g.)
$A$-submodules of $B$. $M$ becomes a commutative multiplicative 
$sl$-monoid with the usual multiplication and 
the meet-semilattice operation $I \wedge J := I + J$ induced
by the partial order $I \leq J \Llra J \sse I$. The neutral
element is $A$, while $\{0\}$ is a zero element as well as 
the last element of the $sl$-monoid $M$. The $sl$-submonoid
$M_+$ of the nonnegative elements of $M$ consists of the
f.g. ideals of $A$.

The subgroup $M^\ast$ of the invertible elements of the 
monoid $M$ is exactly the Abelian multiplicative group 
of those $A$-submodules $I \sse B$ which are $B$-{\em invertible}, 
i.e. $I J = A$ for some $A$-submodule $J$ of $B$. Indeed, every
$B$-invertible $A$-submodule $I$ is f.g., and the $A$-submodule
$J$ satisfying $I J = A$ is unique, i.e. 
$J = I^{- 1} := [A : I] := \{x \in B\,|\,x I \sse A\}$, and
f.g. In particular, $A x \in M^\ast$ and $(Ax)^{- 1} = A x^{- 1}$
for all $x \in B^\ast$, so the factor group $B^\ast/A^\ast$
is identified with a subgroup of $M^\ast$. 
With the notation from Remark 1.7, the lower submonoid 
$\wt{M} := \{I \in M\,|\,\exists J \in M, A \sse I J\}$ of $M$, 
containing the ordered subgroup $M^\ast$ and 
the lower submonoid $M_- := \{I \in M\,|\,A \sse I\}$, 
consists of those f.g. $A$-submodules $I$ of $B$ which
are $B$-{\em regular}, i.e. $I B = B$. 

\bigskip 


\subsection{$l$-valuations on commutative rings}

$\q$ Let $B$ be a commutative unital ring and $\Lam$ a (multiplicative)
Abelian $l$-group, extended as in 1.1. to the commutative $l$-monoid
$\wh{\Lam}$, with its boolean algebra $\partial \wh{\Lam}_+$ of
idempotent elements.

\begin{de}
 A map $w : B \lra \wh{\Lam}$ is called a $l$-{\em valuation}
whenever the following conditions are satisfied.
\begin{enumerate}
 \item[\rm (1)] $w(x y) = w(x) w(y)$ for all $x, y \in B$.
\item[\rm (2)] $w(x + y) \geq w(x) \wedge w(y)$ for all $x, y \in B$.
\item[\rm (3)] $w(1) = \eps$ and $w(0) = \om$.
\end{enumerate}
\end{de}

From the axioms above we deduce that $w(- x) = w(x)$ for all $x \in B$,
$w(B^\ast)$ is a subgroup of $\Lam$, and the map 
$\Eb(B) := \{x \in B\,|\,x^2 = x\} \lra \partial \wh{\Lam}_+, 
x \mapsto w(1 - x)$ is a morphism of boolean algebras.

In particular, if $\Lam = \{\eps\}$ is trivial, so $\wh{\Lam} = \Lam$,
i.e. $\om = \eps$, then the $l$-valuation $w$ is the constant map 
$x \mapsto \eps$. If $\Lam$ is a nontrivial totally ordered Abelian
group, so $\wh{\Lam} = \Lam \sqcup \{\om\}$, then we obtain the usual 
notion of {\em valuation} \cite[VI.3.1.]{Bou}, \cite[I.1.]{Kne}

Some notions and basic facts about valuations 
extend to $l$-valuations as follows. Let $w : B \lra \wh{\Lam}$
be a $l$-valuation, and assume that $\Lam \neq \{\eps\}$.

The set $w^{- 1}(\om)$ is a radical ideal of $B$. 
It is called the {\em support} of $w$ and is denoted by supp$(w)$. 
The subring $A_w := w^{- 1}(\wh{\Lam}_+)$ of $B$ is called the 
{\em $l$-valuation ring} of $w$. Thus supp$(w)$ is a radical ideal 
of the ring $A_w$ too, contained in the {\em conductor} of
$A_w$ in $B$, the biggest ideal $\{x \in B\,|\,B x \sse A_w\}$
of $B$ contained in $A_w$. 

Setting $\overline{B} := B/{\rm supp}(w)$,
with the canonical projection $\pi : B \lra \overline{B}$, there 
exists a unique $l$-valuation $\overline{w} : \overline{B} \lra
\wh{\Lam}$ such that $\overline{w} \circ \pi = w$, with 
supp$(\overline{w}) = \{0\}$, so the factor ring $\overline{B}$ 
is reduced, and $A_{\overline{w}} = \overline{A_w} := 
A_w/{\rm supp}(w)$. In particular, if $w$ is a valuation
then supp$(w)$ is a prime ideal, ${\bf p}_w := \{x \in B\,|\,w(x)
> \eps\}$ is a prime ideal of $A_w$, the {\em center} of $w$,
and $\overline{w}$ extends uniquely to a Krull valuation $\wt{w}$
on the quotient field of the domain $\overline{B}$, so 
$A_{\overline{w}} \sse A_{\wt{w}}, {\bf p}_{\wt{w}} \cap A_{\overline{w}} =
{\bf p}_{\overline{w}} = {\bf p}_w/{\rm supp}(w)$.

Consider the subdirect representation
$\wh{\Lam} \lra \prod_{P \in \cP(\Lam)} \wh{\Lam/P}$,
where $\cP(\Lam)$ consists of all minimal prime convex $l$-subgroups
of $\Lam$, cf. Corollary 1.6.. Composing the 
$l$-valuation $w : B \lra \wh{\Lam}$ with the 
projection $\wh{\Lam} \lra \wh{\Lam/P} = (\Lam/P) \sqcup
\{\om_P\}$ for $P \in \cP(\Lam)$, we obtain a family of
valuations $(w_P : B \lra \wh{\Lam/P})_{P \in \cP(\Lam)}$.
It follows that supp$(w) = \bigcap_{P \in \cP(\Lam)}$
supp$(w_P)$, $A_w = \bigcap_{P \in \cP(\Lam)} A_{w_P}$, 
and $\overline{B}$ is identified to a subdirect product 
of domains through the injective ring morphism 
$\overline{B} \lra  \prod_{P \in \cP(\Lam)} \overline{B_P}$, 
where $\overline{B_P} := B/{\rm supp}(w_P)$ for $P \in \cP(\Lam)$.

For any $l$-valuation $w : B \lra \wh{\Lam}$ with $A := A_w$,
let $M := M(A, B)$ be the commutative $sl$-monoid associated 
as in 1.2 to the commutative ring extension $A \sse B$.
The map $\wh{w} : M \lra \wh{\Lam}$ sending a f.g. 
$A$-submodule $I = \sum_{1 \leq i \leq n} A x_i \sse B$ to
the element $\wh{w}(I) := \bigwedge_{1 \leq i \leq n} w(x_i) 
\in \wh{\Lam}$ is a well-defined morphism of $sl$-monoids with
a last element. Its image 
$$\wh{w}(M) = \{\bigwedge_{1 \leq i \leq n}
\varphi_i\,|\,n \in \N, \varphi_i \in w(B)\}$$ 
is a $sl$-submonoid of $\wh{\Lam}$, generated as meet-semilattice
by $w(B)$, with 
$$\wh{w}(M)_+ = \wh{w}(M_+) = \{\bigwedge_{1 \leq i \leq n}
\varphi_i\,|\,n \in \N, \varphi_i \in w(A)\} \sse \wh{\Lam}_+,$$ 
$$\wh{w}(M)_- = \wh{w}(M_-) = \{\bigwedge_{1 \leq i \leq n} 
(\varphi_i)_-^{- 1}\,|\,n \in \N, \varphi_i \in 
w(B \sm A)\} \sse \Lam_-.$$
The subset $U := \wh{w}^{- 1}(\Lam) \sse M$ is closed under
multiplication, and $M^\ast \sse \wt{M} \sse U$. Notice that
$I \leq J\,({\rm i.e.}\,J \sse I) \Llra  \wh{w}(I) \leq \wh{w}(J)$
for all $I \in M^\ast, J \in M$, i.e. $I = \{x \in B\,|\,w(x)
\geq \wh{w}(I)\}$ for all $I \in M^\ast$. Indeed, a implication 
holds for all $I, J \in M$, while, 
assuming $I \in M^\ast$ and $\wh{w}(I) \leq \wh{w}(J)$,
it follows that $\eps \leq \wh{w}(I^{- 1} J)$, and hence
$I^{- 1} J \sse A$, so $J \sse I$ as desired. Consequently,
$M^\ast$ is torsion free, and the morphism 
$\wh{w} : M \lra \wh{\Lam}$ induces by restriction a 
monomorphism of ordered groups 
$\wh{w}|_{M^\ast} : M^\ast \lra \Lam$.

Notice also that  
$I \in U, J \in M, J \leq I \Lra J \in U$ since $\Lam$ 
is a lower subset of $\wh{\Lam}$, therefore
the monoid of fractions ${\bf M}_w := U^{- 1} M$
is a $sl$-monoid with $[\frac{I}{J}] \wedge [\frac{I'}{J'}] =
[\frac{I J' + I' J}{J J'}]$, and ${\bf M}_w^\ast$ is a 
lower subset of ${\bf M}_w$ and hence a $l$-group.
Denote also by $\wh{w}$ the extended morphism of $sl$-monoids
${\bf M}_w \lra \wh{\Lam}, [\frac{I}{J}] \mapsto \wh{w}(I) \wh{w}(J)^{- 1}$.

\begin{de}
The image $\wh{w}({\bf M}_w) = \{\varphi \in \wh{\Lam}\,|
\,\exists \psi \in \wh{w}(M) \cap \Lam, \varphi \psi \in \wh{w}(M)\}$
is called the {\em value $sl$-monoid} of the $l$-valuation $w$ and is
denoted by $\cM_w$.
The $l$-subgroup $\cM_w^\ast = \cM_w \cap \Lam = \wh{w}({\bf M}_w^\ast)$ 
of $\Lam$ is called the {\em value $l$-group} of $w$.
\end{de}

Notice that $\cM_w^\ast$ contains the $l$-subgroup generated by 
$w(B) \cap \Lam$ as well as the $l$-subgroup generated by the subset
$\{\varphi_- = \varphi(\eps)^{- 1}\,|\,\varphi \in w(B \sm A_w)\}$.
In particular, if $w$ is a valuation then $\cM_w^\ast = \{\varphi
\cdot \psi^{- 1}\,|\,\varphi, \psi \in w(B) \cap \Lam\}$ is the
totally ordered value group of $w$. If in addition $w$ is a Krull
valuation, i.e. $B$ is a field, then ${\bf M}_w^\ast \cong \cM_w^\ast$.

The coarsening and equivalence relations for valuations extends to 
$l$-valuations as follows.

\begin{de}
 A $l$-valuation $w' : B \lra \wh{\Lam'}$ is called {\em coarser}
than the $l$-valuation $w : B \lra \wh{\Lam}$ (or a {\em coarsening}
of $w$) if the following equivalent conditions are satisfied.
\ben
\item[\rm (i)] There exists a morphism 
\footnote{$^)$The morphism $F$ is necessarily surjective.}$^)$ of 
$sl$-monoids $F : \cM_w \lra \cM_{w'}$ such that, for all $x \in B$, 
$w'(x) = F(w(x))$, and $F(\cM_{w}^\ast) = \cM_{w'}^\ast$.
\item[\rm (ii)] There exists a morphism of $sl$-monoids
$f : \wh{w}(M) \lra \wh{w'}(M)$ such that, for all $x \in B$,
$w'(x) = f(w(x))$, and $f(\wh{w}(M) \cap \Lam) = \wh{w'}(M) \cap \Lam'$.
\een
\end{de}

The relation above is a preordering inducing the equivalence
relation :

\begin{de}
 Two $l$-valuations $w, w'$ on the commutative ring $B$ are said
to be {\em equivalent} (for short, $w \sim w'$) if the 
following equivalent conditions are satisfied :
\ben
\item[\rm (i)] There is an isomorphism $F : \cM_w \lra \cM_{w'}$
of $sl$-monoids such that, for all $x \in B$, $w'(x) = F(w(x))$. 
\item[\rm (ii)] For all $x_1, \dots, x_n, y \in B$, 
$\bigwedge_{1 \leq i \leq n} w(x_i) \leq w(y) \Llra
\bigwedge_{1 \leq i \leq n} w'(x_i) \leq w'(y)$,
in particular, $A := A_w = A_{w'}$, and the $sl$-monoids
${\bf M}_w, {\bf M}_{w'}$ are isomorphic over $M = M(A, B)$.
\een
\end{de}

\begin{rems} \em 
 \ben
\item[\rm (1)] If $w'$ is coarser than $w$ then supp$(w) \sse$
supp$(w')$ and $A_w \sse A_{w'}$. If, in addition, $w$ is a
valuation then $w'$ is a valuation too, and supp$(w) =$ supp$(w')$.
\item[\rm (2)] Let $w : B \lra \wh{\Lam}$ be a $l$-valuation.
For any convex $l$-subgroup $H$ of $\cM_w^\ast$, the natural 
projection $\cM_w^\ast \lra \Gam_H := \cM_w^\ast/H$ extends to 
a morphism of $sl$-monoids $h : \cM_w \lra \wh{\Gam_H}$ 
(cf. Remark 1.7 and Proposition 1.8). The map 
$w/H = h \circ w : B \lra \wh{\Gam_H}$ is a $l$-valuation 
with $\cM_{w/H} = h(\cM_w)$, which is, up to 
equivalence, the minimal coarsening $w'$ of $w$ with 
$\cM_{w'}^\ast \cong \Gam_H$. In particular, taking $H = \{\eps\}$,
we obtain the minimal, up to equivalence, coarsening
$\wt{w} := w/\{\eps\}$ of $w$ with $\cM_{\wt{w}}^\ast \cong \cM_w^\ast$,
and hence $A_{\wt{w}} = A_w$.

If $w$ is a valuation then the coarsenings $w'$ of $w$ 
correspond uniquely, up to equivalence, to the convex
subgroups $H$ of the totally ordered Abelian group $\cM_w^\ast$
via $w' = w/H$; in particular, $\wt{w} \sim w$.
\een
\end{rems} 

The following definition provides a generalisation of the basic 
notion of {\em Manis valuation} \cite{M}.

\begin{de}
 Let $w : B \lra \wh{\Lam}$ be a $l$-valuation, $M := M(A_w, B)$ 
the $sl$-monoid associated to the ring extension $A_w \sse B$,
and $\wh{w} : M \lra \wh{\Lam}$ the canonical morphism of
$sl$-monoids. $w$ is said to be a {\em Manis $l$-valuation}
if $\wh{w}(M) \cap \Lam$ is a $l$-group generated by the
subset $\{\varphi_-\,|\,\varphi \in w(B \sm A_w)\}$. 
\end{de}

Assuming that $w$ is a Manis $l$-valuation, it follows that
$\cM_w = \wh{w}(M)$, $\cM_w^\ast = \wh{w}(M) \cap \Lam$, and
$(\cM_w^\ast)_+ = \{\bigvee_{1 \leq i \leq n} (\varphi_i)_-\,|\,
n \in \N, \varphi_i \in w(B \sm A_w)\}$.


\subsection{Pr\" ufer extensions}

$\q$ A commutative unital ring extension $A \sse B$ is said to be a 
{\em Pr\" ufer extension} if $A$ is a $B$-Pr\" ufer ring 
in the sense of Griffin \cite{Grif}.
These extensions relate to Manis valuations in much 
the same way as Pr\" ufer domains with Krull valuations. They admit 
various characterizations as shown for instance in 
\cite[I, Theorem 5.2.]{Kne} We mention only the following two 
useful criteria for a commutative ring extension $A \sse B$ 
to be a Pr\" ufer extension :
\ben
 \item[\rm (P1)] Every subring of $B$ containing $A$ 
is integrally closed in $B$.
 \item[\rm (P2)] For every element $x \in B$ there exists $y \in A$ such
that $x y \in A$ and $x (1 - x y) \in A$.
\een

In particular, a commutative ring $A$ is a {\em Pr\" ufer ring}
if and only if the ring extension $A \sse B$ is Pr\" ufer, where
$B = Q(A)$ is the total ring of quotients of $A$.

Thus the Pr\" ufer extensions are the models of an 
inductive ($\forall \exists$) theory in the first order language 
$(+, -, \cdot, 0, 1)$ of rings augmented with a unary predicate
standing for a subring $A$ of $B$. 

The Pr\" ufer extensions have a good ``multiplicative ideal
theory'' \cite[II]{Kne}. 

Given a commutative ring extension $A \sse B$, let 
$M := M(A, B), M^\ast$, and $\wt{M}$ be as defined in 1.2.
By \cite[II, Theorem 1.13.]{Kne}, the ring extension
$A \sse B$ is Pr\" ufer if and only if $I$ is $B$-invertible
for all $I \in \wt{M}$, i.e. $M^\ast = \wt{M}$.
Assuming that the ring extension
$A \sse B$ is Pr\" ufer, it follows by Remark 1.7.
(see also \cite[II, Corollary 1.14.]{Kne}) that 
$M^\ast$ is an Abelian $l$-group with $I \wedge J = I + J, 
I \vee J = I \cap J = (I^{- 1} + J^{- 1})^{- 1}$ for all 
$I, J \in M^\ast$, so $M^\ast_+ = \{I \in M^\ast\,|\,I \sse A\}$,
and $I_+ = I \cap A, I_- = (A + I)^{- 1} = 
I^{- 1} \cap A, |I| = I \cap I^{- 1} = I_+ I_- = I_+ \cap I_-$
for $I \in M^\ast$. By Remark 1.7. again, the canonical
embedding of the Abelian $l$-group $M^\ast$ into its 
commutative $l$-monoid extension $\wh{M^\ast}$, as defined 
in 1.1, extends uniquely to the morphism of $sl$-monoids 
$\wh{\mathfrak{w}} : M \lra \wh{M^\ast}$, defined by 
$\wh{\mathfrak{w}}(I)(\alpha) = I + \alpha$ for 
$I \in M, \alpha \in M^\ast_+$. Composing the morphism 
$\wh{\mathfrak{w}} : M \lra \wh{M^\ast}$ with
the map $B \lra M, x \mapsto A x$, we obtain the 
$l$-valuation $\mathfrak{w} : B \lra \wh{M^\ast}$, 
called the {\em $l$-valuation associated to the 
Pr\" ufer extension} $A \sse B$. Notice that 
$\mathfrak{w}$ is trivial (the constant map $x \mapsto A$) 
$\Llra A = B$.

We obtain $\{x \in B\,|\,\mathfrak{w}(x) \geq \gam\} = \gam$
for all $\gam \in M^\ast$, in particular, $A_{\mathfrak{w}} = A$,
and $\mathfrak{w}^{- 1}(\gam) = \gam \sm \bigcup_{\alpha \in
M^\ast, \alpha \subsetneq \gam} \alpha$ (may be empty)
for $\gam \in M^\ast$. It follows also that 
supp$(\mathfrak{w}) = \cap_{\alpha \in M^\ast_+} \alpha$ is 
the conductor of $A$ in $B$. Indeed, assuming that $x \in B \sm$
supp$(\mathfrak{w})$, let $\alpha \in M^\ast_+$ be such that $x \not \in \alpha$,
so $\beta := \alpha + Ax \in M^\ast$ and $\beta^{- 1} \subsetneq \alpha^{- 1}$.
Thus $x y \not \in A$ for all $y \in \alpha^{- 1} \sm \beta^{- 1}$,
therefore $B x \not \sse A$ as desired.  

The morphism $\wh{\mathfrak{w}} : M \lra \wh{M^\ast}$ is induced by
the $l$-valuation $\mathfrak{w} : B \lra \wh{M^\ast}$, i.e.
$\wh{\mathfrak{w}}(I) = \wedge_{1 \leq i \leq n} \mathfrak{w}(x_i)$ for 
$I = \sum_{1 \leq i \leq n} A x_i \in M$. Consequently,
the $sl$-submonoid $\wh{\mathfrak{w}}(M)$ of the $l$-monoid $\wh{M^\ast}$
contains the Abelian $l$-group $M^\ast$ and is generated 
as meet-semilattice by its submonoid $\mathfrak{w}(B)$. 
Thus $\cM_{\mathfrak{w}} = \wh{\mathfrak{w}}(M)$
is the value $sl$-monoid, while $\cM_{\mathfrak{w}}^\ast = M^\ast$ 
is the value $l$-group of the $l$-valuation $\mathfrak{w}$ 
associated to the Pr\" ufer extension $A \sse B$. As the $l$-group 
$M^\ast$ is generated by the invertible ideals 
$\mathfrak{w}(x)_- = (A + A x)^{- 1}$ for $x \in 
B \sm A$, $\mathfrak{w}$ is a Manis $l$-valuation (cf. Definition 1.17),
and hence a Pr\" ufer-Manis $l$-valuation according to the next
definition which extends to $l$-valuations the basic notion of 
{\em Pr\" ufer-Manis valuation} \cite[I, 6, Definition 1]{Kne}).

\begin{de}
Let $w : B \lra \wh{\Lam}$ be a $l$-valuation.
\ben
\item[\rm (1)] $w$ is said to be a {\em Pr\" ufer $l$-valuation}
if the ring extension $A_w \sse B$ is Pr\" ufer.
\item[\rm (2)] $w$ is said to be a {\em Pr\" ufer-Manis $l$-valuation}
if $w$ is a Pr\" ufer as well as a Manis $l$-valuation. 
\een
\end{de}

\begin{rem} \em 
Let $w : B \lra \wh{\Lam}$ be a $l$-valuation with $A_w = B$.
Let $M$ denote the $sl$-monoid of the f.g. ideals of $B$. 
Then $w$ is obviously Pr\" ufer. It is Pr\" ufer-Manis $\Llra
\cM_w^\ast = \wh{w}(M) \cap \Lam = \{\eps\}$. As an example 
of such a $l$-valuation, let $\Lam$ be a nontrivial Abelian 
$l$-group, and let $B$ be the underlying
boolean ring of the boolean algebra with support 
$\partial \wh{\Lam}_+$ and opposite order. Then the inclusion
$w : B \lra \wh{\Lam}$ is a $l$-valuation with 
$A_w = B, \cM_w = w(B) = \partial \wh{\Lam}_+, \cM_w^\ast =
\{\eps\}$. In particular, if $\Lam$ is totally ordered
then $B = \F_2$, the field with $2$ elements. 
\end{rem}

\begin{lem}
 Let $A \sse B$ be a Pr\" ufer extension, $M$-the $sl$-monoid
of f.g. $A$-submodules of $B$, $M^\ast$-the $l$-group of
invertible $A$-submodules of $B$, and 
$\mathfrak{w} : B \lra \wh{M^\ast}$-the Pr\" ufer-Manis
$l$-valuation associated as above to the Pr\" ufer extension $A \sse B$.
Let $w : B \lra \wh{\Lam}$ be a Pr\" ufer $l$-valuation 
with $A_w = A$. Then :
\ben
\item[\rm (1)] The morphism of $sl$-monoids 
$\wh{w} : M \lra \wh{\Lam}$ induces by
restriction a monomorphism of Abelian $l$-groups
$M^\ast \lra \Lam$, identifying $M^\ast$ with 
$\wh{w}(M^\ast) = \wh{w}(M)^\ast$, the $l$-group of invertible 
elements of the $sl$-monoid $\wh{w}(M)$.
\item[\rm (2)] $\wh{\mathfrak{w}} = f \circ \wh{w}$, 
$\mathfrak{w} = f \circ w$, and $\cM_{\mathfrak{w}} =
f(\wh{w}(M))$, where $f : \wh{w}(M) \lra \wh{M^\ast}$ is
the canonical morphism of $sl$-monoids extending the
isomorphism $\wh{w}(M)^\ast \cong M^\ast$ of Abelian $l$-groups.
\item[\rm (3)] $\wh{w}(M) \cap \Gam =
\wh{w}(M)^\ast$, where $\Gam$ is the convex hull of 
$\wh{w}(M)^\ast \cong M^\ast$ in $\Lam$. In particular,
$w \sim \mathfrak{w}$ whenever $\Gam = \Lam$.
\een
\end{lem}

\bp
The morphism $M^\ast \lra \Lam,\alpha \mapsto \wh{w}(\alpha)$ 
of Abelian $l$-groups is injective cf. 1.3. To show that 
$\wh{w}(M^\ast) = \wh{w}(M)^\ast$, let $\gam \in \wh{w}(M)^\ast$. 
By assumption, $\gam = \wh{w}(I), \gam^{- 1} = \wh{w}(J)$
for some $I, J \in M$, and hence $\gam = (\gam \wedge \eps)
(\gam^{- 1} \wedge \eps)^{- 1} = \wh{w}(I + A) \wh{w}(J + A)^{- 1} =
\wh{w}(K)$, where $K = (I + A) (J + A)^{- 1} \in M^\ast$, as
desired.

The assertion (2) of the lemma follows in a straightforward way.

To prove the assertion (3), let $\gam := \wh{w}(I) \in \wh{w}(M)
\cap \Gam$. We have to show that $\gam \in \wh{w}(M)^\ast =
\wh{w}(M^\ast)$. By assumption there exists $\mathbf{a} \in M^\ast$
such that $\gam \leq \wh{w}(\bf{a})$, therefore $I + \mathbf{a} \in 
M^\ast$ and $\gam = \wh{w}(I + \mathbf{a}) \in \wh{w}(M^\ast)$ as
desired. Assuming that $\Gam = \Lam$, it follows that the morphism
$f : \wh{w}(M) \lra \wh{M^\ast}$ of $sl$-monoids is injective,
and hence $\cM_w = \wh{w}(M) \cong f(\wh{w}(M)) = \cM_{\mathfrak{w}}$,
so $w \sim \mathfrak{w}$ as required.
\ep 

\begin{co} Let $A \sse B$ be a Pr\" ufer ring extension.
Then $\mathfrak{w} : B \lra \wh{M^\ast}$ is, up
to equivalence, the minimal Pr\" ufer-Manis $l$-valuation
on $B$ with $l$-valuation ring $A$, i.e. $\mathfrak{w}$
is coarser than any Pr\" ufer-Manis $l$-valuation
$w : B \lra \wh{\Lam}$ with $A_w = A$. Moreover for any
such $w$, $\cM_w^\ast \cong \cM_{\mathfrak{w}}^\ast = M^\ast$.
\end{co}

\begin{co} {\em (\cite[I, Proposition 5.1]{Kne})}
 Let $A \sse B$ be a Pr\" ufer extension such that $A \neq B$.
The necessary and sufficient condition for the Pr\" ufer-Manis
$l$-valuation $\mathfrak{w} : B \lra \wh{M^\ast}$ to be a valuation,
i.e. $M^\ast$ is totally ordered, is that the set $B \sm A$ is
multiplicatively closed.
\end{co}

\bp
Assuming that the Abelian group $M^\ast$ is totally ordered,
let $x, y \in B \sm A$ be such that $x y \in A$. As $A + Ax,
A + A y \in M^\ast$, we may assume without loss that 
$x \in A + A y$, so $x = \lam + \mu y$ for some $\lam, \mu \in A$,
and hence $P(x) = 0$, where $P(X) = X^2 - \lam X - \mu x y \in A[X]$.
Since $A$ is integrally closed in $B$, it follows that 
$x \in A$, i.e. a contradiction. Conversely, assume that
$B \sm A$ is multiplicatively closed and there exists 
$\alpha \in M^\ast$ such that $\alpha \not \sse A$ and 
$\alpha^{- 1} \not \sse A$. Choose $x \in \alpha \sm A,
y \in \alpha^{- 1} \sm A$. As $B \sm A$ is multiplicatively
closed by assumption, it follows that $x y \not \in A$, contrary
to the fact that $x y \in \alpha \alpha^{- 1} = A$.
Consequently, $M^\ast$ is totally ordered, and hence 
$\mathfrak{w}$  is a valuation, as required.
\ep

\bigskip 

In the following four sections (2 - 5) we provide an abstract axiomatic 
framework for the residue structures induced by Pr\" ufer extensions.


\section{Commutative regular semirings and quasi-semirings}

\begin{de}
By a \textbf{commutative regular semiring}, 
abbreviated a \textbf{cr-sring}, we understand an algebra 
$(R, +, \bullet,  -, ^{- 1}, \eps)$ of signature $(2, 2, 1, 1, 0)$,
satisfying the following conditions
\ben
\item[\rm (1)] $(R, +)$ and $(R, \bullet)$ are commutative regular semigroups, with
$- x$ and $x^{- 1}$ the corresponding quasi-inverses of any element $x \in R$;
\item[\rm (2)] $E^{+} \cap E^{\bullet} = \{\eps\}$, where
$E^{+} := \{x \in R\,\mid\,2x := x + x = x\} = \{e^{+}(x) := x + (- x)\,\mid\,
x \in R\}, E^{\bullet} := \{x \in R\,\mid\,x^{2} := x \bullet x = x\} =
\{e^{\bullet}(x) := x \bullet x^{- 1}\,\mid\,x \in R\}$;
\item[\rm (3)] \textbf{Distributive law :} $x \bullet y + x \bullet z =
x \bullet (y + z)$ for all $x, y, z \in R$.
\een
\end{de}

\begin{rems} \em
(1) The \textit{commutative}, not necessarily 
unital, \textit{regular rings}
(in the sense of von Neumann) are exactly those cr-srings $R$
for which $(R, +)$ is an Abelian group. In this case, the 
neutral element $\eps = 0$ is the unique element of $E^{+}$, 
while $E^{\bullet}$ is a quasiboolean lattice with the least element $0$
under the operations $x \wedge y = x \bullet y, x \vee y = x + y - x \bullet y,
x \setminus y = x - x \bullet y$. $E^{\bullet}$ is a Boolean algebra $\Llra R$
has a unit $1$, in which case $\urcorner x = 1 - x$. Thus the fields and the
Boolean rings are amongst the simplest examples of cr-srings.

(2) On the opposite side, the 
\textit{Abelian $l$-groups} are identified with
those cr-srings $R$ for which $(R, \bullet)$ is an Abelian group.
In this case, the neutral element $\eps = 1$ is the unique element
of $E^{\bullet}$, while $E^{+} = R$, so $(R, +)$ is a semilattice. Indeed,
we obtain by 2.1.(3) and 2.1.(2) the identities
$$x + x = x \bullet \eps + x \bullet \eps = x \bullet (\eps + \eps) =
x \bullet \eps = x$$
for all $x \in R$. As the group and the semilattice
operations $\bullet$ and $+$ are compatible by 2.1.(3), it follows that 
$(R, \bullet, \wedge, \vee)$ is an Abelian $l$-group, where
$$x \wedge y = x + y, x \vee y  = (x^{- 1} + y^{- 1})^{- 1}$$

Notice that the trivial cr-sring (the singleton) is the only common 
member of the two remarkable subclasses of cr-srings considered above.
\item[\rm (3)] In the presence of the axioms 2.1.(1) and 2.1.(3), 
the axiom 2.1.(2) is equivalent with the conjunction of the 
following two equational axioms
\ben
\item[\rm 2.1.(2')] $\eps \in E^{+}$, i.e. $2 \eps = \eps$, and
\item[\rm 2.1.(2'')] $E^{+} \sse R^{\bullet}_{\eps} := \{x \in R\,\mid\,
e^{\bullet}(x) = \eps\}$, i.e. $e^{\bullet}(e^{+}(x)) = \eps$ for all $x \in R$.
\een
The implications 2.1.(2) $\Lra$ 2.1.(2'), and (2.1.(2') $\bigwedge$ 
2.1(2'')) $\Lra$ 2.1.(2) are obvious, so it remains to prove the
implication 2.1.(2) $\Lra$ 2.1.(2''). Assuming that $x \in E^{+}$, 
we obtain by 2.1.(3)  
$$e^{\bullet}(x) = x^{- 1} \bullet
x = x^{- 1}  \bullet (x + x) = x^{- 1} \bullet x + x^{- 1} \bullet x =
2 e^{\bullet}(x),$$ 
so $e^{\bullet}(x) \in E^{+} \cap E^{\bullet}$, and
hence $e^{\bullet}(x) = \eps$ by 2.1.(2), as desired. 

Thus the class of all cr-srings is a variety of algebras of signature 
(2, 2, 1, 1, 0), i.e. it is closed under homomorphic images, 
subalgebras and direct products. 
\end{rems}

We introduce a larger class of algebras of the same signature as above
by relaxing suitably the distributive law 2.1.(3) and adding three 
new natural axioms as follows. 

\begin{de}
By a \textbf{commutative regular quasi-semiring}, 
abbreviated a \textbf{cr-qsring}, we understand an algebra 
$(R, +, \bullet,  -, ^{- 1}, \eps)$ of signature $(2, 2, 1, 1, 0)$,
satisfying the axioms $2.1.(1),(2)$, and the following new axioms
\ben
\item[\rm (3')]  \textbf{Quasidistributive law :} $x \bullet y + x \bullet z \leq
x \bullet (y + z)$ for all $x, y, z \in R$, where $\leq$ denotes the 
canonical partial order on $(R, +) : x \leq y \Longleftrightarrow x = y + e^{+}(x)$;
\item[\rm (4)] $y \leq z \Lra x \bullet y \leq x \bullet z$ for all $x, y, z \in R$;
\item[\rm (5)] $- (x \bullet y) = x \bullet (- y)$ for all $x, y \in R$;
\item[\rm (6)] $e^{+}(x + y) \leq e^{+}(x \bullet y)$ for all $x, y \in E^{\bullet}$.
\een
\end{de}

The next lemmas collect some basic properties of the cr-qsrings.

\begin{lem} Let $R$ be a cr-qsring. With the notation above, the following
statements hold.
\ben
\item[\rm (1)] The unary operations $x \mapsto - x$ and $x \mapsto x^{- 1}$
commute, i.e. $(- x)^{- 1} = - (x^{- 1})$ for all $x \in R$.
\item[\rm (2)] The multiplication induces an action of the
multiplicative semigroup $(R, \bullet)$ on the additive
idempotent semigroup (semilattice) $(E^{+}, +)$, i.e.
$x \bullet f \in E^{+}$, and
$x \bullet (f + g) = x \bullet f + x \bullet g$ for all
$x \in R, f, g \in E^{+}$. 
\item[\rm (3)] $E^{+} = R^{\bullet}_{\eps} := \{x \in R\,\mid\,
e^{\bullet}(x) = \eps\}$.
\item[\rm (4)] $(E^{+}, \bullet, \wedge, \vee)$ is a 
multiplicative Abelian $l$-group, with the lattice operations 
$e \wedge f = e + f, e \vee f = (e^{- 1} + f^{- 1})^{- 1}$.
\item[\rm (5)] The map $v : R \lra E^{+}, x \mapsto v(x) :=
\eps \bullet x$ is a {\em quasivaluation} (abreviatted a {\em qvaluation}) 
on the cr-qsring $R$ with values
in the multiplicative Abelian $l$-group $E^{+}$, i.e. 
 $v : (R, \bullet) \lra (E^{+}, \bullet)$ is a surjective morphism of
semigroups, and $v(x + y) \geq v(x) \wedge v(y) = v(x) + v(y)$
for all $x, y \in R$. In addition, $v(- x) = v(x)$ for all $x \in R$, and
$v(x) \leq v(y)$ provided $x \leq y$.
\item[\rm (6)] The epimorphism $v : (R, \bullet) \lra (E^{+}, \bullet)$ 
splits, and the embedding $E^{+} \hookrightarrow R$ 
is the unique section $s$ of $v$ satisfying in addition 
$s(f + g) = s(f) + s(g)$ for all $f, g \in E^{+}$.
Moreover the qvaluation $v$ satisfies the following universal property:
for every Abelian $l$-group $(\Lam, +, \wedge, \vee)$ and for every
morphism of semigroups \\
$w : (R, \bullet) \lra (\Lam, +)$ satisfying 
$w(f + g) = w(f) \wedge w(g)$ for all $f, g \in E^{+}$,
there exists a unique morphism of $l$-groups $\varphi : E^{+} \lra \Lam$
such that $\varphi \circ v = w$. In particular, any such map
$w$ is a qvaluation provided it is surjective.
\een
\end{lem}

\bp
(1) follows easily from 2.3.(5).

(2) For $x \in R, f \in E^{+}$, we obtain by
2.3.(3')
$$x \bullet f = x \bullet (f + f) \geq x \bullet f + x \bullet f,$$
so
$$x \bullet f + x \bullet f = x \bullet f + e^{+}(x \bullet f
+x \bullet f) = x \bullet f + e^{+}(x \bullet f) = x \bullet f,$$
and hence $x \bullet f \in E^{+}$, as desired.

For $x \in R, f, g \in E^{+}$, it follows by 2.3.(3') again that
$x \bullet (f + g) \geq x \bullet f + x \bullet g$. On the
other hand, $f, g \in E^{+} \Lra f + g = f \wedge g \leq f$,
and hence $x \bullet (f + g) \leq x \bullet f$ by 2.3.(4).
Similarly, we obtain $x \bullet (f + g) \leq x \bullet g$,
therefore $x \bullet (f + g) \leq x \bullet f \wedge x \bullet
g = x \bullet f + x \bullet g$, since $x \bullet f, x \bullet g
\in E^{+}$ as shown above. Consequently, we obtain the required
equality.

Thus we have obtained a morphism of semigroups $\psi : (R, \bullet) \lra
{\rm End}(E^{+}, +)$ defined by $\psi(x)(f) = x \bullet f$ for $x \in R,
f \in E^{+}$. We shall see later that the image of $\psi$ is contained
in Aut$(E^{+}, +)$.

(3, 4) We already know from (2) above that $E^{+}$ is closed under
multiplication. It remains to show that $E^{+} = R^{\bullet}_{\eps}$ to
conclude that $E^{+}$ is an Abelian group under multiplication, with
$\eps$ as neutral element. Assuming that $f \in E^{+}$, it follows by (2) that
$e^{\bullet}(f) = f^{- 1} \bullet f \in E^{+} \cap E^{\bullet}$, and hence
$e^{\bullet}(f) = \eps$, i.e. $f \in R^{\bullet}_{\eps}$, since 
$E^{+} \cap E^{\bullet} = \{ \eps \}$ by 2.1.(2). Conversely, 
assuming that $x \in R^{\bullet}_{\eps}$, we obtain 
$x = x \bullet e^{\bullet}(x) = x \bullet \eps \in E^{+}$ by (2) 
since $\eps \in E^{+}$. Thus $E^{+} = R^{\bullet}_{\eps}$,
as desired.

Since the group operation $\bullet$ and the semilattice operation $+$ on $E^{+}$
are compatible by (2), we obtain the required structure of Abelian $l$-group
on $E^{+}$. In particular, the multiplicative Abelian group $(E^{+}, \bullet)$
is torsion free.

(5, 6) Setting $v(x) := \eps \bullet x \in E^{+}$ for $x \in R$, 
we obtain for arbitrary $x, y \in R$ : 
$$v(x) \bullet v(y) = (\eps \bullet x) \bullet (\eps \bullet y) =
\eps^{2} \bullet (x \bullet y) = \eps \bullet (x \bullet y) = v(x \bullet y),$$
and
$$v(x + y) = \eps \bullet (x + y) \geq \eps \bullet x + \eps \bullet y =
v(x) + v(y)$$ 
by 2.3.(3').

Notice that $v(- x) = v(x)$ for all $x \in R$ since 
the Abelian group $(E^{+}, \bullet)$ is torsion free, and
$(- x)^{2} = x^{2}$ (by 2.3.(4)) $\Lra v(- x)^{2} = v(x)^{2}$,
while $x \leq y \Lra v(x) = \eps \bullet x \leq \eps \bullet y =
v(y)$ by 2.3.(4).

For $x \in E^{+}$, we obtain $v(x) = \eps \bullet x = x$ by (3).
Thus $v$ is a qvaluation, and the embedding 
$i : (E^{+}, \bullet) \hookrightarrow (R, \bullet)$ is a section 
of the epimorphism $v :(R, \bullet) \lra (E^{+}, \bullet)$ satisfying 
the supplementary property $i(f + g) = i(f) + i(g)$ for all 
$f, g \in E^{+}$. To prove its uniqueness, let 
$s : E^{+} \lra R$ be a map satisfying 
$v \circ s = 1_{E^{+}}$ and
$s(f + g) = s(f) + s(g)$ for $f, g \in E^{+}$. We have to show that
$s(f) = f$ for all $f \in E^{+}$. It follows that $s(f) \in E^{+}$ for all
$f \in E^{+}$ since $s(f) + s(f) = s(f + f) = s(f)$, therefore
$s(f) = v(s(f)) = f$, as required.

To prove the universal property of the qvaluation $v$,
let $\Lam$ and $w : R \lra \Lam$ be as in the statement (6) above.
First notice that if a map $\varphi : E^{+} \lra \Lam$ satisfies 
$\varphi \circ v = w$, then $\varphi(f) = \varphi(v(f)) = w(f)$ for
all $f \in E^{+}$, so $\varphi =  w\mid_{E^{+}}$ is unique. By
assumption, the map $\varphi$ as defined above is a morphism
of $l$-groups, and $\varphi(v(x)) = w(v(x)) = w(\eps \bullet x) =
w(\eps) + w(x) = 0 + w(x) = w(x)$ for all $x \in R$, so 
$\varphi \circ v = w$ as desired. As a composition of a qvaluation
with a morphism of $l$-groups, the map $w$ is a qvaluation too
provided it is surjective.

Finally notice that the morphism $\psi : (R, \bullet) \lra {\rm End}(E^{+}, +)$
as defined in (2) is the composition of the epimorphism 
$v : (R, \bullet) \lra (E^{+}, \bullet)$ with the canonical
monomorphism $\iota : (E^{+}, \bullet) \lra {\rm Aut}(E^{+}, +) \sse
{\rm End}(E^{+}, +)$ defined by $\iota(f)(g) = f \bullet g$ for
$f, g \in E^{+}$. Indeed, we obtain 
$$\psi(x)(f) = x \bullet f =
x \bullet (\eps \bullet f) = (x \bullet \eps) \bullet f =
v(x) \bullet f = \iota(v(x))(f)$$
for $x \in R, f \in E^{+}$, i.e. $\psi = \iota \circ v$ as desired. 
In particular, the image of $\psi$ is a subgroup of Aut$(E^{+}, +)$,
identified with the multiplicative Abelian group $(E^{+}, \bullet)$.
\ep 

\begin{rem} \em
Though very similar, the notions of valuation and qvaluation are
quite different. Thus, the constant map $x \mapsto 0$, the unique
qvaluation on any field $F$, is distinct from the trivial valuation
on $F$. However the qvaluations are naturally related to the valued 
fields and the Pr\" ufer domains \cite{Res}.
\end{rem}

\begin{lem} The surjective maps $e^{+} : R \lra E^{+}, e^{\bullet} : R \lra
E^{\bullet}$, and $v : R \lra E^{+}$ have the following properties.
\ben
\item[\rm (1)] $v \leq e^{+}$, i.e. $v(x) \leq e^{+}(x)$ for all $x \in R$.
In particular, the restriction $e^{+}\mid_{E^{\bullet}}$ takes values
in $E^{+}_{+} := \{f \in E^{+}\,\mid\,f \geq \eps\}$, the positive cone
of the Abelian $l$-group $(E^{+}, \bullet, + = \wedge, \vee)$.
\item[\rm (2)] $e^{+} = v \bullet (e^{+} \circ e^{\bullet})$, i.e.
$e^{+}(x) = v(x) \bullet e^{+}(e^{\bullet}(x)) = x \bullet e^{+}(e^{\bullet}(x))$
for all $x \in R$. In particular, $e^{+}(x^{- 1}) = v(x)^{- 2} \bullet e^{+}(x) =
x^{- 2} \bullet e^{+}(x)$ for all $x \in R$.
\item[\rm (3)] The restriction $e^{+}\mid_{E^{\bullet}} : (E^{\bullet}, \bullet)
\lra (E^{+}_{+}, +)$ is a morphism of semilattices, i.e. $e^{+}(x \bullet y) =
e^{+}(x) + e^{+}(y)$ for all $x, y \in E^{\bullet}$, which preserves 
the common least element $\eps$.
\item[\rm (4)] \textbf{Leibniz's rule :} The additive map $e^{+} : R \lra E^{+}$ 
is a ``derivation'', i.e. \\ $e^{+}(x \bullet y) = x \bullet e^{+}(y) +
y \bullet e^{+}(x)$ for all $x, y \in R$.
\item[\rm (5)] $\{x \in E^{\bullet}\,\mid\,e^{+}(x) = \eps\} = 
\{x + \eps\,\mid\,x \in E^{\bullet}\}$.
\item[\rm (6)] $e^{+}(e^{\bullet}(x)) \leq e^{+}(e^{\bullet}(y))$
provided $x \leq y$.
\een
\end{lem}

\bp
(1) For any $x \in R$, it follows by Lemma 2.4.(5, 6) that
$$e^{+}(x) = v(e^{+}(x)) = v(x + (- x)) \geq v(x) + v(-x) = v(x) + v(x) = v(x),$$
while for $x \in E^{\bullet}$, we obtain $e^{+}(x) \geq v(x) = v(x \bullet
x^{- 1}) = \eps$, as required. 

(2) For any $x \in R$, it follows by 2.3.(3') and 2.3.(5) that
$$v(x)^{- 1} \bullet e^{+}(x) = x^{- 1} \bullet (x + (- x)) \geq x^{- 1} \bullet
x +(- (x^{- 1} \bullet x)) = e^{+}(e^{\bullet}(x))$$
and
$$v(x) \bullet e^{+}(e^{\bullet}(x)) = x \bullet (e^{\bullet}(x) +
(- e^{\bullet}(x))) \geq x + (- x) = e^{+}(x),$$
and hence the desired identity.

(3) Let $x, y \in E^{\bullet}$, in particular, $v(x) = v(y) = \eps$.
It follows by 2.3.(3') and 2.3.(5) again that 
$$e^{+}(x \bullet y) = x \bullet y + x \bullet (- y) \leq x \bullet e^{+}(y) =
v(x) \bullet e^{+}(y) = e^{+}(y),$$
and similarly, $e^{+}(x \bullet y) \leq e^{+}(x)$. As $(E^{+}, +)$ is 
a meet-semilattice, we deduce the inequality 
$e^{+}(x \bullet y) \leq e^{+}(x) + e^{+}(y)$. To get an equality, we have to
appeal for the first time to the axiom 2.3.(6).

(4) Let $x, y \in R$. By (2) and (3) above, we obtain
$e^{+}(x \bullet y) = x \bullet y \bullet e^{+}(e^{\bullet}(x \bullet y)) =
x \bullet y \bullet e^{+}(e^{\bullet}(x) \bullet e^{\bullet}(y)) = x \bullet y
\bullet (e^{+}(e^{\bullet}(x)) + e^{+}(e^{\bullet}(y))) = y \bullet e^{+}(x) +
x \bullet e^{+}(y)$, as required.

(5) Assuming that $x \in E^{\bullet}$ and $e^{+}(x) = \eps$, it follows
that $x = x + \eps$ since $x + \eps \leq x$  and $e^{+}(x + \eps) = e^{+}(x) +  \eps =
\eps = e^{+}(x)$. Conversely, for any $x \in E^{\bullet}$, we obtain 
$e^{+}(x + \eps) = e^{+}(x) + \eps = \eps$ since 
$x \in E^{\bullet} \Lra e^{+}(x) \geq \eps$, so it remains to show that 
$x + \eps \in E^{\bullet}$. As $v(x) = \eps \bullet x = \eps$, it
follows by 2.3.(4) that $(x + \eps)^{2} \geq x + \eps$, and hence,
again by 2.3.(4), $e^{\bullet}(x + \eps) = (x + \eps) \bullet (x + \eps)^{- 1} \leq 
(x + \eps)^{2} \bullet (x + \eps)^{- 1} = x + \eps$, therefore
$e^{\bullet}(x + \eps) = x + (\eps + e^{+}(e^{\bullet}(x + \eps))) = x + \eps$,
i.e. $x + \eps \in E^{\bullet}$, as desired.
\item[\rm (6)] Let $x, y \in R$ be such that $x \leq y$, in particular, 
$e^{+}(x) \vee v(y) \leq e^{+}(y)$ by (1). By (2) it follows that 
$e^{+}(e^{\bullet}(x)) = e^{+}(x) \bullet v(x)^{- 1} = e^{+}(x) \bullet 
v(y + e^{+}(x))^{- 1} \leq e^{+}(x) \bullet (v(y) + e^{+}(x))^{- 1} = 
(e^{+}(x) \vee v(y)) \bullet v(y)^{- 1} \leq e^{+}(y) \bullet v(y)^{- 1} =
e^{+}(e^{\bullet}(y))$, as required.
\ep

A characterization of those cr-qsrings which are cr-srings is given by 
the next lemma.

\begin{lem}
The following statements are equivalent for an algebra 
$(R, +, \bullet, -, ^{- 1}, \eps)$ of signature $(2, 2, 1, 1, 0)$.
\ben
\item[\rm (1)] $R$ is cr-sring.
\item[\rm (2)] $R$ is a cr-qsring satisfying the following
equivalent conditions.
\ben
\item[\rm (i)] $e^{+} = v$;
\item[\rm (ii)] $e^{+}(x) = \eps$ for all $x \in E^{\bullet}$;
\item[\rm (iii)] $e^{+}(x \bullet y) = e^{+}(x) \bullet e^{+}(y)$ for
all $x, y \in R$;
\item[\rm (iv)] $v(x + y) = v(x) + v(y)$ for all $x, y \in R$.
\een
\een 
\end{lem}

\bp
First notice that in a cr-qsring $R$, the conditions (i) - (iv)
above are equivalent since (i) $\Llra$ (ii) by Lemma 2.6.(2),
(i) $\Lra$ (iii) and (i) $\Lra$ (iv) are obvious, 
(iii) $\Lra$ (ii) as $x \in E^{\bullet} \Lra e^{+}(x) =
e^{+}(x^{2}) = e^{+}(x)^{2}$ (by assumption) $\Lra e^{+}(x) = \eps$
for all $x \in E^{\bullet}$, while (iv) $\Lra$ (i) as
$e^{+}(x) = v(e^{+}(x)) = v(x - x) = v(x) + v(- x)$ (by assumption),
so $e^{+}(x) = 2 v(x) = v(x)$ for all $x \in R$, as desired.

(2) $\Lra$ (1) : We have to show that $R$ satisfies 2.1.(3).
As $x \bullet y + x \bullet z \leq x \bullet (y + z)$ by 2.3.(3'),
the desired equality is a consequence of the identity
$$e^{+}(x \bullet y + x \bullet z) = e^{+}(x \bullet (y + z)),$$
which is immediate by (iii) and the identity $f \bullet g + f \bullet h =
f \bullet (g + h)$ for $f, g, h \in E^{+}$.

(1) $\Lra$ (2) : Assuming that $R$ is a cr-sring, we have to show
that $R$ satisfies the axioms 2.3.(4), (5), (6), and the equivalent
conditions (i) - (iv) above. First notice that for all $x \in R, f \in E^{+}$,
it follows by 2.1.(3) that $x \bullet f = x \bullet (f + f) =
x \bullet f + x \bullet f$, i.e. $x \bullet f \in E^{+}$. Consequently,
assuming $y \leq z$, i.e. $y = z + f$ for some $f \in E^{+}$, it
follows again by 2.1.(3) that $x \bullet y = x \bullet z + x \bullet f$,
so $x \bullet y \leq x \bullet z$ as $x \bullet f \in E^{+}$. Thus 2.3.(4)
is true on $R$. On the other hand, 2.3.(5) is an immediate consequence of
2.1.(3), while 2.3.(6) is obvious since $R$ satisfies the condition (ii) :
assuming that $x \in E^{\bullet}$, it follows
by 2.1.(3) that $e^{+}(x)^{2} = (x - x) \bullet (x - x) = 2 (x^{2} - x^{2}) 
= x - x = e^{+}(x)$, therefore $e^{+}(x) \in E^{+} \cap E^{\bullet} =
\{\eps\}$ by 2.1.(2).
\ep

\begin{rems} \em
(1) Since all the axioms involving the binary 
relation $\leq$ are equivalent with identities, while, 
as in the case of cr-srings (see Remarks 2.2.(3)),
the axiom 2.1.(2) can be replaced by a conjunction of identities,
the class of all cr-qsrings is a variety of algebras of signature 
(2, 2, 1, 1, 0), and the cr-srings form a subvariety of it.

(2) For any cr-qsring $R$, the subset 
$$\wt{R} := \{x \in R\,\mid\,e^{+}(x) = v(x)\} =
\{x \in R\,\mid\,e^{+}(e^{\bullet}(x)) = \eps\}$$  
is the support of the largest cr-sring contained in $R$.
By Lemma 2.7., it suffices to show that $\wt{R}$ is a substructure
of $R$. The closure under the unary operations $-$ and $^{- 1}$ is obvious,
the closure under the multiplication follows by Lemma 2.6.(3), while
for $x, y \in \wt{R}$, it follows by Lemma 2.6.(1) and Lemma 2.5.(5) that
$$v(x + y) \leq e^{+}(x + y) = e^{+}(x) + e^{+}(y) =
v(x) + v(y) \leq v(x + y),$$
therefore $x + y \in \wt{R}$. Notice
that $E^{+} \sse \wt{R}$, and $\wt{R}^{+}_{\eps} := \{x \in \wt{R}\,\mid\,
e^{+}(x) = \eps\}$ is the largest regular ring contained in $R$. 
\end{rems}
\bigskip


\section{Directed commutative regular quasi-semirings}

$\q$ We introduce an interesting class of cr-qsrings whose 
intersection with the class of all cr-srings is the
class of all Abelian $l$-groups.

\begin{pr} The following properties are equivalent for
a cr-qsring $R$.
\ben
\item[\rm (1)] $\eps \leq x$ for all $x \in E^{\bullet}$.
\item[\rm (2)] $v \leq 1_{R}$, i.e. $v(x) \leq x$ for all $x \in R$.
\item[\rm (3)] $E^{+} \cap {x\downarrow} = {v(x)\downarrow}$ for
all $x \in R$, where ${x\downarrow} := \{y \in R\,\mid\,y \leq x\}$.
\item[\rm (4)] $x + v(x - y) = y + v(x - y)$ for all $x, y \in R$.
\item[\rm (5)] There exists the meet $x \wedge y$ with respect to the
partial order $\leq$ for any pair $(x, y)$ of elements of $R$, and
the semilattice operation $\wedge$ is compatible with the
multiplication, i.e. $x \bullet (y \wedge z) = (x \bullet y) \wedge
(x \bullet z)$ for all $x, y, z \in R$.
\item[\rm (6)] The partial orders $\leq$ and $\preceq$ coincide on
$E^{\bullet}$, where $x \preceq y \Llra x = y \bullet e^{\bullet}(x)$
is the partial order induced by multiplication.
\een
\end{pr}

\bp
(1) $\Lra$ (2) : Let $x \in R$. By assumption, $\eps \leq e^{\bullet}(x)$,
and hence $v(x) = x \bullet \eps \leq x \bullet e^{\bullet}(x) = x$ by
2.3.(4).

(2) $\Lra$ (3) : Let $x \in R, f \in E^{+} \cap x\downarrow$. By
Lemma 2.4.(5, 6), it follows that $f = v(f) \leq v(x)$. Conversely, 
let $y \in v(x)\downarrow$, i.e. $y \leq v(x)$, in particular, 
$y = v(x) + e^{+}(y) \in E^{+}$ since $v(x) \in E^{+}$. 
As $v(x) \leq x$ by assumption, we deduce also that
$y \in x\downarrow$ as desired.

(3) $\Lra$ (1) is obvious since $v(x) = \eps$ for all $x \in E^{\bullet}$.

(2) $\Lra$ (4) : Let $x, y \in R$. By assumption $v(x - y) \leq x - y$,
so $v(x - y) = v(y - x) = x - y + v(x - y)$. Consequently,
$y + v(x - y) = y + x - y + v(x - y) = x + (e^{+}(y) + v(x - y)) =
x + v(x - y)$ since $v(x - y) \leq e^{+}(x - y) = e^{+}(x) + e^{+}(y) \leq
e^{+}(y)$.

(4) $\Lra$ (5) : Let $x, y \in R$. By assumption, $z := x + v(x - y) =
y + v(x - y)$, so $z \leq x$ and $z \leq y$. It remains to show that 
$u \leq x, u \leq y \Lra u \leq z$ to conclude that $z = x \wedge y$
is the meet of the pair $(x, y)$ with respect to the partial order $\leq$.
For such an element $u$, we obtain $- u \leq - y$, so $e^{+}(u) = u - u \leq
x - y$, and hence $e^{+}(u) = v(e^{+}(u)) \leq v(x - y)$ by Lemma 2.4.(5, 6). 
Consequently, $u = u + e^{+}(u) \leq x + v(x - y) = z$,
as required.

It remains to verify the compatibility with the multiplication of the
semilattice operation $\wedge$. Let $x, y, z \in R$. The inequality
$x \bullet (y \wedge z) \leq x \bullet y \wedge x \bullet z$ follows by
2.3.(4). On the other hand, $x \bullet y - x \bullet z \leq x \bullet (y - z)$
(by 2.3.(4, 5)) implies by Lemma 2.4.(5) and 2.3.(4) the opposite inequality 
$$x \bullet y \wedge x \bullet z = x \bullet y + v(x \bullet y - x \bullet z) \leq
x \bullet y + v(x \bullet (y - z)) = $$ 
$$x \bullet y + x \bullet v(y - z) \leq x \bullet (y + v(y - z)) = 
x \bullet (y \wedge z).$$

(5) $\Lra$ (6) : Since for all $x, y \in E^{\bullet},\,x \bullet y \in E^{\bullet}$ 
is the meet of the pair $(x, y)$ with respect to the partial order 
$\preceq$, we have to show that $x \bullet y = x \wedge y$, the meet of 
the same pair with respect to the partial order $\leq$, which exists 
by assumption.

Let us show that $x \wedge y \in E^{\bullet}$ provided 
$x, y \in E^{\bullet}$. Since $x, y \in E^{\bullet}$ and the 
operations $\wedge$ and $\bullet$ are compatible by assumption, 
it follows that $(x \wedge y)^{2} = x \wedge y \wedge (x \bullet y)$,
in particular, $(x \wedge y)^{2} \leq x \wedge y$. Moreover the
last inequality becomes an equality since, using Lemma 2.6.(4) and 
again the compatibility above, we obtain
$$e^{+}((x \wedge y)^{2}) = (x \wedge y) \bullet e^{+}(x \wedge y) =
(x \bullet e^{+}(x \wedge y)) \wedge (y \bullet e^{+}(x \wedge y)) =
e^{+}(x \wedge y),$$ 
as $x \bullet f = f$ for $x \in E^{\bullet}, f \in E^{+}$.

Consequently, $x \wedge y \in E^{\bullet}$, and $x \wedge y = x \wedge y \wedge
(x \bullet y)$, i.e. $x \wedge y \leq x \bullet y$, so it remains to show
that $e^{+}(x \bullet y) \leq e^{+}(x \wedge y)$ to obtain the required
equality $x \wedge y = x \bullet y$. As 
$e^{+}\mid_{E^{\bullet}} : (E^{\bullet}, \bullet) \lra (E^{+}, +)$ is a 
morphism of semilattices by Lemma 2.6.(3), it suffices to check
that $x \bullet y \preceq x \wedge y$. Using again the compatibility
of the operations $\wedge$ and $\bullet$, we obtain 
$(x \bullet y) \bullet (x \wedge y) = x^{2} \bullet y \wedge x \bullet y^{2} =
x \bullet y$, i.e. $x \bullet y \preceq x \wedge y$, as desired.

(6) $\Lra$ (1) : Let $x \in E^{\bullet}$. Since $v(x) = \eps \bullet x =
\eps$, i.e. $\eps \preceq x$, it follows by assumption that $\eps \leq x$,
as required.
\ep

\begin{de} A cr-qsring $R$ is said to be {\em directed} if $R$ satisfies
the equivalent conditions $(1) -(6)$ from {\em Proposition 3.1.} 
\end{de}

\begin{rems} \em
(1) A directed cr-qsring $R$ is {\em trivial}, i.e. $R = \{\eps\}$,
if and only if $\Ep = \{\eps\}$.

(2) By Lemma 2.7., Proposition 3.1. and Remarks 2.2.(2),
it follows that the Abelian $l$-groups are identified with those
cr-srings which are also directed cr-qsrings.

(3) By Remarks 2.8.(1) and Proposition 3.1.(4), the class 
of all directed cr-qsrings is a subvariety of the variety of
cr-qsrings.
\end{rems}

The next lemmas collect some basic properties of the directed 
cr-qsrings.

\begin{lem} Let $R$ be a directed cr-qsring. Then the following assertions hold.
\ben
\item[\rm (1)] $e^{+}(x \wedge y) = v(x - y) = \eps \bullet (x - y)$,
and $v(x \wedge y) = v(x) \wedge v(y) = v(x) + v(y)$ for all $x, y \in R$.
In particular, $v(x) = x \wedge e^{+}(x)$ for all $x \in R$.
\item[\rm (2)] $v(x - y) = (v(x) + v(y)) \bullet e^{+}(e^{\bullet}(x \wedge y))$
for all $x, y \in R$. 

In particular,
$v(x \pm y ) \leq v(x) \bullet e^{+}(e^{\bullet}(y)) =
v(x) \bullet v(y)^{- 1} \bullet e^{+}(y)$ for all $x, y \in R$.
\item[\rm (3)] The semilattice operation $\wedge$ is compatible with
the operations $+$ and $-$, i.e. \\ 
$x + (y \wedge z) = (x + y) \wedge (x + z)$
and $(- y) \wedge (- z) = - (y \wedge z)$ for all $x, y, z \in R$.
\item[\rm (4)] The following assertions are equivalent for $x, y \in R$.
\ben
\item[\rm (i)] $e^{+}(x \wedge y) = e^{+}(x) + e^{+}(y)$;
\item[\rm (ii)] $e^{+}(x \wedge y) = x - y$;
\item[\rm (iii)] $x - y \in E^{+}$.
\een 
In particular, $e^{+}(x \wedge y) = e^{+}(x) + e^{+}(y) = x - y = y - x$
provided the pair $(x, y)$ is bounded above with respect to the 
partial order $\leq$, i.e. $x \leq z$ and $y \leq z$ for some $z \in R$.
Moreover $x \wedge y = x + y - z$, $x + y = x \wedge y + z$,
and $(x + y) \wedge z \leq x \wedge y$ provided $x \leq z$ and $y \leq z$;
in particular, $x \wedge y = x + y \Llra x + y \leq z$.
\item[\rm (5)] $e^{+}(x \wedge y \wedge z) = e^{+}(x \wedge z) +
e^{+}(y \wedge z)$ for all $x, y, z \in R$.
\een
\end{lem}

\bp
(1) By Proposition 3.1.(4, 5), $x \wedge y = x + v(x - y)$,
and hence $e^{+}(x \wedge y) = e^{+}(x) + v(x - y) = v(x - y)$
since $v(x - y) \leq e^{+}(x - y) = e^{+}(x) + e^{+}(y) \leq e^{+}(x)$.
To show that $v(x \wedge y) = v(x) + v(y)$, notice that 
$$(x \wedge y \leq x) \bigwedge (x \wedge y \leq y) \Lra 
v(x \wedge y) \leq v(x) \wedge v(y) = v(x) + v(y)$$ 
by Lemma 2.4.(5), while $v(x) \leq x$ and $v(y) \leq y$ (by
Proposition 3.1.(2)) imply $v(x) + v(y) = v(x) \wedge v(y) \leq
x \wedge y$, therefore  $v(x) + v(y) = v(v(x) + v(y)) \leq 
v(x \wedge y)$, again by Lemma 2.4.(5).

The identity $v(x) = x \wedge e^{+}(x)$ is a consequence of
the inequality $v(x) \leq x \wedge e^{+}(x)$ and of the
equality $e^{+}(x \wedge e^{+}(x)) = v(x - e^{+}(x)) = v(x) = e^{+}(v(x))$.

(2) The identity is immediate by (1) and Lemma 2.6.(2),
while the inequality is a consequence of the identity since
$v(x) + v(y) \leq v(x)$ and 
$$x \wedge y \leq y \Lra e^{+}(e^{\bullet}(x \wedge y)) \leq 
e^{+}(e^{\bullet}(y))$$
by Lemma 2.6.(6).

(3) The inequality $x + (y \wedge z) \leq (x + y) \wedge
(x + z)$ is obvious. To get an identity, we note that 
$$e^{+}((x + y) \wedge (x + z)) = v((x + y) - (x + z)) = 
v(y - z + e^{+}(x)) \leq e^{+}(x) + v(y - z) = e^{+}(x + (y \wedge z))$$
since 
$v(y - z + e^{+}(x)) \leq e^{+}(y - z + e^{+}(x)) \leq e^{+}(x)$ and
$$y - z + e^{+}(x) \leq y - z \Lra v(y - z + e^{+}(x)) \leq v(y - z).$$

The compatibility of $\wedge$ with $-$ is obvious.

(4) (ii) $\Lra$ (i) is obvious.

(iii) $\Lra$ (ii) : We get $e^{+}(x \wedge y) = v(x - y) = x - y$
since $x - y \in E^{+}$ by assumption, and $v$ is the identity on
$E^{+}$.

(i) $\Lra$ (iii) : As $x \wedge y \leq x, y$, it follows by assumption
that 
$$x \wedge y = x + e^{+}(x \wedge y) = x + e^{+}(x) + e^{+}(y) =
x + e^{+}(y),$$ 
and similarly, $x \wedge y = y + e^{+}(x)$. Consequently,
$$x - y = (x + e^{+}(y)) - (y + e^{+}(x)) = e^{+}(x \wedge y),$$ 
so $x - y \in E^{+}$ as desired.

Assuming that $x \leq z$ and $y \leq z$ for some $z \in R$, we
get 
$$x - y = (z + e^{+}(x)) - (z + e^{+}(y)) = e^{+}(x) + e^{+}(y) +
e^{+}(z) \in E^{+},$$ 
i.e. the elements $x$ and $y$ satisfy condition
(iii) above, so (i) and (ii) are satisfied too. Consequently,
$$x + y - z = (z + e^{+}(x)) + (z + e^{+}(y)) - z = z + e^{+}(x \wedge y) =
x \wedge y,$$ 
therefore $x \wedge y + z = x + y + e^{+}(z) = x + y$,
and 
$$(x + y) \wedge z = z + v(x + y - z) = z + v(x \wedge y) \leq
z + e^{+}(x \wedge y) = x \wedge y.$$ 
In particular, $x + y \leq z$ implies $x + y \leq x \wedge y$, 
and hence $x + y = x \wedge y$ since $e^{+}(x + y) = e^{+}(x \wedge y)$.

Finally notice that (5) follows from (4) since 
$x \wedge y \wedge z = (x \wedge z) \wedge (y \wedge z)$,
and the pair $(x \wedge z, y \wedge z)$ is bounded above by $z$.
\ep

\begin{lem} Let $R$ be a directed cr-qsring. Then the following assertions hold.
\ben
\item[\rm (1)] The semilattice morphism $e^{+}\mid_{E^{\bullet}} :
(E^{\bullet}, \bullet = \wedge) \lra (E^{+}_{+}, +)$ is a monomorphism,
and $y \in E^{\bullet}$ whenever $\eps \leq y \leq x$ and $x \in E^{\bullet}$.
\item[\rm (2)] $e^{+}(x \bullet y) = x - y = y - x$ for all 
$x, y \in E^{\bullet}$.
\item[\rm (3)] $x \bullet y \leq x$ for all $x \in R, y \in E^{\bullet}$,
and for all $x \in R, \{y \in R\,\mid\,x \bullet y \geq x\} = \eb(x)\uparrow :=
\{y \in R\,\mid\,\eb(x) \leq y\}$, and $\{y \in R\,\mid\,x \bullet y = x\} =
\{y \in \eb(x)\uparrow\,\mid\,v(y) = \eps\}$.
In particular, $x \bullet y = x \Llra e^{\bullet}(x) \leq y$, for 
$x \in R, y \in E^{\bullet}$, and 
$x \leq y \Lra e^{\bullet}(x) \leq e^{\bullet}(y)$, for $x, y \in R$.
\item[\rm (4)] $x \bullet y \bullet (x^{- 1} \wedge y^{- 1}) = x \wedge y$
for all $x, y \in R$. In particular, $x \wedge y \leq (x^{- 1} \wedge y^{- 1})^{- 1}$
(with equality if and only if $v(x) = v(y)$),
$(x \wedge y) \bullet (x^{- 1} \wedge y^{- 1})^{- 1} \leq x \bullet y$, 
$v(x \bullet y \bullet (x^{- 1} - y^{- 1})) = v(x - y)$, and
$\eb(x \wedge y) = \eb(x^{- 1} \wedge y^{- 1})$ for all $x, y \in R$.
In addition, $\eb(x \wedge y) = \eb(x \bullet y) = \eb(x) \wedge \eb(y)$
provided $x - y \in \Ep$.
\een
\end{lem}

\bp
(1) Let $x, y \in E^{\bullet}$ be such that $e^{+}(x) \leq e^{+}(y)$,
so $e^{+}(x \bullet y) = e^{+}(x) + e^{+}(y) = e^{+}(x)$. As $x \geq x \wedge y =
x \bullet y$ by Proposition 3.1.(6), it follows that $x = x \bullet y$, as required.

Now assume that $x \in E^{\bullet}$ and $\eps \leq y \leq x$. Applying $v$ to the
previous chain of inequalities, we get $v(y) = \eps$, therefore $e^{+}(e^{\bullet}(y)) =
e^{+}(y) \leq e^{+}(x)$ by Lemma 2.6.(2), so $e^{\bullet}(y) \leq x$ as shown above.
Consequently, $e^{\bullet}(y) = x + e^{+}(e^{\bullet}(y)) = x + e^{+}(y) = y$, and
hence $y \in E^{\bullet}$ as desired. In other words, $E^{\bullet}$ is identified
through the embedding $e^{+}\mid_{E^{\bullet}}$ with an ideal of the poset $(E^{+}_{+}, \leq)$.

(2) Let $x, y \in E^{\bullet}$. As $x \bullet y = x \wedge y$, and
$e^{+}(x \bullet y) = e^{+}(x) + e^{+}(y)$, the identity $e^{+}(x \bullet y) = x - y$
follows by Lemma 3.4.(4), (i) $\Lra$ (ii).

(3) Let $x \in R, y \in E^{\bullet}$. We get 
$x \bullet y = x \bullet (e^{\bullet}(x) \bullet y) = x \bullet (e^{\bullet}(x) \wedge y)
\leq x \bullet e^{\bullet}(x) = x$. If $e^{\bullet}(x) \leq y$ then $x \leq x \bullet y$
by multiplication with $x$, and hence $x \bullet y = x$. Conversely, assuming 
that $x \bullet y = x$, it follows by multiplication with $x^{- 1}$ that
$e^{\bullet}(x) \wedge y = e^{\bullet}(x) \bullet y = e^{\bullet}(x)$, so 
$e^{\bullet}(x) \leq y$.

Assuming now that $x \leq y$, it follows by Lemma 2.6.(2, 4) that 
$$\epl(x \bullet \eb(y)) = \epl(x) + x \bullet \epl(\eb(y)) = \epl(x) + (y + \epl(x))
\bullet \epl(\eb(y)) \geq$$ 
$$\epl(x) + y \bullet \epl(\eb(y)) + \epl(x) \bullet
\epl(\eb(y)) = \epl(x) + \epl (y) = \epl(x),$$ 
and hence $x \bullet \eb(y) = x$, i.e.
$\eb(x) \leq \eb(y)$, since $x \bullet \eb(y) \leq x$, as shown above. 

(4) By Proposition 3.1.(5) and statement (3) above,
it follows that $x \bullet y \bullet (x^{- 1} \wedge y^{- 1}) =
x \bullet \eb(y) \wedge y \bullet \eb(x) \leq x \wedge y$, 
in particular, $\eb(x \wedge y) = \eb(x^{- 1} \wedge y^{- 1}) \leq 
\eb(x \bullet y)$ by (3) again. We have
to show that $\epl(x \bullet y \bullet (x^{- 1} \wedge
y^{- 1})) = \epl(x \wedge y)$. By Lemma 2.6.(4) and Lemma 3.4.(1, 2),
we get 
$$\epl(x \bullet y \bullet
((x^{- 1} \wedge y^{- 1}) = \epl(x \bullet y) \bullet (x^{- 1} \wedge
y^{- 1}) + x \bullet y \bullet \epl(x^{- 1} \wedge y^{- 1}) =$$
$$\epl(x \bullet y) \bullet (v(x) \vee v(y))^{- 1} + v(x \bullet y \bullet
(x^{- 1} - y^{- 1})) =$$
$$(v(x) + v(y)) \bullet (\epl(\eb(x \bullet y)) +
\epl(\eb(x^{- 1} \wedge y^{- 1}))) =$$
$$(v(x) + v(y)) \bullet \epl(\eb(x \wedge y))
= v(x - y) = \epl(x \wedge y),$$
as required. In particular, we obtain the identity $v(x \bullet y 
\bullet (x^{- 1} - y^{- 1})) = v(x - y)$. By multiplication with 
$(x^{- 1} \wedge y^{- 1})^{- 1}$, it follows that 
$$(x \wedge y) \bullet (x^{- 1} \wedge y^{- 1})^{- 1} = x \bullet
y \bullet \eb(x \wedge y) \leq x \bullet y,$$ 
as desired. On the other hand, 
$x \wedge y \leq (x^{- 1} \wedge y^{- 1})^{- 1}$ since 
$$x \wedge y \wedge (x^{- 1} \wedge y^{- 1})^{- 1} = x \wedge y \wedge
(x \bullet y \bullet (x \wedge y)^{- 1}) = 
((x \wedge y)^2 \bullet (x \wedge y)^{- 1})
\wedge (x \bullet y \bullet (x \wedge y)^{- 1}) = $$
$$(x \wedge y)^{- 1} \bullet 
((x \wedge y)^2 \wedge (x \bullet y)) = (x \wedge y)^{- 1} \bullet 
(x \wedge y)^2 = x \wedge y.$$
Finally, assuming that $x - y \in \Ep$, we obtain $\eb(x \wedge y) =
\eb(x \bullet y)$ since \\ $\eb(x \wedge y) \leq \eb(x \bullet y)$ and
$$\epl(\eb(x \wedge y)) = \frac{\epl(x \wedge y)}{v(x \wedge y)} =
\frac{\epl(x) + \epl(y)}{v(x) + v(y)} = \frac{(\epl(x) + \epl(y))
\bullet (v(x) \vee v(y))}{v(x \bullet y)} \geq$$ 
$$\frac{\epl(x) \bullet v(y) + \epl(y) \bullet v(x)}
{v(x \bullet y)} = \frac{\epl(x \bullet y)}
{v(x \bullet y)} = \epl(\eb(x \bullet y))$$
by Lemmas 2.6.(4) and 3.4.(4).
\ep

\begin{lem}
Let $R$ be a directed cr-qsring. Then the following assertions hold.
\ben
\item[\rm (1)] There exists the join $x \vee y$ provided the 
elements $x, y \in R$ are {\em incident}, i.e. the pair $(x, y)$ 
is bounded above with respect to the partial order $\leq$. 
Set by convention $x \vee y = \infty$ whenever 
the elements $x$ and $y$ are not incident.
\item[\rm (2)] $e^{+}(x \vee y) = e^{+}(x) \vee e^{+}(y)$ and
$v(x \vee y) = v(x) \vee v(y)$ provided $x \vee y \neq \infty$.
\item[\rm (3)] The partial binary operation $\vee$ is compatible
with the operations $+, -$, and $\bullet$, i.e. 
for all $x, y, z \in R$ such that $y \vee z \neq \infty$, the following
identities hold.
\ben
\item[\rm (i)] $(x + y) \vee (x + z) = x + (y \vee z)$;
\item[\rm (ii)] $(- y) \vee (- z) = - (y \vee z)$;
\item[\rm (iii)] $(x \bullet y) \vee (x \bullet z) = x \bullet (y \vee z)$.
\een
\item[\rm (4)] $(x \wedge y) \vee (x \wedge z) = x \wedge (y \vee z)$
for all $x, y, z \in R$, with $y \vee z \neq \infty$.
In particular, for all $x \in R, x\downarrow$
is a distributive lattice with a last element.
\item[\rm (5)] For all $x, y \in R$, satisfying $x \vee y \neq \infty$,
$$x \bullet y = (x \wedge y) \bullet (x \vee y),\q
\eb(x \wedge y) = \eb(x \bullet y) = \eb(x) \wedge \eb(y),$$
$$\eb(x \vee y) = \eb(x) \vee \eb(y),\q{\rm and}\q
(x^{- 1} \wedge y^{- 1})^{- 1} = (x \vee y) \bullet 
\eb(x \bullet y) \leq x \vee y.$$
In particular, $x \vee y \in E^{\bullet}$ provided 
$x, y \in E^{\bullet}$ and $x \vee y \neq \infty$.
\item[\rm (6)] For all $x, y \in R$, the following assertions are equivalent.
\ben
\item[\rm (i)] $(x^{- 1} \wedge y^{-1})^{- 1} = x \vee y$.
\item[\rm (ii)] $x - y \in \Ep$ and $\eb(x) = \eb(y)$.
\item[\rm (iii)] $\eb(x \wedge y) = \eb(x) = \eb(y)$.
\een
\item[\rm (7)] For all $x, y \in R$ satisfying $x \wedge y \in \Ep$,
i.e. $x \wedge y = v(x) + v(y)$, the following assertions are equivalent.
\ben
\item[\rm (i)] $x \vee y \neq \infty$.
\item[\rm (ii)] $x - y \in \Ep$.
\item[\rm (iii)] $x - y = x \wedge y$.
\een
In particular, for all $x \in R, y \in \Ep, x \vee y \neq \infty
\Llra \epl(x) + y \leq v(x)$, and $x \vee y \neq \infty$ for
all $x, y \in \Eb$ which are orthogonal, i.e. $x \bullet y = \eps$. 
\een
\end{lem}

\bp
(1) Let $x, y, z \in R$ be such that $x \leq z$ and
$y \leq z$. We show that the element $u := z + (e^{+}(x)
\vee e^{+}(y))$ is the join of the pair $(x, y)$.
Obviously, $e^{+}(u) = e^{+}(z) + (e^{+}(x) \vee e^{+}(y)) =
e^{+}(x) \vee e^{+}(y)$. As $x \leq z$, we get $x = z + e^{+}(x) =
z + e^{+}(x) + (e^{+}(x) \vee e^{+}(y)) = u + e^{+}(x)$, so
$x \leq u$, and similarly, $y \leq u$. Assuming that $x \leq t$
and $y \leq t$ for some $t \in R$, in particular, $e^{+}(u) \leq
e^{+}(t)$, it remains to check that $u \leq t$. Since $x, y \leq z$,
it follows that $x, y \leq z \wedge t$, so we may assume without 
loss that $t \leq z$. Then we get
$t + e^{+}(u) = z + e^{+}(t) + e^{+}(u) = z + e^{+}(u) = u$, 
and hence $u \leq t$ as desired.

(2) Let $x, y \in R$ be such that $x \vee y \neq \infty$.
We already know from (1) that $e^{+}(x \vee y) = e^{+}(x) \vee e^{+}(y)$,
so it remains only to evaluate $v(x \vee y)$. As $x, y \leq x \vee y$,
it follows that $v(x) \vee v(y) \leq v(x \vee y)$. On the other hand,
setting $u = x \vee y$, we get $v(x) = v(u + e^{+}(x)) \geq v(u) + 
e^{+}(x)$, and similarly, $v(y) \geq v(u) + e^{+}(y)$, therefore
$v(x) \vee v(y) \geq v(u) + (e^{+}(x) \vee e^{+}(y)) = v(u) + e^{+}(u) =
v(u)$, so $v(u) = v(x) \vee v(y)$ as desired.

(3) Let $x, y, z \in R$ be such that $u := y \vee z \neq \infty$.
As $y, z \leq u$, it follows that $x + y, x + z \leq x + u$ and
$x \bullet y, x \bullet z \leq x \bullet u$, and hence 
$(x + y) \vee (x + z) \leq x + u$ and $(x \bullet y) \vee (x \bullet z)
\leq x \bullet u$. Since
$e^{+}((x + y) \vee (x + z)) = e^{+}(x + y) \vee e^{+}(x + z) =
e^{+}(x) + (e^{+}(y) \vee e^{+}(z)) = e^{+}(x + u)$, we conclude that
$(x + y) \vee (x + z) = x + (y \vee z)$. On the other hand, we obtain
$x \bullet y = x \bullet (u + e^{+}(y)) \geq x \bullet u + x \bullet
e^{+}(y)$, and similarly, $x \bullet z \geq x \bullet u + x \bullet e^{+}(z)$,
therefore $x \bullet y \vee x \bullet z \geq (x \bullet u + x \bullet
e^{+}(y)) \vee (x \bullet u + x \bullet e^{+}(z)) = x \bullet u +
(x \bullet e^{+}(y) \vee x \bullet e^{+}(z)) = x \bullet u + 
x \bullet e^{+}(u) = x \bullet u$ since $e^{+}(x \bullet u) \leq 
x \bullet e^{+}(u)$.

The proof of the identity (ii) is straightforward.

(4) Let $x, y, z \in R$ be such that 
$u := y \vee z \neq \infty$. As $y, z \leq u$, it follows that
$(x \wedge y) \vee (x \wedge z) \leq x \wedge u$. We get equality
since $e^{+}(x \wedge u) = v(u - x) = v((y - x) \vee (z - x)) =
v(y - x) \vee v(z - x) = e^{+}(x \wedge y) \vee e^{+}(x \wedge z) =
e^{+}((x \wedge y) \vee (x \wedge z))$.

(5) Let $x, y \in R$ be such that $u := x \vee y \neq 
\infty$. The identity $x \bullet y = (x \wedge y) \bullet u$ is
an immediate consequence of the compatibility of the 
multiplication with $\wedge$ and $\vee$, while the identity
$\eb(x \wedge y) = \eb(x \bullet y)$ follows by Lemma 3.5.(4),
since $x \vee y \neq \infty \Lra x - y \in \Ep$.

On the other hand, $\eb(x), \eb(y) \leq \eb(u)$ by Lemma 3.5.(3),
therefore \\
$t := \eb(x) \vee \eb(y) \neq \infty$ and $t \leq \eb(u)$.
To obtain the identity $t = \eb(u)$, we have to show that
$\epl(\eb(u)) \leq \epl(t)$. By Lemma 2.6.(2, 4) and statement (2) 
of the lemma, we get 
$$\epl(\eb(u)) = \epl(u) \bullet v(u)^{- 1} = (\epl(x) \vee \epl(y)) \bullet
(v(x) \vee v(y))^{- 1} =$$
$$\epl(x) \bullet \epl(y) \bullet ((\epl(x) + \epl(y))
\bullet (v(x) \vee v(y)))^{- 1} \leq \epl(x) \bullet \epl(y) \bullet
\epl(x \bullet y)^{- 1} = $$
$$\epl(\eb(x)) \bullet \epl(\eb(y)) \bullet \epl(\eb(x \bullet y))^{- 1}
= \epl(\eb(x)) \vee \epl(\eb(y)) \leq \epl(t),$$
as required.

Using Lemma 3.5.(3, 4), we obtain
$$(x^{- 1} \wedge y^{- 1})^{- 1} = x \bullet y \bullet (x \wedge y)^{- 1} =
(x \vee y) \bullet (x \wedge y) \bullet (x \wedge y)^{- 1} =
(x \vee y) \bullet \eb(x \bullet y) \leq x \vee y.$$
Finally, for $x, y \in \Eb$ with $x \vee y \neq \infty$,
we get $\eb(x \vee y) = \eb(x) \vee \eb(y) = x \vee y$, and hence
$x \vee y \in \Eb$.

(6)
(i) $\Lra$ (ii) : Assume that $(x^{- 1} \wedge y^{- 1})^{- 1} = x \vee y$.
As $x \vee y \neq \infty$, we get $x - y \in \Ep$. On the
other hand, it follows from (5) that $(x \vee y) \bullet \eb(x \bullet y) =
x \vee y$, therefore $\eb(x \bullet y) \geq \eb(x \vee y) = \eb(x) \vee \eb(y)$
by Lemma 3.5.(3), so $\eb(x) = \eb(y)$ as desired.

(ii) $\Lra$ (iii) : Assuming that $x - y \in \Ep$ and $\eb(x) = \eb(y)$,
it follows by Lemmas 2.6.(2) and 3.4.(1, 4) that 
$$\epl(\eb(x)) \bullet (v(x) + v(y)) = \epl(x) + \epl(y) =
\epl(x \wedge y) = \epl(\eb(x \wedge y)) \bullet (v(x) + v(y)),$$
therefore $\eb(x) = \eb(x \wedge y)$ by Lemma 3.5.(1).

(iii) $\Lra$ (i) : Assuming that $\eb(x) = \eb(y) = \eb(x \wedge y)$,
it follows by Lemma 3.5.(4) that
$$x \wedge (x^{- 1} \wedge y^{- 1})^{- 1} = x \bullet (x^{- 1} \wedge y^{- 1})^{- 1}
\bullet (x^{- 1} \wedge y^{- 1}) = x \bullet \eb(x \wedge y) = 
x \bullet \eb(x) = x,$$
i.e. $x \leq (x^{- 1} \wedge y^{- 1})^{- 1}$, and, similarly,
$y \leq (x^{- 1} \wedge y^{- 1})^{- 1}$. Consequently, 
$x \vee y \neq \infty$, so $x \vee y = (x^{- 1} \wedge y^{- 1})^{- 1}$
by (5).

(7) (i) $\Lra$ (ii) and (ii) $\Llra$ (iii) are obvious, so
it remains to prove the implication (ii) $\Lra$ (i). It suffices 
to show that $z := x \vee \epl(y) \neq \infty$ and $u := y \vee \epl(x)
\neq \infty$, since then we obtain $x = x + \epl(x) \leq z + u,
y = \epl(y) + y \leq z + u$, and $\epl(z + u) = \epl(z) = \epl(u) =
\epl(x) \vee \epl(y)$, and hence $x \vee y = z + u \neq \infty$ as
desired. Thus we are reduced to the case when $y \in \Ep$ and 
$\epl(x) + y = v(x) + y \leq v(x)$. Let us show that $x \vee y = t^{- 1}$,
where $t := x^{- 1} + \del$, with $\del := \frac{\epl(\eb(x))}{v(x) \vee y}$.
We get $\epl(t) = \del, v(t) = v(x)^{- 1} + \del = (v(x) \vee y)^{- 1}, 
\epl(\eb(t)) = \epl(\eb(x))$, and hence 
$\eb(t) = \eb(x), y \leq v(t^{- 1}) \leq t^{- 1}$ and 
$\epl(t^{- 1}) = \epl(x) \vee y$. To conclude that $t^{- 1} = x \vee y$,
it remains to note that  
$$x \wedge t^{- 1} = x \bullet t^{- 1} \bullet (x^{- 1} \wedge t) =
x \bullet t^{- 1} \bullet t = x \bullet \eb(t) = x \bullet \eb(x) = x.$$
\ep

\begin{rems} \em
(1) Since the directed cr-qsrings form a
variety of algebras of signature $(2, 2, 1, 1, 0)$, any
substructure $R'$ of a directed cr-qsring $R$ is a directed
cr-qsring, in particular, $x \wedge y \in R'$ provided 
$x, y \in R'$. However, for $x \in R, y \in R', x \leq y$ does 
not imply $x \in R'$, and also, for $x, y \in R', \infty \neq
x \vee y \in R$ does not imply $x \vee y \in R'$ (see Remark 4.9.)

(2) As in the case of commutative regular rings, the map 
$E^{\bullet} \times E^{\bullet} \lra E^{\bullet}, (x, y)
\mapsto x + y - x \bullet y$ is well defined
in any directed cr-qsring $R$, but by contrast with the 
former case, where the map above defines the join $x \vee y$
of any pair of idempotents $(x, y)$, in the latter case we
obtain $x + y - x \bullet y = (x - x \bullet y) +
(y - x \bullet y) + x \bullet y = x \bullet y = x \wedge y$
for all $x, y \in E^{\bullet}$ since $x - x \bullet y = 
y - x \bullet y = e^{+}(x \bullet y)$ by Lemma 3.5.(2). 
Notice also that in an arbitrary 
cr-qsring $R$, we have only the inequality 
$x + y - x \bullet y \leq (x + y - x \bullet y)^{2}$ for
all $x, y \in E^{\bullet}$, so, in general, $E^{\bullet}$ 
is not necessarily closed under the map above.
\end{rems}


\section{The metric structure of a directed commutative
regular quasi-semiring}

$\q$ Let $R$ be a directed cr-qsring. As we have shown in
Section 2, the subset $E^{+}$ of the idempotents of
the commutative regular semigroup $(R, +)$ has a natural
structure of Abelian $l$-group with multiplication $\bullet$
as the group operation, neutral element $\eps$,
partial order $\leq$, and addition $+$ 
as the corresponding meet-semilattice operation $\wedge$.
Let us denote by $\Lam$ this structure of Abelian
$l$-group, and by $\Lam_{+} = \{\alpha \in \Lam\,\mid\,\alpha \geq \eps\}$,
the monoid of the nonnegative elements of $\Lam$. For
any $\alpha \in \Lam$, set $\alpha_{+} := \alpha \vee \eps =
(\alpha^{- 1} + \eps)^{- 1}, \alpha_{-} := (\alpha^{- 1})_{+}$,
and $\mid \alpha \mid := \alpha \vee \alpha^{- 1} = (\alpha +
\alpha^{- 1})^{- 1} = \alpha_{+} \bullet \alpha_{-}$.
In this section we shall define a distance map on the
directed cr-qsring $R$ with values in $\Lam_{+}$, and we
shall investigate the main properties of this metric
structure. 
\smallskip

For any pair $(x, y)$ of elements of $R$, set 
$$[x, y] := \{z \in R\,\mid\,x \wedge y \leq z = (x \wedge z) \vee
(y \wedge z)\}$$
Notice that $[x, x] = \{x\}, [x, y] = [y, x], x, y, x \wedge y \in [x, y]$,
and $x \vee y \in [x, y]$ provided $x \vee y \neq \infty$. Notice
also that $[x, y] = \{z\,\mid\,x \leq z \leq y\}$ (the interval)
whenever $x \leq y$. Thus $[x, y]$ coincides with the interval
$[x \wedge y, x \vee y]$ provided $x \vee y \neq \infty$.

Call {\em cell} or {\em simplex} any subset of $R$ of the form
$[x, y]$. Given a cell $C \sse R$, any $x \in R$ for which 
there exists $y \in R$ such that $C = [x, y]$ is called an
{\em end} of the cell $C$. The (non-empty) subset of all ends of
a cell $C$ is denoted by $\partial C$ and called the {\em boundary} of
$C$.

To provide equivalent descriptions for the cells as introduced above,
we define two maps $\lam : R \times R \lra \Lam_{+}$ and
$d : R \times R \lra \Lam_{+}$ as follows :
$$\lam(x, y) := \frac{e^{+}(x)}{e^{+}(x \wedge y)} = \frac{e^{+}(x)}{v(x - y)} =
e^{+}(x) \bullet (x - y)^{- 1}$$
and
$$d(x, y) = \lam(x, y) \bullet \lam(y, x) = 
\frac{e^{+}(x) \bullet e^{+}(y)}{(e^{+}(x \wedge y))^{2}}$$

Notice that $\lam(x, y) = (x \bullet y^{- 1})_+$ and $d(x, y) =
\mid x \bullet y^{- 1} \mid$ provided $x, y \in E^+$.

The next lemmas collect some basic properties of the maps $\lam$ and $d$.

\begin{lem}
\ben
\item[\rm (1)] $\lam(x, y) = \eps \Llra x \leq y$.
\item[\rm (2)] $\lam(x, y) = \lam(x, x \wedge y) = d(x, x \wedge y)$.
\item[\rm (3)] The ternary map $\wh{\lam} : R^{3} \lra \Lam_{+}$,
defined by 
$$\wh{\lam}(x, y, z) = \lam(x, y) \bullet \lam(y, z) \bullet \lam(z, x) =
\frac{e^{+}(x) \bullet e^{+}(y) \bullet e^{+}(z)}
{e^{+}(x \wedge y) \bullet e^{+}(y \wedge z) \bullet e^{+}(z \wedge x)}$$
is symmetric in the variables $x, y, z$.
\item[\rm (4)] $\lam(x, y) \leq \lam(x, z) \bullet \lam(z, y)$ for all
$x, y, z \in R$.
\item[\rm (5)] $\lam(x, y) = \lam(v(x), v(y)) \bullet d(\eb(x), 
\eb(x \wedge y))$. In particular, 
$\lam(x, y) = \lam(v(x), v(y)) = v(x \bullet y^{- 1})_+$
if and only if $\eb(x) = \eb(x \wedge y)$.
\item[\rm (6)] $d(x, y) = d(v(x), v(y)) \bullet d(\eb(x), \eb(y))
\bullet d(\eb(x \bullet y), \eb(x \wedge y))^2$. In particular,
$d(x, y) = d(v(x), v(y)) = \mid v(x \bullet y^{- 1}) \mid$ 
if and only if $\eb(x) = \eb(y) = \eb(x \wedge y)$, and
$d(x, y) = d(\epl(x), \epl(y)) = \mid \epl(x) \bullet \epl(y)^{- 1} \mid$
provided $x, y \in \Eb$.
\een 
\end{lem}

\bp
The statemements (1)-(3) are obvious. Notice that (4) is equivalent 
with the inequality
$$e^{+}(x \wedge z) \bullet e^{+}(y \wedge z) \leq e^{+}(z) \bullet
e^{+}(x \wedge y)$$
Indeed, setting $u := (x \wedge z) \vee (y \wedge z) \leq z$, it follows
by Lemma 3.4.(4) and Lemma 3.6.(2) that
$${\rm LHS} = (e^{+}(x \wedge z) \vee e^{+}(y \wedge z)) \bullet
(e^{+}(x \wedge z) + e^{+}(y \wedge z)) = e^{+}(u) \bullet 
e^{+}(x \wedge y \wedge z) \leq {\rm RHS}$$
as desired.

(5) By Lemmas 2.6.(2), 3.4.(1, 2), 3.5.(3), we get
$$\lam(x, y) = \frac{v(x) \bullet \epl(\eb(x))}{(v(x) + v(y)) \bullet \epl(\eb(x \wedge
y))} = \lam(v(x), v(y)) \bullet d(\eb(x), \eb(x \wedge y)),$$
as required. The necessary and sufficient condition to have the
equality $\lam(x, y) = \lam(v(x), v(y))$ is
immediate by Lemma 3.5.(1). Finally, notice that (5) $\Lra$ (6).
\ep 

\begin{lem}
\ben
\item[\rm (1)] Any directed cr-qsring $R$ is a $\Lam$-{\em metric 
space} with the $\Lam$-{\em valued distance map} $d : R \times R \lra
\Lam$ satisfying
\ben
\item[\rm (i)] $d(x, y) = \eps \Llra x = y$,
\item[\rm (ii)] $d(x, y) = d(y, x)$ for all $x, y \in R$, and
\item[\rm (iii)] \textbf{Triangle inequality} : $d(x, y) \leq
d(x, z) \bullet d(z, y)$ for all $x, y, z \in R$. 
\een
\item[\rm (2)] $d(x, y) = d(x, x \wedge y) \bullet d(y, x \wedge y)$
for all $x, y \in R$.
\een
\end{lem}

\bp
The triangle inequality follows by Lemma 4.1.(4), while the rest of 
assertions are obvious.
\ep  .

\begin{lem} Let $R$ be a directed cr-qsring. Then the following assertions hold.
\ben
\item[\rm (1)] For all $x, y \in R$,
$$[x, y] = \{z \in R\,\mid\,\lam(x, z) \bullet \lam(z, y) =
\lam(x, y)\} = \{z \in R\,\mid\,d(x, z) \bullet d(z, y) = d(x, y)\}.$$
\item[\rm (2)] The ternary relation $B(x, y, z) \Llra y \in [x, z]$ is
a {\em betweenness relation}, i.e. the following hold for all $x, y, z, u
\in R$.
\ben
\item[\rm (i)] $B(x, x, y)$,
\item[\rm (ii)] $B(x, y, x) \Lra x = y$, and
\item[\rm (iii)] $B(x, y, z), B(x, u, y) \Lra B(z, u, x)$.
\een
\item[\rm (3)] The betweenness relation $B$ is compatible with the 
semilattice operation $\wedge$, i.e. $B(x, y, z) \Lra B(x \wedge u,
y \wedge u, z)$ for all $x, y, z, u \in R$.
\item[\rm (4)] For all $a \in R$, the map $a\downarrow \lra \Lam_{+},
x \mapsto d(a, x) = \frac{\epl(a)}{\epl(x)}$ is an antiisomorphism
of distributive lattices, with the inverse $\alpha \in \Lam_{+}
\mapsto a + \frac{\epl(a)}{\alpha} \in a\downarrow$.
\een
\end{lem}

\bp
(1) Let us denote by $M$ the set involving the map $\lam$, and by $N$ the
set involving $d$. Let $M'$ be the set obtained from $M$ by interchanging
the elements $x$ and $y$, and notice that $M \cap M' = N$ by Lemma 4.1.(4)
and the definition of $d$. Consequently, it suffices to prove that
$[x, y] = M$ since $[x, y] = [y, x]$. However the equality $[x, y] = M$ 
follows easily with the same argument as in the proof of Lemma 4.1.(4).

(2) is straightforward, while (3) and (4) follow by Lemma 3.6.(4).
\ep

\begin{lem}
\ben
\item[\rm (1)] $d(\varphi(x), \varphi(y)) \leq d(x, y)$, 
where $\varphi(x)$ is defined by one of the following formulas, 
for some $a \in R$ :
\ben
\item[\rm (i)] $\varphi(x) = a + x$;
\item[\rm (ii)] $\varphi(x) = a \bullet x$;
\item[\rm (iii)] $\varphi(x) = a \wedge x$;
\item[\rm (iv)] $\varphi(x) = a \vee x$ provided $a \vee x \neq \infty$.
\een
\item[\rm (2)] $\lam(- x, - y) = \lam(x, y)$ and 
$d(- x, - y) = d(x^{- 1}, y^{- 1}) = d(x, y)$ for all $x, y \in R$.
\een
\end{lem}

\bp
(1) Since $x \geq y \Lra \varphi(x) \geq \varphi(y)$, it suffices 
to prove the inequality $d(\varphi(x), \varphi(y))
\leq d(x, y)$ in the case $x \geq y$, and hence we have to check
that $\epl(\varphi(x)) \bullet \epl(y) \leq \epl(\varphi(y)) \bullet
\epl(x)$ under the assumption $x \geq y$. The verification in
the cases (i) and (iv) is straightforward. In the case (ii), 
by Lemma 2.6.(2, 4, 6), we obtain
$$\epl(a \bullet x) \bullet \epl(y) = \epl(x) \bullet \epl(y)
\bullet v(a) \bullet (\eps + \frac{\epl(\eb(a))}{\epl(\eb(x))}) \leq$$ 
$$\epl(x) \bullet \epl(y) \bullet v(a) \bullet 
(\eps + \frac{\epl(\eb(a))}{\epl(\eb(y))}) = \epl(a \bullet y) \bullet \epl(x),$$
as desired. In the case (iii), $\epl(a \wedge y) \vee \epl(y) \leq \epl(x)$
implies by Lemma 3.4.(5) that 
$$\epl(a \wedge x) \bullet \epl(y) \leq (\epl(a \wedge x) + \epl(y)) \bullet
\epl(x) = (\epl(a \wedge x) + \epl (x \wedge y)) \bullet \epl(x) =$$
$$\epl(a \wedge x \wedge y) \bullet \epl(x) = \epl(a \wedge y) \bullet \epl(x),$$
as required.

(2) The identities $\lam(- x, - y) = \lam(x, y)$ and $d(- x, - y) = d(x, y)$ 
are obvious since $(- x) \wedge (- y) = -(x \wedge y)$ by Lemma 3.4.(3). By
Lemma 2.4.(5), Lemma 2.6.(2), Lemma 3.4.(1) and Lemma 3.5.(4), we get
$$\lam(x^{- 1}, y^{- 1}) = \frac{\epl(x^{- 1})}{\epl(x^{- 1} \wedge y^{- 1})} =
\frac{v(x)^{- 2} \bullet \epl(x)}{v(x \bullet y)^{- 1} \bullet \epl(x \wedge y)} =
v(x^{- 1} \bullet y) \bullet \lam(x, y),$$
therefore $d(x^{- 1}, y^{- 1}) = d(x, y)$ as desired. 
\ep

\begin{rem} \em 
By Lemma 4.4.(2), the commuting involutions $x \mapsto - x$ and
$x \mapsto x^{- 1}$ are both automorphisms of the $\Lam$-metric
space $(R, d)$. Notice that the involution $x \mapsto - x$ is an
automorphism of the structure $(R, \wedge, B)$, while the involution
$x \mapsto x^{- 1}$ is only an automorphism of the structure
$(R, B)$, and also an isomorphism from the structure $(R, \wedge, B)$
onto the structure $(R, \sqcap, B)$, where the semilattice operation
$\sqcap$, defined by $x \sqcap y := (x^{- 1} \wedge y^{- 1})^{- 1}$,
is the meet of any pair $(x, y)$ with respect to the partial order 
$x \sqsubseteq y \Llra x^{- 1} \leq y^{- 1} \Llra x \leq x^2 \bullet y^{- 1}$. 
The semilattice operations $\wedge$ and $\sqcap$ are both 
compatible with the betweenness relation $B$, and, by Lemmas 3.5.(4) and 3.6.(6), 
$x \wedge y = x \sqcap y \Llra v(x) = v(y)$, and $x \sqcap y = x \vee y
\Llra \eb(x \wedge y) = \eb(x) = \eb(y)$.
Notice also that the corresponding partial join operation $x \sqcup y$ is
defined if and only if $x^{- 1} \vee y^{- 1} \neq \infty$, and
$x \sqcup y = (x^{- 1} \vee y^{- 1})^{- 1}$.
\end{rem}


\subsection{Congruences on directed commutative regular quasi-semirings}

\begin{lem}
The map $\,\equiv\,\, \mapsto \{\alpha \in \Lam_{+}\,\mid\,\alpha \equiv \eps\}$
is an isomorphism from the lattice of all congruences of the directed
cr-qsring $R$ onto the lattice of all convex submonoids of $\Lam_{+}$.
Its inverse sends a convex submonoid $S \sse \Lam_{+}$ to the congruence
relation $\{(x, y) \in R \times R\,\mid\,d(x, y) \in S\}$. 
\end{lem}

\bp
Let $\equiv$ be a congruence of the directed cr-qsring $R$. 
The relation induced by $\equiv$ on the substructure $\Lam$ is
obviously a congruence on the commutative $l$-group $\Lam$,
and hence $\{\alpha \in \Lam_{+}\,\mid\,\alpha \equiv \eps\}$ is
a convex submonoid of $\Lam_{+}$. First we
have to show that for all $x, y \in R$, $x \equiv y \Llra d(x, y) 
\equiv \eps$, i.e. $\epl(x) \bullet \epl(y) \equiv \epl(x \wedge y)^{2}$.
The implication $\Lra$ is a consequence of the obvious implications
$x \equiv y \Lra x \equiv x \wedge y$, $x \equiv y \Lra \epl(x) \equiv
\epl(y)$, and $(x \equiv z, y \equiv z) \Lra x \bullet y
\equiv z^{2}$. Conversely, assuming that $d(x, y) \equiv \eps$, it
follows that $\lam(x, y) \equiv \lam(y, x) \equiv \eps$ since
$\eps \leq \lam(x, y), \lam(y, x) \leq d(x, y)$. Consequently,
$\epl(x) \equiv \epl(x \wedge y) \equiv \epl(y)$, and hence
$x = x + \epl(x) \equiv x + \epl(x \wedge y) = x \wedge y$, and,
similarly, $y \equiv x \wedge y$, so $x \equiv y$ as desired.

Next, assuming that $S$ is a convex submonoid of $\Lam_{+}$,
we deduce that the binary relation 
$\{(x, y) \in R \times R\,\mid\,d(x, y) \in S\}$
is a congruence of the directed cr-qsring $R$ thanks to 
Lemmas 4.2. and 4.4.
\ep

\begin{lem}
 Let $S = \cS(\Eb)$ be the smallest convex submonoid of $\Lam_+$
containing $\epl(\Eb)$, $G = \cG(\Eb) = \{\gam \in \Lam\,|\,
| \gam | \in S \}$ the corresponding convex $l$-subgroup
of $\Lam$, and $\equiv\, =\, \equiv_{\Eb}$ the congruence on the
directed cr-qsring $R$ defined by $S$. Then the following 
assertions hold.
\ben
\item[\rm (1)] $S = \bigcup_{x \in \Eb, n \geq 1} [\eps, \epl(x)^n]$
provided $x \vee y \neq \infty$ for all $x, y \in \Eb$.
\item[\rm (2)] For all $x, y \in R$, $x \equiv y \Llra v(x \bullet
y^{- 1}) \in G$.
\item[\rm (3)] $G$ is the smallest convex $l$-subgroup of $\Lam$
for which the surjective map \\ $\wt{v} : R \lra \Lam/G$ induced
by the $q$-valuation $v : R \lra \Lam$ is a morphism of cr-qsrings.
\item[\rm (4)] $G$ is the smallest convex $l$-subgroup of $\Lam$
for which ${\rm Ker}(\wt{v}) := \{x \in R\,\mid\,v(x) \in G\}$ is 
a directed cr-qsring, a substructure of $R$. Moreover $x \vee y \in 
\Ker(\wt{v})$ provided $x, y \in \Ker(\wt{v})$ and $x \vee y \neq
\infty$.
\item[\rm (5)] $\equiv$ is the smallest congruence on the directed
cr-qsring $R$ for which the quotient $R/\equiv$ is an Abelian $l$-group.
\een 
\end{lem}

\bp
(1) Let $\alpha, \beta \in \Lam_+$ be such that $\alpha \leq \epl(x)^n,
\beta \leq \epl(y)^m$ for $x, y \in \Eb, n, m \geq 1$. By assumption,
$z := x \vee y \in \Eb$, therefore $\alpha \bullet \beta \leq \epl(z)^{n + m}$,
as required.

(2) For all $x, y \in \Eb$, we get $d(x, y) \leq d(x, \eps) \bullet
d(y, \eps) = \epl(x) \bullet \epl(y) \in S$, and hence $d(x, y) \in S$,
and $x \equiv \eps$ for all $x \in \Eb$. Then, by Lemma 4.1.(6), it
follows that for all $x, y \in R$, 
$$x \equiv y \Llra d(x, y) \in S \Llra d(v(x), v(y)) = \mid v(x \bullet
y^{- 1}) \mid \in S \Llra v(x \bullet y^{- 1}) \in G,$$
as desired.

(3) By (2), the quotient directed cr-qsring $R/{\equiv}$ is isomorphic
to the quotient Abelian $l$-group $\Lam/G$, and the surjective
map $\wt{v} : R \lra \Lam/G$ induced by $v$ is obviously a morphism
of cr-qsrings. On the other hand, let $H$ be a convex $l$-subgroup
of $\Lam$ with the property that the surjective map 
$\wt{v}_H : R \lra \Lam/H$ induced by $v$ is a morphism of cr-qsrings.
As $v(x) = \eps \in H$ provided $x \in \Eb$, it follows that 
$\wt{v}_H(\epl(x)) = \wt{v}_H(x - x) = \wt{v}_H(x) - \wt{v}_H(x) =
\eps\,{\rm mod}\,H$, i.e. $\epl(x) = v(\epl(x)) \in H$ for all
$x \in \Eb$, and hence $G \sse H$, as required.

(4) Obviously, $\eps \in \Ker(\wt{v})$, and $\Ker(\wt{v})$ is
closed under multiplication and the unary operations $x \mapsto - x$
and $x \mapsto x^{- 1}$, so it remains to note that $v(x - y) =
\epl(x \wedge y) = \epl(\eb(x \wedge y)) \bullet (v(x) + v(y)) \in G$,
i.e. $x - y \in \Ker(\wt{v})$, provided $x, y \in \Ker(\wt{v})$.
Note that $\Ep(\Ker(\wt{v})) = G, \Eb(\Ker(\wt{v})) = \Eb$,
and the convex monoid $S = G_+$ is obviously generated by $\epl(\Eb)$.
On the other hand, let $H$ be a convex $l$-subgroup of $\Lam$ such
that $\Ker(\wt{v}_H) = \{x \in R\,\mid\,v(x) \in H\}$ is a 
substructure of $R$. With the same argument as in the proof of (3),
it follows that $\epl(x) = v(\epl(x)) \in H$ for all $x \in \Eb$,
and hence $G \sse H$ as desired. Notice also that for any such $H$, 
in particular for $G$, $x \vee y \in \Ker(\wt{v}_H)$ provided 
$x, y \in \Ker(\wt{v}_H)$ and $x \vee y \neq \infty$, since 
$v(x \vee y) = v(x) \vee v(y) \in H$ by Lemma 3.6.(2).

Finally notice that (5) is a reformulation of (3).
\ep


\subsection{Superrigid directed commutative regular quasi-semirings}

$\q$ Call {\em rigid} a directed cr-qsring $R$ for which 
$\cS(\Eb) = \Lam_+$, i.e., by Lemma 4.7., $R$ does 
not have proper Abelian $l$-group quotients.
In particular, $R$ is rigid provided the injective map 
$\epl \mid_{\Eb} : \Eb \lra \Lam_+$ is surjective. Call such an $R$
{\em superrigid}, and for any $\alpha \in \Lam_+$, denote by $1_\alpha$ 
the unique element of $\Eb$ satisfying $\epl(1_\alpha) = \alpha$. 
Moreover, for any $\alpha \in \Lam$, put $1_\alpha := 1_{\alpha_+} +
\alpha = \displaystyle\lim_{\stackrel{\gam \rightarrow \infty}{}}
(1_\gam + \alpha)$, and notice that $\epl(1_\alpha) = \alpha$,
$v(1_\alpha) = \alpha_-^{- 1} \leq \eps, \eb(1_\alpha) = 1_{\alpha_+}$,
$1_\alpha^{- 1} = 1_{\alpha_+} \vee \alpha_- \geq \eps$ (by Lemma 3.6.(7)), 
and the map $\alpha \mapsto 1_\alpha$ is an isomorphism 
$(\Lam, +, \vee) \lra (\{1_\alpha\,|\,\alpha \in \Lam\}, \wedge,
\vee)$ of distributive lattices.

By Lemma 4.7., any directed cr-qsring can be seen as an {\em extension 
of an Abelian $l$-group by a rigid directed cr-qsring}.
On the other hand, the category of Abelian $l$-groups 
is equivalent to a full subcategory of 
the category of all superrigid directed cr-qsrings, as follows :

\begin{pr}
The forgetful functor, sending a directed cr-qsring $R$ to the
associated Abelian $l$-group $\Lam = (\Ep, \bullet, +, \vee)$, induces by
restriction an equivalence from the full subcategory of those
superrigid directed cr-qsrings $R$ which satisfy the supplementary condition
$$\forall x, y \in R, \eb(x + y) = \eb(x + v(y)) \vee \eb(y + v(x))$$
to the category of Abelian $l$-groups.
\end{pr}

\bp
First notice that in any directed cr-qsring $R$, for all $x, y \in R$,
$x + v(y) \leq x + y$ since $v(y) \leq y$, therefore
$(x + v(y)) \vee (y + v(x)) \neq \infty$, and 
$(x + v(y)) \vee (y + v(x)) \leq x + y$, in particular,
$$\eb((x + v(y)) \vee (y + v(x))) = \eb(x + v(y)) \vee \eb(y + v(x)) \leq
\eb(x + y),$$
where the elements $\eb(x + v(y)) \leq \eb(x)$ and $\eb(y + v(x)) \leq \eb(y)$ 
are orthogonal, i.e. $\eb(x + v(y)) \bullet \eb(y + v(x)) = \eps$.

Assuming that the inequality above becomes an identity 
for all $x, y \in R$, the following hold.
\ben
\item[\rm (i)] $\eb(x + x) = \eb(x + v(x)) = \eb(v(x)) = \eps$, so 
$2 x := x + x \in \Ep$, and hence $- x = x$ for all $x \in R$.
\item[\rm (ii)] For all $x, y \in R$, $x + x \wedge y = (x + x) 
\wedge (x + y) = \epl(x) \wedge (x + y) = v(x + y) = \epl(x \wedge y)$,
and $x \wedge (x + y) = x + v(x + x + y) = x + \epl(x) + v(y) =
x + v(y)$.
\item[\rm (iii)] The map $\psi : R \lra \Lam \times \Lam_+, x \mapsto 
(v(x), \epl(\eb(x)))$ is injective, identifying the semilattice
$(R, \wedge)$ with a subsemilattice of the semilattice with
support $\Lam \times \Lam_+$ defined by
$$(\gam, \del) \wedge (\gam', \del') := (\gam + \gam', 
\dfrac{\gam \bullet \del + \gam' \bullet \del'}
{(\gam + \gam' \bullet \del') \vee (\gam' + \gam \bullet \del)}),$$
with the associated partial order 
$$(\gam, \del) \leq (\gam', \del') \Llra \gam = \gam' + \gam \bullet \del, 
\gam \bullet \del \leq \gam' \bullet \del' \Llra \gam = \gam' + \gam \bullet \del,
\del \leq \del'.$$
Notice that $\psi(\Ep) = \{(\gam, \eps)\,\mid\,\gam \in \Lam\}$, and
$\psi(\Eb)$ is identified with a convex subset of $\Lam_+$.
\een

Assuming in addition that $R$ is superrigid, we have to show
that $\psi$ maps isomorphically
$R$ onto the directed cr-qsring $\cR_\Lam$ with support 
$\{(\gam, \del) \in \Lam \times \Lam_+\,\mid\,
\mid \gam \mid + \del = \eps\}$, and the operations
$$(\gam, \del) + (\gam', \del') := (\dfrac{\gam \bullet \del + \gam' \bullet
\del'}{\mid \frac{\gam}{\gam'} \mid + (\frac{\gam}{\gam'})_+ \bullet \del +
(\frac{\gam}{\gam'})_- \bullet \del'}\,,\mid \frac{\gam}{\gam'} \mid + 
(\frac{\gam}{\gam'})_+ \bullet \del + (\frac{\gam}{\gam'})_- \bullet \del'\,),$$
$$(\gam, \del) \bullet (\gam', \del') := (\gam \bullet \gam', \del + \del'),
- (\gam, \del) := (\gam, \del), (\gam, \del)^{- 1} := (\gam^{- 1}, \del),
\epsilon := (\eps, \eps).$$
One checks easily that the operations as defined above make $\cR_\Lam$
a superrigid directed cr-qsring, generated by the union
$\Ep(\cR_\Lam) \cup \Eb(\cR_\Lam) = \{(\gam, \eps)\,\mid\,\gam \in \Lam\}
\cup \{(\eps, \del)\,\mid\,\del \in \Lam_+\}$. More precisely, for
any pair $(\gam, \del) \in \cR_\Lam$, we obtain $(\gam_+, \del),
(\gam_-, \del) \in \cR_\Lam$, and
$$(\gam, \del) = (\gam_+, \del) \bullet (\gam_-, \del)^{- 1} = 
((\dfrac{\del}{\gam_+}, \eps) + (\eps, \gam_+ \bullet \del))^{- 1}
\bullet ((\dfrac{\del}{\gam_-}, \eps) + (\eps, \gam_- \bullet \del)).$$ 

To get the desired isomorphism $R \cong \cR_\Lam$, it remains to show
that for all $x \in R$, the elements $\gam := v(x)$ and $\del := \epl(\eb(x))$
of $\Lam$ are orthogonal, i.e. $\mid \gam \mid + \del = \eps$. Since
$R$ is superrigid by assumption, there exists uniquely $y \in \Eb$ such that
$\psi(y) = (\eps, \gam_+ \bullet \del)$. It follows that 
$\psi(x \wedge y) = (\gam_-^{- 1}, (\dfrac{\del}{\gam_-})_+)$, and
$\psi(x + x \wedge y) = ((\dfrac{\del}{\mid \gam \mid})_+ \bullet \gam_-^{- 1},
(\dfrac{\mid \gam \mid + \del}{\gam_-})_+)$. Since $x + x \wedge y \in
\Ep$ by the property (ii) above, it follows that 
$(\dfrac{\mid \gam \mid + \del}{\gam_-})_+ = \eps$, i.e. 
$\mid \gam \mid + \del \leq \gam_-$. Applying the same argument to
the element $x^{- 1} = \psi^{- 1}((\gam^{- 1}, \del))$, we get
$\mid \gam \mid + \del \leq \gam_+$, and hence $\mid \gam \mid + \del
= \eps$, as required. To end the proof, it remains to note that
any morphism $\varphi :\Lam \lra \Lam'$ of Abelian $l$-groups
extends uniquely to a morphism $\wt{\varphi} : \cR_\Lam \lra \cR_{\Lam'},
(\gam, \del) \mapsto (\varphi(\gam), \varphi(\del))$ of directed
cr-qsrings.
\ep 

\begin{rem} \em
Using the superrigid directed cr-qsrings $\cR_\Lam$ as defined
in Proposition 4.8., we can provide examples of directed cr-qsrings 
$R$ with elements $x, y \in R$ satisfying $x - y \in \Ep$ but 
$x \vee y = \infty$ (see Remarks 3.7.(1)). Let 
$\alpha \in \Lam_+$. Then the subset
$\{(\gam, \del) \in \cR_\Lam\,\mid\,\del \leq \alpha\}$ is 
obviously a substructure of $\cR_\Lam$. Now, assuming that
$\alpha \neq \eps$, notice that the following are equivalent.
\ben
\item[\rm (1)] The subset $S_\alpha := \{(\gam, \del) \in \cR_\Lam
\,\mid\,\del < \alpha\}$ is a substructure of $\cR_\Lam$.
\item[\rm (2)] $\alpha$ is not a proper disjoint join in
$\Lam_+$, i.e. $\del + \del' \neq \eps$ whenever $\alpha =
\del \vee \del'$ with $\del, \del' \in \Lam_+ \sm \{\eps\}$;
equivalently, $\{\gam \in \Lam\,\mid\,\mid \gam \mid = \alpha\} =
\{\alpha, \alpha^{- 1}\}$.
\een
Indeed, (1) $\Lra$ (2) since, assuming that $\alpha = 
\del \vee \del'$ with $\del, \del' \neq \eps, \del + \del' = \eps$,
it follows that $(\del', \del), (\del, \del') \in S_\alpha$ but
$(\del', \del) + (\del, \del') = (\eps, \alpha) \not \in S_\alpha$. 
Conversely, (2) $\Lra$ (1) since, assuming that $S_\alpha$ is not a
substructure of $\cR_\Lam$, it follows that there exist
$(\gam, \del), (\gam', \del') \in S_\alpha$ such that 
$(\gam, \del) + (\gam', \del') \not \in S_\alpha$, i.e. 
$\alpha = \rho \vee \rho'$ is a proper disjoint join, where
$\rho = (\dfrac{\gam}{\gam'})_+ \bullet \del + (\dfrac{\gam}{\gam'})_-
\leq \del < \alpha$ and $\rho' = (\dfrac{\gam}{\gam'})_- \bullet
\del' + (\dfrac{\gam}{\gam'})_+ \leq \del' < \alpha$.

Assuming that $\alpha$ is not a proper disjoint join, 
so $R := S_\alpha$ is a directed cr-qsring as 
a substructure of $\cR_\Lam$, and 
assuming in addition that $\alpha = \sigma \vee \tau$
with $\eps < \sigma, \tau < \alpha$, it follows that 
$x := (\eps, \sigma), y := (\eps, \tau) \in \Eb(R)$,
in particular, 
$x - y = x + y = (\si + \tau, \eps) \in \Ep(R) = \Ep(\cR_\Lam)$,
but the elements $x$ and $y$ are not incident in $R$, though
they are incident in $\cR_\Lam$, with $x \vee y = (\eps, \alpha)$. 
For instance, taking $\Lam$ the free Abelian $l$-group \cite{W}
with two generators $\zeta$ and $\eta$, and $R := S_\alpha$, where 
$\alpha := \,|\,\zeta\,|\, \bullet\,|\, \eta\,|$, the elements 
$x := (\eps, \,|\, \zeta \,|\, \bullet\, \eta_+)$, 
$y := (\eps, \,|\, \zeta \,|\, \bullet\, \eta_-) \in \Eb(R)$ satisfy 
the required conditions. For other examples see 
Remarks 5.3.(6) and 5.6.(5).
\end{rem}


\subsection{Locally linear directed commutative regular quasi-semirings}

$\q$ As an immediate consequence of definitions and
Lemmas 4.6. and 4.3.(4), we obtain various characterizations
for subdirectly irreducible directed cr-qsrings as follows..

\begin{co}
The following are equivalent for a directed cr-qsring $R$.
\ben
\item[\rm (1)] $R$ is subdirectly irreducible.
\item[\rm (2)] The Abelian $l$-group $\Lam$ is totally ordered.
\item[\rm (3)] The poset $(R, \leq)$ is an {\em order-tree}, i.e.
there exists the meet $x \wedge y$  for any pair $(x, y) \in R \times R$,
and the poset $x\downarrow$ is totally ordered for all $x \in R$.
\item[\rm (4)] For all $x, y \in R$ satisfying $x - y \in \Ep$, 
either $x \leq y$ or $y \leq x$.
\item[\rm (5)] For all $x, y, z \in R$, either $x \wedge y = 
y \wedge z$ or $x \wedge y = z \wedge x$ or $z \wedge x = y \wedge z$.
In particular, the {\em median} $m(x, y, z) := (x \wedge y) \vee (y \wedge z) 
\vee (z \wedge x)$ of the triple $(x, y, z)$ is well defined, and
moreover $m(x, y, z) \in \{x \wedge y, y \wedge z, z \wedge x\}$ is
the single element of the intersection $[x, y] \cap [y, z] \cap [z, x]$.
\item[\rm (6)] For all $x, y \in R$, $[x, y] = [x \wedge y, x] \cup
[x \wedge y, y]$.
\item[\rm (7)] For all $x, y \in R$, $[x, y] = [x, z] \cup [z, y]$
provided $z \in [x, y]$.
\item[\rm (8)] Any cell of $R$ has at most two ends.
\een
\end{co}

\begin{de}
A directed cr-qsring $R$ is called {\em locally linear} 
(abbreviated an {\em lcr-qsring}) if $R$
satisfies the equivalent conditions of {\em Corollary 4.10.}
\end{de}

By Corollary 4.10., any lcr-qsring has
an underlying structure of $\Lam$-tree \cite{MS}, \cite{AB}.
A natural extension of this arboreal structure to a larger
class of directed cr-qsrings will be discussed in Section 5.

\begin{rems} \em
(1) The totally ordered Abelian groups are 
identified with those 
lcr-qsrings $R$ for which $\Eb = \{\eps\}$. The simplest
lcr-qsring $R$ for which $\Eb$ is not a singleton is 
obtained by adding to the totally ordered multiplicative
group $\Lam = \Ep = \{\gam^n\,\mid\,n \in \Z\} \cong \Z$,
with $\gam > \eps := \gam^0$ the smallest positive element, 
a single element $\del$ satisfying
$\del = - \del = \del^2 = \eb(\del), \del + \del = \epl(\delta) =
\gam, \del + \gam^n = \del$ for $n \geq 1$ and $\gam^n$
otherwise, and $\delta \bullet \gam^n = \gam^n$, in particular,
$v(\del) = \eps$. Thus $R = \Ep \cup \{\del\}$, and 
$\Eb = \{\eps, \del\}$. More generally, the lcr-qsrings
$R$, with $\Lam = \Ep \cong \Z$ as above and $\,|\, \Eb \,|\, = 2$ are
up to isomorphism classified by the pairs $(F, l)$ consisting 
of a field $F$ and a natural number $l$, the simplest case
above corresponding to the pair $(\F_2, 0)$, where $\F_2$ is
the field with two elements. Notice that all lcr-qsrings of this type
are rigid but not superrigid. For any pair $(F, l)$ as above,
we define a lcr-qsring $L = L(F, l)$ as follows. For $l = 0$,
put $L = F^\ast \sqcup \Lam = F \sqcup (\Lam \sm \{\gam\})$, 
identifying the neutral element $0 \in F$ with the generator
$\gam$ of $\Lam$, and extend the operations on
the field $F$ and the ordered multiplicative group $\Lam$ by
setting for all $x \in F^\ast, n \in \Z$, $x + \gam^n = 
x$ if $n \geq 1$, $x + \gam^n = \gam^n$ if $n \leq 0$, and
$x \bullet \gam^n = \gam^n$. For $l \geq 1$, put $L = G \sqcup
\Lam$, where $G := \Z \times F^\ast$ has the natural structure
of Abelian group under the law $(n, x) \bullet (m, y) = (n + m, x y)$,
with the neutral element $\del := (0, 1)$, and $(n, x)^{- 1} =
(- n, x^{- 1})$. In addition, define $- (n, x) = (n, - x)$, 
\begin{center}
$(n, x) + (m, y) = \left\{ \begin{array}{lcl}
(n, x + y) & \mbox{if} & n = m, x + y \neq 0  \\
\gam^{l n + 1} & \mbox{if} & n = m, x + y = 0 \\
(n, x) & \mbox{if} & n < m,
\end{array}
\right.$
\end{center}

\begin{center}
 $(n, x) + \gam^m = \left\{ \begin{array}{lcl}
(n, x) & \mbox{if} & m > l n \\
\gam^m & \mbox{if} & m \leq l n,
\end{array}
\right.$
\end{center}
and $(n, x) \bullet \gam^m = \gam^{l n + m}$. In particular,
$\epl((n, x)) = \gam^{l n + 1}, v((n, x)) = \gam^{l n}$,
$(n, x) \wedge \gam^m = \gam^{\min(n l, m)}, d((n, x), \gam^m) =
\gam^{\mid n l - m \mid + 1}$, and $(n, x) \wedge (m, y) = 
\gam^{\min (n, m) l}, d((n, x), (m, y)) = \gam^{\mid n - m \mid l + 2}$
provided $(n, x) \neq (m, y)$.

Conversely, let $R$ be a lcr-qsring with $\Lam = \Ep \cong \Z$ and
$\Eb = \{\eps, \del\}$, so $R = \Ep \sqcup G$, 
where $G := \{x \in R\,\mid\,\eb(x) = \del\}$.
Thus $(G, \bullet, ^{- 1})$ is an Abelian group with the
neutral element $\del$, and the qvaluation $v : R \lra \Lam = \Ep$
induces a morphism $\varphi : G \lra \Z$ such that
$v(x) = \gam^{\varphi(x)}$ for all $x \in G$. Put 
$G_0 := \Ker(\varphi) = \{x \in G\,\mid\,v(x) = \eps\}$,
and let $l \in \N$ be such that $\varphi(G) = l \Z$,
in particular, $G = G_0$ for $l = 0$. It follows easily that
$x \in G_0 \Lra - x \in G_0$, $x + y \in G_0$ whenever
$x, y \in G_0, x \neq - y$, and $\epl(x) = \epl(\del) =
\gam$ for all $x \in G_0$, and hence $F := G_0 \sqcup \{\gam\}$
becomes a field with $\gam$ as the neutral element for addition,
and $\del$ as the neutral element for multiplication.
If $l = 0$ then $R$ is canonically isomorphic with $L(F, 0)$ as
defined above. If $l \geq 1$, for any $\pi \in G$ satisfying
$\varphi(\pi) = l$, the map $G \lra \Z \times F^\ast, g \mapsto
(\frac{\varphi(g)}{l}, \pi^{- \frac{\varphi(g)}{l}} \bullet g)$
extends to an isomorphism $R \cong L(F, l)$ which is the identity
on $\Lam$ and sends the {\em local uniformizer} $\pi$ to the pair 
$(1, \del)$.

(2) The simplest superrigid lcr-qsrings are of the
form $\cR_\Lam$ (see Proposition 4.8.), where $\Lam$ is a
totally ordered Abelian group, so $\cR_\Lam = \Ep \cup \Eb =
\Lam \times \{\eps\} \cup \{\eps\} \times \Lam_+$.
\end{rems}
\smallskip

According to a classical result of Garrett Birkhoff, we obtain

\begin{pr}
\ben
\item[\rm (1)] Every directed cr-qsring has a subdirect representation
for which all the factors are lcr-qsrings.
\item[\rm (2)] The variety of all directed cr-qsrings is generated
by the class of all lcr-qsrings.
\een
\end{pr}

Any directed cr-qsring $(R, +, \bullet, -, ^{- 1}, \eps)$ 
is canonically embedded as a subalgebra of signature 
(2, 2, 1, 1, 0) into the direct product $\wt{R} := \prod_{P} R_{P}$
of lcr-qsrings, where $P$ ranges over 
the minimal prime convex submonoids of $\Lam_{+}$, in 1-1 
correspondence with the ultrafilters of the distributive lattice 
$(\Lam_{+}, +, \vee)$ with the least element $\eps$, as well as with
the minimal prime convex $l$-subgroups of $\Lam$, and
$R_P := R/\equiv_P$ is the quotient of $R$ by the congruence 
$\equiv_P$ associated to $P$ as in Lemma 4.6. The corresponding 
Abelian $l$-group of the product $\wt{R}$ is the direct product 
$\wt{\Lam} := \prod_P \Lam_P$ of the maximal totally ordered 
factors of $\Lam$, in which the Abelian $l$-group $\Lam$ is 
canonically embedded as a subdirect product.

A basic consequence of Proposition 4.13. is provided by the
statement (2) of the next result.

\begin{co}
Let $R$ be a directed cr-qsring.
\ben
\item[\rm (1)] For all $x, y, z \in R$, $\frac{d(x, z) \bullet d(y, z)}
{d(x, y)} = (x, y)_z^2$, where 
$$(x, y)_z := \frac{\wh{\lam}(x, y, z)}{d(x, y)} =
\frac{\epl(z) \bullet \epl(x \wedge y)}
{\epl(x \wedge z) \bullet \epl(y \wedge z)} \in \Lam_+.$$
The following hold.
\ben
\item[\rm (i)] $(x, y)_z = \eps \Llra z \in [x, y]$; 
\item[\rm (ii)] $(x, y)_z \geq d(x, z) \Llra (x, y)_z = d(x, z) 
\Llra x \in [y, z]$.
\een
\item[\rm (2)] For all $x, y, z, u \in R$, $(x, y)_u + (x, z)_u \leq
(y, z)_u$, i.e. the ternary map \\
$R^3 \lra \Lam_+, (x, y, z) \mapsto (x, y)_u + (x, z)_u$ is symmetric.
\item[\rm (3)] For all $x, y \in R$, the map $d(x, -) : [x, y] \lra
[\eps, d(x, y)], z \mapsto d(x, z)$ is a \\ 
$\Lam$-isometry (not necessarily surjective), 
i.e. $d(z, u) =\, \mid \frac{d(x, z)}{d(x, u)} \mid$ for 
$z, u \in [x, y]$. In particular, the map $d(x, -) : [x, y] \lra [\eps, d(x, y)]$
is a monomorphism of posets, where the partial order 
$\prec_x$ on $[x, y]$ is defined by $z \prec_x u \Llra z \in [x, u]$.
\item[\rm (4)] For all $x, y \in R$, the map 
$[x, y] \lra [x, x \wedge y] \times [x \wedge y, y], z \mapsto
(z\wedge x, z \wedge y)$ is a monomorphism of posets with 
respect to the order $\prec_x$ on $[x, y]$ and the order
$(a, b) \prec (a', b') \Llra a' \leq a, b \leq b'$ on the 
product $[x, x\wedge y] \times [x \wedge y, y]$. Its image
consists of those pairs $(a, b)$ for which $a \vee b \neq \infty$.
The composition of the injective monotone map above
with the isomorphism of bounded distributive lattices
$$[x, x\wedge y] \times [x \wedge y, y] \lra [\eps, \lam(x, y)] \times
[\lam(x, y), d(x, y)], (a, b) \mapsto (d(x, a), d(x, b))$$ 
equals the composition of the injective monotone map $d(x, -)$ from $(3)$
with the monomorphism of bounded distributive lattices
$$[\eps, d(x, y)] \lra [\eps, \lam(x, y)] \times [\lam(x, y), d(x, y)],
\gam \mapsto (\gam + \lam(x, y), \gam \vee \lam(x, y)).$$
\item[\rm (5)] If $[x, y] = [z, u]$ then $x \wedge y = z \wedge u,
\lam(x, z) = \lam(u, y), \lam(x, u) = \lam(z, y), d(x, z) = d(y, u),
d(x, u) = d(y, z)$, and $d(x, y) = d(z, u)$. In particular,
the {\em diameter} $d(x, y)$ is an invariant of the cell $C = [x, y]$.
\item[\rm (6] For all $x, y, z \in R$, the intersection
$[x, y] \cap [y, z] \cap [z, x]$ has at most one element.
\een
\end{co}

\bp
The assertion (1) is immediate by Lemma 4.1.(3).

(2) By (1) we have to show that 
$$\cL := \epl(x \wedge y) \bullet \epl(z \wedge u) + \epl(x \wedge z)
\bullet \epl(y \wedge u) \leq \cR := \epl(y \wedge z) \bullet \epl(x \wedge u)$$
for all $x, y, z, u \in R$. By Proposition 4.13., it suffices 
to check that the inequality $\cL \leq \cR$ holds in the product $\wt{R}$, i.e.
it is true in every factor $R_P$. Consequently, we may assume
from the beginning that $R$ is a lcr-qsring.
Setting $a := x \wedge y \wedge z$, we distinguish by 
Corollary 4.10.(5) the following three cases.

(i) $a = x \wedge y = x \wedge z$ : 
By Lemmas 3.4.(4, 5) and 3.6.(1, 2), we get 
$$\cL = \epl(a) \bullet (\epl(z \wedge u) + \epl(y \wedge u)) =
\epl(a) \bullet \epl(y \wedge z \wedge u) = \epl(a \wedge u) \bullet
\epl(a \vee (y \wedge z \wedge u)) \leq \cR,$$
as required.

(ii) $a = x \wedge y = y \wedge z$ : Then, again by the lemmas above, we obtain
$$A := \epl(x \wedge y) \bullet \epl(z \wedge u) = \epl(a) \bullet
\epl(z \wedge u) = \epl(a \wedge u) \bullet \epl(a \vee (z \wedge u))$$
and 
$$\cR = \epl(a) \bullet \epl(x \wedge u) = \epl(a \wedge u) \bullet
\epl(a \vee (x \wedge u)),$$
therefore $\cL \leq A \leq \cR$ provided $z \wedge u \leq x \wedge u$.
Otherwise, $x \wedge u = x \wedge z$ by Corollary 4.10.(5), and hence
$A = \epl(a) \bullet \epl(z \wedge u), \cR = \epl(a) \bullet \epl(x \wedge u)$, 
and $B := \epl(x \wedge z) \bullet \epl(y \wedge u) = \epl(a)
\bullet \epl((x \wedge u) \vee (y \wedge u))$, therefore 
$\cL = A + B = \epl(a) \bullet \epl((z \wedge u) \wedge ((x \wedge u) \vee 
((y \wedge u))) = \cR$, as desired. 

(iii) $a = x \wedge z = y \wedge z$ : We proceed as in the case (ii).

(3) is a consequence of (2). Indeed, let $x, y, z, u \in R$ be 
such that $z, u \in [x, y]$, i.e. $(y, z)_x = d(x, z)$ and 
$(y, u)_x = d(x, u)$ by (1),(ii). The inequality 
$d(z, u) \geq \mid \frac{d(x, z)}{d(x, u)} \mid$ is obvious by
the triangle inequality. On the other hand, it follows by (2) that
$$d(x, z) + d(x, u) = (y, z)_x + (y, u)_x \leq (z, u)_x,$$
therefore $d(z, u) \leq \frac{d(x, z) \bullet d(x, u)}{(d(x, z) + d(x, u))^2} =
\mid \frac{d(x, z}{d(x, u)} \mid$ as desired.

The rest of the statements follow by straightforward computation.
\ep


\section{Median and locally faithfully full directed 
commutative regular quasi-semirings}

$\q$ We are now ready to introduce an interesting class of directed
cr-qsrings containing all locally linear directed cr-qsrings,
as well as arbitrary direct products of them.

\begin{lem}
The following are equivalent for a directed cr-qsring $R$.
\ben
\item[\rm (1)] For all $x, y, z \in R$, the elements
$x \wedge y, y \wedge z$, and $z \wedge x$ are bounded above,
and hence the element $m(x, y, z) := (x \wedge y) \vee (y \wedge z) 
\vee (z \wedge x)$ is well defined.
\item[\rm (2)] Every finite family $(x_i)_{1 \leq i \leq n}, n \geq 1$, 
of pairwise incident elements of $R$ is bounded above, and hence the join
$\bigvee \{x_i\,\mid\, i = 1, \cdots, n\}$ is well defined.
\item[\rm (3)] For all $x, y, z \in R$, the intersection 
$[x, y] \cap [y, z] \cap [z, x]$ consists of a single element.
\item[\rm (4)] For all $x, y, z \in R$, the intersection 
$[x, y] \cap [y, z] \cap [z, x]$ is nonempty.
\item[\rm (5)] For all $x, y, z \in R$, there exists $u \in R$
such that $[x, y] \cap [x, z] = [x, u]$ and $u \in [y, z]$.
\een
\end{lem}

\bp
The implications (2) $\Lra$ (1), (3) $\Lra$ (4), (4) $\Lra$ (1),
and (5) $\Lra$ (4) are obvious. (1) $\Lra$ (2) follows by induction,
(4) $\Lra$ (3) by Corollary 4.12.(6), and (1) $\Lra$ (5) with 
$u = m(x, y, z)$, by straightforward verification.
\ep

\begin{de}
 A directed cr-qsring $R$ satisfying the equivalent conditions from
{\em Lemma 5.1.} is said to be {\em median}.
\end{de}

\begin{rems} \em 
(1) Assuming that $R$ is a median directed cr-qsring,
it follows by Lemmas 4.3.(2) and 5.1.(1, 3) that $(R, m : R^3 \lra R)$
is a {\em median set (algebra)} \cite{BH}, \cite[Proposition 3.1.]{MG1}.
Moreover, since $m(x, y, z) \wedge u = m(x \wedge u, y \wedge u, z)$ 
for all $x, y, z, u \in R$, the algebra $(R, m, \wedge)$ of 
signature (3, 2) is a {\em directed median set} \cite{DF}. By 
Proposition 3.1.(5) and Lemmas 3.4.(3) and 3.6.(3),
for any $a \in R$, the translations $x \mapsto a + x$ and $x \mapsto a \bullet x$
are endomorphisms of the directed median set $(R, m, \wedge)$.
By Lemma 4.4., Remark 4.5. and Lemma 5.1.(3), the map 
$x \mapsto - x$ is an involutive automorphism of $(R, m, \wedge)$, 
while the map $x \mapsto x^{- 1}$ is an involutive automorphism of the
median set $(R, m)$, i.e. $m(x, y, z)^{- 1} = m(x^{- 1}, y^{- 1}, z^{- 1})$
for all $x, y, z \in R$, and also an isomorphism between the directed
median sets $(R, m, \wedge)$ and $(R, m, \sqcap)$.
 
(2) With the same assumption as above, for all $x, y \in R$,
the poset $([x, y], \prec_x)$, as defined in Corollary 4.14.(3), is
a bounded distributive lattice with the meet operation $(z, u) \in [x, y]
\mapsto m(z, x, u)$, and the join operation $(z, u) \mapsto m(z, y, u)$.
The $\Lam$-isometry $d(x, -) : [x, y] \lra [\eps, d(x, y)], z \mapsto
d(x, z)$ is a monomorphism of bounded distributive lattices.

(3) In a median directed cr-qsring $R$, the element
$(x, y)_z \in \Lam_+$ as defined in Corollary 4.14.(1) is exactly
the distance $d(m(x, y, z), z)$, and hence the inequality from 
Corollary 4.14.(2) can be alternatively proved as follows.
Setting $a := m(x, y, u), b := m(x, z, u), c := m(y, z, u)$, we obtain 
$$[u, a] \cap [u, b] = [u, m(a, b, u)] = [u, m(c, x, u)] = [u, x] \cap [u, c].$$
Since $a, b \in [u, x], m(a, b, u) \in [u, c]$, and the map 
$d(u, -) : [u, x] \lra [\eps, d(u, x)]$
is a morphism of bounded distributive lattices, it follows that
$$(x, y)_u + (x, z)_u = d(u, a) \wedge d(u, b) = d(u, m(a, b, u)) \leq
d(u, c) = (y, z)_u,$$
as desired.

(4) Extending the signature $(+, \bullet, -, ^{- 1}, \eps)$
with a ternary operation standing for the median $m$, it follows by
Lemma 5.1.(4) that the class of all median directed cr-qsrings is a 
variety defined by the finitely many equational axioms for 
directed cr-qsrings extended with axioms expressing the fact 
that the ternary map $m$ is symmetric and $m(x, y, z) \in [x, y]$ 
for all $x, y, z$. By Corollary 4.10., its subdirectly irreducible
members are the lcr-qsrings, identified with those
directed median cr-qsrings which satisfy the suplementary very
restrictive universal axiom
$$m(x, y, z) \in \{x \wedge y, y \wedge z, z \wedge x\}.$$ 

(5) The directed median cr-qsrings are the models of
the inductive theory $T_m$ in the first order language $L$ with
signature $(+, \bullet, -, ^{- 1}, \eps)$, obtained from the
universal theory $T$ of directed cr-qsrings by adding the 
$\forall \exists$-sentence 
$$\forall x, y, z\, \exists u\,(x \wedge y \leq u) \bigwedge
(y \wedge z \leq u) \bigwedge (z \wedge x \leq u).$$
By Proposition 4.13., $T = (T_m)_{\forall}$, the universal theory of
$T_m$, i.e. the directed cr-qsrings form the class of $L$-substructures
of the models of $T_m$. 

(6) For any Abelian $l$-group $\Lam$, the superrigid 
directed cr-qsring $\cR_\Lam$ as defined in Proposition 4.8.
is median. As shown in Remark 4.9., for any 
$\alpha \in \Lam_+ \sm \{\eps\}$,
the necessary and sufficient condition for the subset 
$S_\alpha := \{(\gam, \del) \in \cR_\Lam\,\mid\,\del < \alpha\}$
to be a substructure of $\cR_\Lam$ in the signature
$(+, \bullet, -, ^{- 1}, \eps)$ is that $\alpha$ is not a
proper disjoint join in $\Lam_+$. The necessary and sufficient 
condition for such a substructure $S_\alpha$ to be also closed
under the median operation is that the cardinality of any
subset of $[\eps, \alpha] \sm \{\eps, \alpha\}$ consisting of
pairwise disjoint elements is at most $2$. Indeed, assuming
that there exist $\eps < \si_i < \alpha, i = 1, 2, 3$, such
that $\si_i + \si_j = \eps$ for $i \neq j$, it follows
that $m(x_1, x_2, x_3) = (\eps, \alpha) \in \cR_\Lam \sm S_\alpha$,
where $x_i := (\eps, \alpha \si_i^{- 1}), i = 1, 2, 3$.
Conversely, assuming that $S_\alpha$ is not median, 
there exist $x_i := (\gam_i, \del_i) \in S_\alpha, i = 1, 2, 3$,
such that $x := m(x_1, x_2, x_3) \in \cR_\Lam \sm S_\alpha$,
therefore 
$$\alpha = \epl(\eb(x)) = \bigvee_{1 \leq i < j \leq 3} \epl(\eb(x_i \wedge x_j)) 
\leq m(\del_1, \del_2, \del_3) \leq \alpha.$$
Consequently, $\del_i \vee \del_j = \alpha$ for $1 \leq i < j \leq 3$,
so $\{\alpha \del_i^{- 1}\,\mid\,i = 1, 2, 3\}$ is a subset of 
$[\eps, \alpha] \sm \{\eps, \alpha\}$ consisting of three
pairwise disjoint elements as desired.

The lexicographic extensions of totally ordered Abelian groups
by Abelian $l$-groups provide examples of such substructures
$S_\alpha$ which are not median. For instance, take $\Lam$ 
the free Abelian (multiplicative)
group with $4$ generators $\si_1, \si_2, \si_3$ and $\alpha$. With
respect to the partial order 
$$\si_1^{n_1} \si_2^{n_2} \si_3^{n_3} \alpha^n \leq
\si_1^{m_1} \si_2^{m_2} \si_3^{m_3} \alpha^m \Llra\,{\rm eiher}\,
n < m\,{\rm or}\,n = m\,{\rm and}\,n_i \leq m_i, i = 1, 2, 3,$$ 
$\Lam$ becomes an Abelian $l$-group with
$$\Lam_+ = \{\si_1^{n_1} \si_2^{n_2} \si_3^{n_3} \alpha^n\,\mid\,
{\rm either}\,n > 0\,{\rm or}\,n = 0\,{\rm and}\,n_i \geq 0, i = 1, 2, 3\}.$$
As $\alpha$ is not a proper disjoint join in $\Lam_+$, $S_\alpha$
is a substructure of $\cR_\Lam$ in the signature $(+, \bullet,
-, ^{- 1}, \eps)$ but $S_\alpha$ is not closed under the
median operation since $m(x_1, x_2, x_3) = (\eps, \alpha) \in
\cR_\Lam \sm S_\alpha$, where $x_i := (\eps, \alpha \si_i^{-1}) \in
S_\alpha, i = 1, 2, 3$.
\end{rems}

Another interesting class of directed cr-qsrings, which turns out to be
a proper subclass of the class of all directed median cr-qsrings, is
introduced by the next statement.

\begin{lem}
 The following are equivalent for a directed cr-qsring $R$.
\ben
\item[\rm (1)] For all $x, y \in R$, $x \vee y \neq \infty \Llra
x - y \in \Ep$.
\item[\rm (2)] For all $x, y \in R$, there exists a (unique)
element $z \in R$ such that $d(x, z) = d(y, x \wedge y)$ and
$d(y, z) = d(x, x \wedge y)$, in particular, $z \in [x, y]$.
\item[\rm (3)] The $\Lam$-metric space $(R, d)$ is {\em locally full},
i.e. for all $x, y \in R$, the map \\
$d(x, -) : [x, y] \lra [\eps, d(x, y)], z \mapsto d(x, z)$ is surjective.
\item[\rm (4)] The $\Lam$-metric space $(R, d)$ is 
{\em locally faithfully full}, i.e. for all $x, y \in R$, the
map $d(x, -) : [x, y] \lra [\eps, d(x, y)]$ is bijective.
\een
\end{lem}

\bp
The implications (4) $\Lra$ (3) and (3) $\Lra$ (2) are obvious, 
while (3) $\Lra$ (4) follows by Corollary 4.14.(3).

(1) $\Lra$ (3) : For $x, y \in R$, let $\gam \in [\eps, d(x, y)]$.
We have to find an element $z \in [x, y]$ such that $d(x, z) = 
\gam$. By Lemma 4.3.(4), the elements 
$$z_1 := x + \dfrac{\epl(x)}{\gam + \lam(x, y)} = 
x + \dfrac{\epl(x) \bullet \epl(x \wedge y)}{\epl(x) + (\gam \bullet
\epl(x \wedge y))}$$ 
and 
$$z_2 := y + \dfrac{\epl(y) \bullet (\gam \vee \lam(x, y))}
{d(x, y)} = y + \dfrac{\gam \bullet \epl(x \wedge y)^2}{\epl(x) + (\gam 
\bullet \epl(x \wedge y))}$$ 
are uniquely determined by the conditions 
$z_1 \in [x, x \wedge y], d(x, z_1) =
\gam + d(x, x \wedge y)$ and $z_2 \in [x \wedge y, y],
d(x, z_2) = \gam \vee d(x, x \wedge y)$ respectively. 
As $z_1 - z_2 = x - y + \epl(x \wedge y) = v(x - y) \in \Ep$,
it follows by assumption that $z := z_1 \vee z_2 \neq \infty$.
The element $z$ has the required property thanks to
Corollary 4.14.(4).

(2) $\Lra$ (1) : Let $x, y \in R$ be such that $x - y \in \Ep$,
i.e. $\epl(x \wedge y) = \epl(x) + \epl(y)$.
We have to show that $x \vee y \neq \infty$. By assumption 
there exists $z \in [x, y]$ such that $d(x, z) = d(y, x \wedge y)$.
As $x \wedge y \leq z = (x \wedge z) \vee (y \wedge z)$, it
follows that 
$$d(x, z) := \dfrac{\epl(x) \bullet \epl(z)}{\epl(x \wedge z)^2} =
\dfrac{\epl(x) \bullet (\epl(x \wedge z) \vee \epl(y \wedge z))}
{\epl(x \wedge z)^2} = \dfrac{\epl(x) \bullet \epl(y \wedge z)}
{\epl(x \wedge z) \bullet \epl(x \wedge y)} = $$
$$d(y, x \wedge y) := \dfrac{\epl(y)}{\epl(x \wedge y)},$$
therefore $\dfrac{\epl(y \wedge z)}{\epl(x \wedge z)} =
\dfrac{\epl(y)}{\epl(x)}$. Consequently, 
$$(\dfrac{\epl(y \wedge z)}{\epl(x \wedge z)})_+ 
:= \dfrac{\epl(y \wedge z)}{\epl(x \wedge z) + \epl(y \wedge z)} = 
\dfrac{\epl(y \wedge z)}{\epl(x \wedge y)} =$$
$$(\dfrac{\epl(y)}{\epl(x)})_+ := \dfrac{\epl(y)}{\epl(x) + \epl(y)} = 
\dfrac{\epl(y)}{\epl(x \wedge y)},$$
and hence $\epl(y \wedge z) = \epl(y)$, and similarly, $\epl(x \wedge z) =
\epl(x)$. Thus $x \wedge z = x, y \wedge z = y$, and hence 
$z = x \vee y \neq \infty$ as desired.
\ep

\begin{de}
 A directed cr-qsring $R$ satisfying the equivalent conditions
from {\em Lemma 5.4.} is said to be {\em locally faithfully full}
(abbreviated {\em lff}).
\end{de}

\begin{rems} \em 
(1) Every lff directed cr-qsring $R$ is median. 
Indeed, for arbitrary $x, y, z \in R$, the element 
$(x, y)_z \in \Lam_+$ as defined in Corollary 4.14.(1) is
bounded above by the distance $d(x, z)$, and hence, by
Lemma 5.4.(4), there exists uniquely $u \in [x, z]$ such
that $d(z, u) = (x, y)_z$. According to Remarks 5.3.(3),
we obtain $u = m(x, y, z)$, the median of the triple $(x, y, z)$.

(2) In a lff directed cr-qsring $R$, for all 
$x, y \in R$, the map $d(x, -) : [x, y] \lra [\eps, d(x, y)]$
is an {\em isomorphism} of $\Lam$-metric spaces as well as an
{\em isomorphism} of bounded distributive lattices.

(3) Extending the signature $(+, \bullet, -, ^{- 1}, \eps)$
with a binary operation $\vee'$, the class of all lff directed
cr-qsrings becomes a variety defined by the finitely many
equational axioms for directed cr-qsrings extended with the
following two equational axioms
\ben
\item[\rm (i)] $x \vee' y = y \vee' x$

and
\item[\rm (ii)] $\epl(x) \bullet (x - y) \bullet \epl(x \vee' y) =
\epl(y) \bullet (x - (x \vee' y))^2$,
\een
defining $x \vee' y$ as the unique element $z$ satisfying 
the conditions $d(x, z) = d(y, x \wedge y)$ and $d(y, z) = d(x, x \wedge y)$
from Lemma 5.4.(2). Notice that, as shown in the proof of the
implication (2) $\Lra$ (1) of Lemma 5.4., $x \vee' y = x \vee y
\Llra x - y \in \Ep$.

By Corollary 4.10., the subdirectly irreducible members of
the variety of all lff directed cr-qsrings are the lcr-qsrings,
identified with those lff directed cr-qsrings which satisfy
the supplementary very restrictive universal axiom
$$(x \vee' y \leq x) \bigvee (x \vee' y \leq y).$$

(4) The lff directed cr-qsrings are the models
of the inductive theory $T_{lff}$ in the first order language 
$L$ with signature $(+. \bullet, -, ^{- 1}, \eps)$,
obtained from the universal theory $T$ of directed
cr-qsrings by adding the $\forall \exists$-sentence
$$\forall x, y\,\exists z\,(2 (x - y) \neq x - y)
\bigvee  ((x \leq  z) \bigwedge (y \leq z)).$$
By Propostion 4.13., $T = (T_{lff})_{\forall}$,
the universal theory of $T_{lff}$, i.e. the directed 
cr-qsrings form the class of $L$-substructures of the models 
of $T_{lff}$.

(5) For any Abelian $l$-group $\Lam$, the superrigid 
directed cr-qsring $\cR_\Lam$ as defined in Proposition 4.8.
is lff, in particular median, where 
$(\gam, \del) \vee (\gam', \del') = (\gam \vee \gam',
\del \vee \del')$ for all $(\gam, \del), (\gam', \del') \in \cR_\Lam$
satisfying $(\gam, \del) + (\gam', \del') \in \Ep(\cR_\Lam)$,
i.e. $\gam + \gam' \bullet \del' = \gam' + \gam \bullet \del$.
For all $\alpha \in \Lam_+ \sm \{\eps\}$,
the necessary and sufficient condition for the subset 
$S_\alpha := \{(\gam, \del) \in \cR_\Lam\,\mid\,\del < \alpha\}$
to be lff, i.e. a substructure of $\cR_\Lam$ in the extended signature
$(+, \bullet, \vee', -, ^{- 1}, \eps)$ is that $\si + \tau \neq \eps$
for all $\eps < \si, \tau < \alpha$, in particular, $\alpha$ is not a
proper disjoint join in $\Lam_+$. Indeed, assuming that there
exists a pair of disjoint elements $(\si, \tau)$ in 
$[\eps, \alpha] \sm \{\eps, \alpha\}$, it follows that 
$x \vee y = (\eps, \alpha) \in \cR_\Lam \sm S_\alpha$, where
$x := (\eps, \alpha \si^{- 1}), y = (\eps, \alpha \tau^{- 1}) \in 
S_\alpha$. The converse is obvious. 

As an example of a median substructure $S_\alpha$ which is not lff,
take $\Lam$ the free Abelian (multiplicative) group with three
generators $\si_1, \si_2$ and $\alpha$, becoming an $l$-group
under the partial order
$$\si^{n_1} \si_2^{n_2} \alpha^n \leq \si_1^{m_1} \si_2^{m_2}
\alpha^m \Llra\,{\rm either}\,n < m\,{\rm or}\,n = m\,{\rm and}\,
n_i \leq m_i, i = 1, 2.$$
In this case, $S_\alpha$ is a median substructure of $\cR_\Lam$
according to Remarks 5.3.(6), but $S_\alpha$ is not lff since
$x_1 \vee x_2 = (\eps, \alpha) \in \cR_\Lam \sm S_\alpha$,
where $x_i := (\eps, \alpha \si_i^{- 1}) \in S_\alpha, i = 1, 2$.
\end{rems}

\bigskip


\section{$l$-valuations and Pr\" ufer extensions associated 
to superrigid directed commutative regular quasi-semirings}

$\q$ In this section we associate a $l$-valued commutative ring
to a superrigid directed cr-qsring, and we discuss the case
when the associated commutative ring is a Pr\" ufer extension
of its $l$-valuation subring. 

\subsection{The $l$-valued commutative ring associated to a
superrigid directed commutative regular quasi-semiring}

$\q$ For a given directed cr-qsring $R$, let 
$R_\alpha := \{x \in R\,|\,\epl(x) = \alpha\}$ for any
$\alpha \in \Lam := \Ep(R)$. The family $(R_\alpha, +)_{\alpha \in \Lam}$ 
is an inverse system of Abelian groups with the connecting morphisms 
$\pi_{\alpha\,\beta} : R_\beta \lra R_\alpha, x \mapsto x + \alpha$ 
for $\alpha \leq \beta$. We denote by 
$B := B(R) = \displaystyle\lim_{\stackrel{\lla}{\alpha \in \Lam}} 
R_\alpha \cong \displaystyle\lim_{\stackrel{\lla}{\alpha \in \Lam_+}} R_\alpha$, 
the inverse limit of the inverse system above.
Thus the Abelian group $B$ consists of all maps 
$\varphi : \Lam_+ \lra R$ satisfying $\epl(\varphi(\alpha)) =
\alpha$ for $\alpha \in \Lam_+$, and $\varphi(\alpha) \leq
\varphi(\beta)$ for $\alpha \leq \beta$. The operations 
$+, -$ are defined pointwise, while the 
neutral element $0$ is the embedding 
$\Lam_+ \lra R, \alpha \mapsto \alpha$.
The injective map 
$\varphi \in B \mapsto \wt{\varphi} : \Lam \lra R$,
where $\wt{\varphi}(\alpha) := \varphi(\alpha_+) + \alpha$
for $\alpha \in \Lam$, so $\wt{\varphi}|_{\Lam_+} = \varphi$,
maps $B$ onto the Abelian group consisting of the sections 
$\psi : \Lam \lra R$ of the surjective map 
$\epl : R \lra \Lam, x \mapsto \epl(x)$, satisfying 
$\psi(\alpha \wedge \beta) = \psi(\alpha) \wedge
\psi(\beta), \psi(\alpha \vee \beta) = \psi(\alpha)
\vee \psi(\beta)$ for all $\alpha, \beta \in \Lam$. 

The Abelian group $(B, +, -, 0)$ becomes a commutative ring 
(not necessarily unital) with respect to the multiplication 
$B \times B \lra B,(\varphi, \psi) \mapsto \varphi \cdot \psi$, 
defined by
$$(\varphi \cdot \psi)(\alpha) := \varphi(\alpha \bullet 
v(\psi(\eps))^{- 1}) \bullet \psi(\alpha \bullet 
v(\varphi(\eps))^{- 1}) + \alpha = $$
$$\displaystyle\lim_{\stackrel{\gam \rightarrow \infty}{}} 
(\varphi(\gam) \bullet \psi(\gam) + \alpha),\,{\rm for}\,\alpha \in \Lam_+.$$
The ring $B$ is unital if and only if the directed cr-qsring $R$
is superrigid, and in this case the unit $1 \in B$ is the 
injective map $\Lam_+ \lra R$, sending $\alpha \in \Lam_+$ to 
the unique element $1_\alpha \in \Eb(R)$ satisfying 
$\epl(1_\alpha) = \alpha$. 
We assume in the following that $R$ is superrigid.

For every $\varphi \in B$, the map 
$v \circ \varphi : \Lam_+ \lra \Lam$
belongs to the commutative $l$-monoid $\wh{\Lam}$ as defined in
1.1, and the map $w : B \lra \wh{\Lam}, \varphi \mapsto
w(\varphi) := v \circ \varphi$ is a $l$-valuation on the commutative
ring $B$, with supp$(w) = \{0\}$, so the ring $B$ is reduced, 
and $A := A_w = w^{- 1}(\wh{\Lam}_+) = \{\varphi \in B\,|\,
\varphi(\eps) = \eps\}$, the $l$-valuation ring of $w$.
Call the map $w : B \lra \wh{\Lam}$ the {\em $l$-valuation
associated to the superrigid directed cr-qsring} $R$.

For every $\alpha \in \Lam_+$, the subgroup $T_\alpha :=
\{x \in R_\alpha\,|\,\eps \leq x\} = \{x \in R_\alpha\,|\,v(x) 
\in \Lam_+\}$ of the Abelian group $(R_\alpha, +)$ becomes
a commutative ring with multiplication $x \cdot y := x \bullet y +
\alpha$ for $x, y \in T_\alpha$, and unit $1_\alpha$. For 
$\alpha \leq \beta$, the connecting map $\pi_{\alpha\,\beta} :
R_\beta \lra R_\alpha$ induces by restriction a ring morphism 
$T_\beta \lra T_\alpha$, and 
$A = \displaystyle\lim_{\stackrel{\lla}{\alpha \in \Lam_+}} T_\alpha
\cong  \displaystyle\lim_{\stackrel{\lla}{\alpha \in \Lam_+}} A/I_\alpha$,
where $I_\alpha := \Ker(A \lra T_\alpha, \varphi \mapsto 
\varphi(\alpha)) = \{\varphi \in B\,|\,w(\varphi) \geq \alpha\}$
is an ideal of $A = I_\eps$. On the other hand, for every
$\alpha \in \Lam_+$, $R_\alpha$ becomes an $A$-module
under the action $\varphi \cdot x := \varphi(\alpha \bullet
v(x)_-) \bullet x + \alpha$ for $\varphi \in A, x \in R_\alpha$,
in particular, $\varphi \cdot x = \varphi(\alpha) \cdot x$ for
$x \in T_\alpha$. Thus the canonical map $\pi_\alpha : B \lra R_\alpha,
\varphi \mapsto \varphi(\alpha)$, for $\alpha \in \Lam_+$,
is a morphism of $A$-modules, and 
$B \cong \displaystyle\lim_{\stackrel{\lla}{\alpha \in \Lam_+}} R_\alpha \cong
\displaystyle\lim_{\stackrel{\lla}{\alpha \in \Lam_+}} B/I_\alpha$ as
$A$-modules.

\begin{lem}
Let $R$ be a superrigid directed cr-qsring $R$, $w : B \lra
\wh{\Lam}$ its associated $l$-valuation, and $A$ the
$l$-valuation ring. Then $w^{- 1}(\Lam) = B^\ast$, and 
$w^{- 1}(\eps) = A^\ast$, so $w(B) \cap \Lam \cong B^\ast/A^\ast$.
In particular, $B$ is a field and $A$ is a valuation ring of $B$
whenever $\Lam$ is totally ordered, i.e. $R$ is a superrigid lcr-qsring.
\end{lem}

\bp
We have to show that $w^{- 1}(\Lam) \sse B^\ast$. Let $\varphi \in B$
with $w(\varphi) = \gam \in \Lam$. Define the map 
$\psi : \Lam_+ \lra R$ by $\psi(\alpha) = \varphi(\alpha 
\bullet \gam_+^2)^{- 1} + \alpha$ for $\alpha \in \Lam_+$.
As $\alpha \bullet \gam_+^2 \geq \gam_+$ for all $\alpha
\in \Lam_+$, it follows that 
$$\epl(\psi(\alpha)) = \dfrac{\alpha \bullet \gam_+^2}
{v(\varphi(\alpha \bullet \gam_+^2))^2} + \alpha = \dfrac{\alpha \bullet
\gam_+^2}{\gam^2} + \alpha = \alpha \bullet \gam_-^2 + \alpha = \alpha.$$ 
To conclude that $\psi \in B$, it suffices to show that
$\varphi(\alpha \bullet \gam_+^2)^{- 1} \leq \varphi
(\beta \bullet \gam_+^2)^{- 1}$ for $\alpha \leq \beta$.
By Lemma 3.5.(4), it follows that
$$\varphi(\alpha \bullet \gam_+^2)^{- 1} \wedge
\varphi(\beta \bullet \gam_+^2)^{- 1} = \varphi(\beta \bullet
\gam_+^2)^{- 1} \bullet \eb(\varphi(\alpha \bullet
\gam_+^2)) \leq \varphi(\alpha \bullet \gam_+^2)^{- 1}.$$
Applying $\epl$ to the both members of the last inequality,
we obtain by Lemma 2.6.(4) the same value $\alpha \bullet
\gam_-^2$, so the inequality becomes an equality as
desired. 

Finally note that $w(\psi) = \gam^{- 1}$ and
$$(\varphi \psi)(\alpha) = \varphi(\alpha \bullet \gam_+)
\bullet \psi(\alpha \bullet \gam_-) + \alpha = \varphi(\alpha
\bullet \gam_+) \bullet (\varphi(\alpha \bullet |\gam| 
\bullet \gam_+)^{- 1} + \alpha \bullet \gam_-) + \alpha = 1_\alpha$$
for all $\alpha \in \Lam_+$, and hence $\varphi = \psi^{- 1} \in
B^\ast$, as required. 
\ep

\begin{lem}
With the same data as in {\em Lemma 6.1}, let 
$\Eb(B) := \{\varphi \in B\,|\,\varphi^2 = \varphi\}$ be
the boolean algebra of idempotent elements of the
commutative ring $B$. Then the $l$-valuation 
$w : B \lra \wh{\Lam}$ induces by restriction an anti-isomorphism 
of boolean algebras $\Eb(B) \lra \partial \wh{\Lam}_+$; we denote  
its inverse by $\eta$.
\end{lem}

\bp
The inclusion $w(\Eb(B)) \sse \Eb(\wh{\Lam}) = \partial
\wh{\Lam}_+$ is obvious, in particular $\Eb(B) \sse A$.
First let us show that the map
$w|_{\Eb(B)} : \Eb(B) \lra \partial \wh{\Lam}_+$ is injective.
Let $\varphi \in \Eb(B)$, i.e. $\varphi(\alpha)^2 + \alpha =
\varphi(\alpha)$, in particular
$v(\varphi(\alpha)) \vee \dfrac{\alpha}{v(\varphi(\alpha))} =
\alpha$, for all $\alpha \in \Lam_+$. By multiplication 
with $\varphi(\alpha)^{- 1}$, it follows that
$$\varphi(\alpha) + \dfrac{\alpha}{v(\varphi(\alpha))} \leq
\varphi(\alpha)^{- 1} \bullet (\varphi(\alpha)^2 + \alpha) =
\eb(\varphi(\alpha)) = 1_{\frac{\alpha}{v(\varphi(\alpha))}}.$$
The inequality above becomes an identity since
$\epl(\varphi(\alpha) + \dfrac{\alpha}{v(\varphi(\alpha))}) =
\dfrac{\alpha}{v(\varphi(\alpha))} = \epl(1_{\frac{\alpha}
{v(\varphi(\alpha))}})$, and hence 
$1_{\frac{\alpha}{v(\varphi(\alpha))}} \leq \varphi(\alpha)$.
As $v(\varphi(\alpha)) \leq \varphi(\alpha)$ too, it follows
that \\
$v(\varphi(\alpha)) \vee 1_{\frac{\alpha}{v(\varphi(\alpha))}}
\leq \varphi(\alpha)$. Applying $\epl$ to the both members of
the inequality above, we get the same value $\alpha$, therefore
$\varphi(\alpha) = v(\varphi(\alpha)) \vee 1_{\frac{\alpha}
{v(\varphi(\alpha))}} = w(\varphi)(\alpha) \vee 1_{\frac{\alpha}
{w(\varphi)(\alpha)}}$, so $\varphi \in \Eb(B)$ is completely
determined by $w(\varphi)$, i.e. the map $w|_{\Eb(B)} : \Eb(B) \lra
\partial \wh{\Lam}_+$ is injective. 

On the other hand, we define a map $\eta : \partial
\wh{\Lam}_+ \lra \Eb(B)$ as follows. Let $\theta \in 
\partial \wh{\Lam}_+, \alpha \in \Lam_+$. As 
$\epl(1_{\frac{\alpha}{\theta(\alpha)}}) + \theta(\alpha) = 
\eps$, so $(\frac{\alpha}{\theta(\alpha)^2})_+ = \frac{\alpha}
{\theta(\alpha)}, (\frac{\alpha}{\theta(\alpha)^2})_- = 
\theta(\alpha)$, it follows that 
$\eta(\theta)(\alpha) := 1_{\frac{\alpha}{\theta(\alpha)^2}}^{- 1} =
\theta(\alpha) \vee 1_{\frac{\alpha}{\theta(\alpha)}}$,
where $1_{\frac{\alpha}{\theta(\alpha)^2}} := 1_\frac{\alpha}
{\theta(\alpha)} + \frac{\alpha}{\theta(\alpha)^2}$(cf. 4.2).
One checks easily that the map $\eta(\theta) : \Lam_+ \lra R$
belongs to $\Eb(B)$, so we have obtained a map 
$\eta : \partial \wh{\Lam}_+ \lra \Eb(B)$ such that 
$w \circ \eta = 1_{\partial \wh{\Lam}_+}$. Consequently,
the map $w|_{\Eb(B)} : \Eb(B) \lra \partial \wh{\Lam}_+$ is bijective
and $\eta$ is its inverse. Moreover it is an anti-isomorphism
of boolean algebras since $w(\varphi \cdot \psi) = w(\varphi)
\bullet w(\psi) = w(\varphi) \vee w(\psi)$ for $\varphi, \psi
\in \Eb(B)$, and 
$$w(1 - \varphi)(\alpha) = v(1_\alpha - \varphi(\alpha))
= v(1_\alpha - (1_{\frac{\alpha}{v(\varphi(\alpha))}} \vee 
v(\varphi(\alpha)))) =$$ 
$$\epl(1_\alpha \wedge (1_{\frac{\alpha}{v(\varphi(\alpha))}}
\vee v(\varphi(\alpha)))) = \epl(1_{\frac{\alpha}{v(\varphi(\alpha))}})
= \dfrac{\alpha}{v(\varphi(\alpha))} =
(\neg w(\varphi))(\alpha)$$
for $\varphi \in \Eb(B), \alpha \in \Lam_+$.
\ep

\begin{pr}
Let $R$ be a superrigid directed cr-qsring. With the notation
above, for any $\theta \in \partial \wh{\Lam}_+$, let 
$R^\theta := \{x \in R\,|\,\wt{\theta}(\epl(x)) = \wt{\theta}
(v(x)) = \eps\} = \{x \in R\,|\,d(x, \eps) \in \Ker(\wt{\theta})_+\}$;
in particular, $R^\eps = R, R^\om = \{\eps\}$. Then the following 
assertions hold.
\begin{enumerate}
 \item[\rm (1)] $R^\theta$ is a substructure of $R$ with
$\Ep(R^\theta) = \Ker(\wt{\theta})$, $\Eb(R^\theta) =
\{1_\alpha\,|\,\alpha \in \Ker(\wt{\theta})_+\}$. Call
$R^\theta$ the {\em superrigid directed cr-qsring induced by}
$\theta$, and denote by \\ $B^\theta := B(R^\theta)$ and
$w^\theta : B^\theta \lra \wh{\Ker(\wt{\theta})}$ the
commutative ring and the $l$-valuation associated as above
to $R^\theta$.
\item[\rm (2)] $(w^\theta)^{- 1}(\Ker(\wt{\theta})) = 
(B^\theta)^\ast$ and $(w^\theta)^{- 1}(\eps) = (A^\theta)^\ast$, 
where $A^\theta \sse B^\theta$ denotes the $l$-valuation ring of
$w^\theta$.
\item[\rm (3)] $w^{- 1}(\wh{\Lam} \theta) = B \eta(\theta) :=
\{\varphi \cdot \eta(\theta)\,|\,\varphi \in B\} = 
\{\varphi \in B\,|\,\varphi \cdot \eta(\theta) = \varphi\}$.
\item[\rm (4)] The restriction map $\varphi \mapsto \varphi|_
{\Ker(\wt{\theta})_+}$ yields a ring isomorphism 
$B \eta(\theta) \lra B^\theta$, and $w^\theta(\varphi|_
{\Ker(\wt{\theta})_+}) = w(\varphi)|_{\Ker(\wt{\theta})_+}$
for all $\varphi \in B \eta(\theta)$. In particular, 
$w^{- 1}(\Lam \theta) = (B \eta(\theta))^\ast = B^\ast \eta(\theta)
\cong (B^\theta)^\ast$ and
$w^{- 1}(\theta) = (A \eta(\theta))^\ast = A^\ast \eta(\theta)
\cong (A^\theta)^\ast$.
\end{enumerate}
\end{pr}

\bp
As $d(x, \eps) = \frac{\epl(x)}{\epl(x \wedge \eps)^2} =
|v(x)| \bullet \epl(\eb(x))$ for all $x \in R$, and
$\Ker(\wt{\theta})$ is a convex $l$-subgroup of $\Lam$,
it follows that $d(x, \eps) \in \Ker(\wt{\theta})_+ \Llra
\wt{\theta}(\epl(x)) = \wt{\theta}(v(x)) = \eps$.

(1) We have to show that $R^\theta$ is closed under the
operations $+, -, \bullet, ^{- 1}$. As $d(- x, \eps) =
d(x^{- 1}, \eps) = d(x, \eps)$ for all $x \in R$ by
Lemma 4.4.(2), it follows that $R^\theta$ is closed under
the unary operations $-$ and $^{- 1}$. Let $x, y \in R^\theta$.
Applying the endomorphism $\wt{\theta}$ to the relations
$\epl(x + y) = \epl(x) \wedge \epl(y)$ and $v(x) \wedge v(y) \leq 
v(x + y) \leq \epl(x + y)$, we obtain $x + y \in R^\theta$,
while $x \bullet y \in R^\theta$ follows by applying $\wt{\theta}$
to the identities $v(x \bullet y) = v(x) \bullet v(y)$ and
$\epl(x \bullet y) = \epl(x) \bullet v(y) \wedge v(x) \bullet \epl(y)$.
Since $d(\alpha, \eps) = |\alpha|$ for $\alpha \in \Lam = \Ep$,
and $d(1_\alpha, \eps) = \alpha$ for $\alpha \in \Lam_+$, it
follows that $\Ep(R^\theta) = \Ker(\wt{\theta})$ and
$\Eb(R^\theta) = \{1_\alpha\,|\,\alpha \in \Ker(\wt{\theta})_+)$,
and hence the directed cr-qsring $R^\theta$ is superrigid.

(2) Apply Lemma 6.1 to the superrigid directed cr-qsring $R^\theta$.

(3) Let $\varphi \in B$. As supp$(w) = \{0\}$, it follows
by Lemma 1.10.(1), (i) $\Llra$ (ii), and Lemma 6.2 that
$$\varphi \in w^{- 1}(\wh{\Lam} \theta) \Llra  w(\varphi \cdot
(1 - \eta(\theta))) = w(\varphi) \cdot \neg \theta = \om
\Llra \varphi \cdot (1 - \eta(\theta)) = 0 \Llra \varphi \in 
B \eta(\theta).$$

(4) Let $\varphi \in B \eta(\theta), \alpha \in \Ker(\wt{\theta})_+$.
Since $\epl(\varphi(\alpha)) = \alpha \in \Ker(\wt{\theta})$
and $v(\varphi(\alpha)) = w(\varphi)(\alpha) \in \Ker(\wt{\theta})$
by Lemma 1.10.(1), (i) $\Llra$ (iii), it follows that $\varphi(\alpha)
\in R^\theta$, therefore $\varphi|_{\Ker(\wt{\theta})_+} \in B^\theta$.
The map $B \eta(\theta) \lra B^\theta, \varphi \mapsto \varphi|_
{\Ker(\wt{\theta})_+}$ is obviously a ring morphism. The morphism
above is injective. Indeed, assuming that 
$\varphi|_{\Ker(\wt{\theta})_+} = 0$, i.e. $\varphi(\alpha)
= \alpha$ for all $\alpha \in \Ker(\wt{\theta})_+$, it follows
by Lemma 1.10.(1), (i) $\Llra$ (vi), that 
$$w(\varphi)(\alpha) = \theta(\alpha) \cdot 
w(\varphi)(\frac{\alpha}{\theta(\alpha)}) = \theta(\alpha) \cdot 
v(\varphi(\frac{\alpha}{\theta(\alpha)}))  = \theta(\alpha) \cdot 
\frac{\alpha}{\theta(\alpha)} = \alpha$$
for all $\alpha \in \Lam_+$, so $w(\varphi) = \om$ and
hence $\varphi = 0$ as required. To show that the morphism
above is surjective, let $\psi \in B^\theta$, and consider
for any $\alpha \in \Lam_+$ the pair 
$$(\psi(\frac{\alpha}{\theta(\alpha)}) \in R^\theta \sse R, 
\theta(\alpha) \bullet v(\psi(\eps)) \in \Lam = \Ep).$$
Since 
$$\epl(\psi(\frac{\alpha}{\theta(\alpha)})) + \theta(\alpha)
\bullet v(\psi(\eps)) = \frac{\alpha}{\theta(\alpha)} +
\theta(\alpha) \bullet v(\psi(\eps)) = v(\psi(\eps)) \leq 
v(\psi(\frac{\alpha}{\theta(\alpha)})),$$
it follows by Lemma 3.6.(7) that the element $\varphi(\alpha) :=
\psi(\frac{\alpha}{\theta(\alpha)}) \vee \theta(\alpha) \bullet
v(\psi(\eps)) \in R$ is well defined. As $\epl(\varphi(\alpha))
= \frac{\alpha}{\theta(\alpha)} \vee \theta(\alpha) \bullet
v(\psi(\eps)) = \alpha$ and $\alpha \leq \beta \Lra \varphi(\alpha)
\leq \varphi(\beta)$, we deduce that the map $\varphi : \Lam_+ \lra
R$ is an element of the ring $B$ such that 
$\varphi|_{\Ker(\wt{\theta})_+} = \psi$. Moreover it follows
by (3) and Lemma 1.10.(1), (i) $\Llra$ (v), that 
$\varphi \in B \eta(\theta)$ as desired. 
\ep

\begin{rem} \em 
 The correspondence $R \mapsto (B, w : B \lra \wh{\Lam})$
extends to a covariant functor defined on the category of
superrigid directed cr-qsrings with morphisms \\ $f: R \lra R'$
satisfying the condition that the Abelian $l$-group 
$\Lam' := \Ep(R')$ is the convex hull of the image 
$f(\Lam) \sse \Lam'$. For such a morphism $f : R \lra R'$,
define the ring morphism $B(f) : B(R) \lra B(R')$ by
$B(f)(\varphi)(\alpha') :=$ 
$$\displaystyle\lim_{\stackrel
{\alpha \rightarrow \infty}{}} (f(\varphi(\alpha)) + \alpha') =
f(\varphi(\alpha)) + \alpha'\,{\rm for\,\,some\,(for\,\,all)}\,
\alpha \in \Lam_+\,{\rm such\,that}\,f(\alpha) \geq \alpha',$$
for $\varphi \in B(R), \alpha' \in \Lam'_+$. On the other
hand, let $\wh{f} : \wh{\Lam} \lra \wh{\Lam'}$ be the morphism
of commutative $l$-monoids induced by the morphism of
Abelian $l$-groups $f|_\Lam : \Lam \lra \Lam'$ (cf. Lemma 1.2.(1)).
It follows that the pair $(B(f), \wh{f}) : (B(R), w) \lra (B(R'), w')$ 
is a morphism of $l$-valued commutative rings, i.e. 
$w' \circ  B(f) = \wh{f} \circ w$.
\end{rem}

\begin{rem} \em
 There exist nontrivial superrigid directed cr-qsrings $R$ 
with the associated $l$-valuation $w : B \lra \wh{\Lam}$
such that $A_w = B$. For instance, let $R = \cR_\Lam$ 
as defined in the proof of Proposition 4.8, where $\Lam$ is
a nontrivial Abelian $l$-group. It follows that the associated
$l$-valuation $w : B \lra \wh{\Lam}$ is injective with
image $w(B) = \partial \wh{\Lam}_+$, identifying the ring
$B = A_w$ with the underlying boolean ring of the boolean algebra
with support $\partial \wh{\Lam}_+$ and opposite order.
In particular, $B = \F_2$, the field with $2$ elements,
provided $\Lam$ is totally ordered.
\end{rem}

\subsection{Pr\" ufer commutative regular quasi-semirings}

$\q$ Let $R$ be a superrigid directed cr-qsring, and $B := B(R)$, 
$w : B \lra \wh{\Lam}, A := A_w$, its associated commutative ring, 
$l$-valuation and $l$-valuation subring respectively. Let $M := M(A, B)$ be
the commutative $sl$-monoid of all f.g. $A$-submodules of $B$,
associated to the ring extension $A \sse B$ (cf. 1.2). 
Let $\wt{M} := \{I \in M\,|\,I B = B\} = \{I \in M\,|\,
\exists J \in M, A \sse I J\}$ be the monoid of $B$-regular f.g.
$A$-submodules of $B$, and $M^\ast$ be the subgroup
of invertible $A$-submodules of $B$. As shown in 1.3, the
$l$-valuation $w$ induces a canonical morphism of $sl$-monoids
$\wh{w} : M \lra \wh{\Lam}$, defined by $\wh{w}(I) := \bigwedge
_{1 \leq i \leq n} w(\varphi_i)$, where
$\{\varphi_i\,|\,i = 1, \cdots, n\}$ is an arbitrary system of 
generators of the f.g. $A$-submodule $I$ of $B$.
Recall that $M^\ast \sse \wt{M} \sse \wh{w}^{- 1}(\Lam)$,
the morphism $\wh{w}$ induces by restriction a 
monomorphism of ordered groups $M^\ast \lra \Lam$, and
$I_{\wh{w}(\mathbf{a})} := \{\varphi \in B\,|\,w(\varphi) \geq
\wh{w}(\mathbf{a})\} = \mathbf{a}$ for all $\mathbf{a} \in M^\ast$.

The next lemma permits us to characterize those superrigid
directed cr-qsrings $R$ for which the associated ring extension 
$A \sse B$ is Pr\" ufer, i.e. $M^\ast = \wt{M}$, so
$w$ is a Pr\" ufer $l$-valuation. Call them {\em Pr\" ufer
cr-qsrings}.

\begin{lem}
 Let $R$ be a superrigid directed cr-qsring, and let 
$\varphi, \psi \in B := B(R)$. Then the following 
assertions are equivalent.
\begin{enumerate}
 \item[\rm (1)] The elements $\psi, \varphi \psi$ and
$\varphi (1 - \varphi \psi)$ belong to the $l$-valuation
ring $A := A_w$.
\item[\rm (2)] $\psi(\gam^2) = \varphi(\gam)^{- 1} + \gam^2$,
i.e. $\psi(\gam^2) \leq \varphi(\gam)^{- 1}$,
where $\gam := w(\varphi)_- = v(\varphi(\eps))^{- 1}$.
\end{enumerate}
\end{lem}

\bp
$(1) \Lra (2)$ : By assumption, $\psi(\eps) = (\varphi
\cdot \psi)(\eps) = \eps$, and $\varphi(\eps) = (\varphi^2 \cdot \psi)(\eps)$.
Consequently, we obtain $(\varphi \cdot \psi)(\eps) = \varphi(\eps)
\bullet \psi(\gam) + \eps = \eps$, therefore $\varphi(\eps) \bullet
\psi(\gam) = \eps$ since $\epl(\varphi(\eps) \bullet \psi(\gam)) \leq 
v(\varphi(\eps)) \bullet \gam = \eps$. It follows that 
$v(\psi(\gam)) = v(\varphi(\eps))^{- 1} = \gam = \epl(\psi(\gam))$, 
so $\psi(\gam) = \gam \in \Ep$. Further we obtain 
$(\varphi \cdot \psi)(\gam) = \varphi(\gam) \bullet \psi(\gam^2) +
\gam$, and hence $(\varphi \cdot \psi)(\gam) = \varphi(\gam) \bullet
\psi(\gam^2)$ since $v(\varphi(\gam)) = v(\varphi(\eps)) = \gam^{- 1}$
implies $\epl(\varphi(\gam) \bullet \psi(\gam^2)) \leq 
v(\varphi(\gam)) \bullet \gam^2 = \gam$. Consequently,
$$\varphi(\eps) = (\varphi^2 \cdot \psi)(\eps) = \varphi(\eps) \bullet 
(\varphi \cdot \psi)(\gam) + \eps = \varphi(\eps) \bullet 
\varphi(\gam) \bullet \psi(\gam^2)$$ 
since $\epl(\varphi(\eps) \bullet (\varphi \cdot \psi)(\gam)) 
\leq v(\varphi(\eps)) \bullet \gam = \eps$. As 
$$\epl(\eb(\psi(\gam^2)) = \frac{\gam^2}{v(\psi(\gam^2))} \leq 
\frac{\gam^2}{v(\psi(\gam))} = \gam = \epl(\eb(\varphi(\eps)))
\leq \epl(\eb(\varphi(\gam)),$$
we obtain the desired inequality by multiplying with 
$\varphi(\eps)^{- 1} \bullet \varphi(\gam)^{- 1}$ 
the both terms of the identity above 
$$\varphi(\gam)^{- 1} \geq \eb(\varphi(\eps)) \bullet \varphi(\gam)^{- 1} = 
\eb(\varphi(\eps)) \bullet \eb(\varphi(\gam)) \bullet \psi(\gam^2) =
\psi(\gam^2).$$

$(2) \Lra (1)$ : We obtain $w(\psi) \geq v(\psi(\gam^2)) =
v(\varphi(\gam)^{- 1}) + \gam^2 = \gam + \gam^2 = \gam \geq \eps$,
$w(\varphi \psi) = w(\varphi) \cdot w(\psi) \geq v(\varphi(\eps))
\bullet \gam = \eps$, therefore $\psi, \varphi \psi \in A$.
Further we obtain 
$$(\varphi \psi)(\gam) = \varphi(\gam)
\bullet \psi(\gam^2) + \gam = \varphi(\gam) \bullet \psi(\gam^2) =
\varphi(\gam) \bullet (\varphi(\gam)^{- 1} + \gam^2) \geq $$
$$\eb(\varphi(\gam)) + \varphi(\gam)) \bullet \gam^2 = 
1_{\gam^2} + \gam = 1_\gam,$$ 
and hence $w(1 - \varphi \psi)) \geq
v(1_\gam - (\varphi \psi)(\gam)) \geq v(1_\gam - 1_\gam) = \gam$.
Thus $w(\varphi (1 - \varphi \psi)) = w(\varphi) w(1 - \varphi \psi)
\geq \gam^{- 1} \bullet \gam = \eps$, therefore $\varphi (1 - \varphi
\psi) \in A$ as desired. 
\ep

\begin{co}
 Let $R$ be a superrigid directed cr-qsring, with the associated 
commutative ring extension $A \sse B$. The necessary and sufficient
condition for $R$ to be a Pr\" ufer cr-qsring is that for all $\varphi \in B$ 
there exists $\psi \in B$ such that $\psi(\gam^2) \leq \varphi(\gam)^{- 1}$, 
where $\gam := w(\varphi)_- = v(\varphi(\eps))^{- 1}$. 
\end{co}

\bp
The statement is a consequence of Lemma 6.6 and 1.4. (P 2).
\ep 

Notice that, by constrast with Pr\" ufer ring extensions, the
property of a superrigid directed cr-qsring to be Pr\" ufer
is not an elementary property.

A particular class of Pr\" ufer cr-qsrings is obtained by
considering the following completeness property
for directed cr-qsrings.

\begin{de}
A directed cr-qsring $R$, with $\Lam = \Ep(R)$, is said to be
{\em spherically complete} (for short, {\em complete})
if for every element $x \in R$, 
the following equivalent conditions are satisfied.
\ben
\item[\rm (1)] There exists a map $\varphi : \Lam \lra R$ 
such that $\epl(\varphi(\alpha)) = \alpha$ for all 
$\alpha \in \Lam$, $\varphi(\alpha) \leq \varphi(\beta)$ 
provided $\alpha \leq \beta$, and $\varphi(\epl(x)) = x$.
\item[\rm (2)] There exists a coherent family of balls
with center $x$ and arbitrary radia, i.e. a map $\psi : \Lam \lra R$
such that $d(\psi(\alpha), x) = |\alpha|$ for all 
$\alpha \in \Lam$, in particular, $\psi(\eps) = x$, and
$\psi(\alpha) \leq \psi(\beta)$ provided $\alpha \leq \beta$.
\een
{\em [For $\varphi$ satisfying (1), define $\psi(\alpha) =
\varphi(\alpha \bullet \epl(x))$ for $\alpha \in \Lam$.
Conversely, for $\psi$ satisfying (2), define $\varphi(\alpha)
= \psi(\frac{\alpha}{\epl(x)}) = \psi((\frac{\alpha}{\epl(x)})_+) + 
\alpha$ for $\alpha \in \Lam$.]}
\end{de}

\begin{co}
A superrigid directed cr-qsring $R$ is Pr\" ufer whenever it
is complete. 
\end{co}

\bp
Let $A \sse B$ be the ring extension associated to $R$.
Assuming that $R$ is complete, it follows that for all
$\alpha \in \Lam := \Ep(R)$, the canonical map 
$\pi_\alpha : B \lra R_\alpha, \varphi \mapsto
\wt{\varphi}(\alpha) := \varphi(\alpha_+) + \alpha$ is
onto. In particular, for every $\varphi \in B$, with
$\gam := w(\varphi)_- = v(\varphi(\eps))^{- 1}$, 
there is $\psi \in B$ such that $\psi(\gam^3) =
\pi_{\gam^3}(\psi) = \varphi(\gam)^{- 1} \in R_{\gam^3}$,
and hence $\psi(\gam^2) \leq \psi(\gam^3) = \varphi(\gam)^{- 1}$.
Consequently, $R$ is a Pr\" ufer cr-qsring by Corollary 6.7.
\ep


\subsection{Pr\" ufer-Manis commutative regular quasi-semirings}

$\q$ Let $R$ be a Pr\" ufer cr-qsring. According to 1.4, 
the canonical embedding $M^\ast \lra \wh{M^\ast}$ of the
Abelian $l$-group $M^\ast$ into its commutative $l$-monoid
extension $\wh{M^\ast}$, as defined in 1.1, extends to a morphism
of $sl$-monoids $\wh{\mathfrak{w}} : M \lra \wh{M^\ast}$, defined by
$\wh{\mathfrak{w}}(I)(\mathbf{a}) := I + \mathbf{a}$ for $I \in M,
\mathbf{a} \in M^\ast_+$, while the map $\mathfrak{w} : B \lra \wh{M^\ast},
\varphi \mapsto \wh{\mathfrak{w}}(A \varphi)$, induced by $\wh{\mathfrak{w}}$, 
is the Pr\" ufer-Manis $l$-valuation associated to the Pr\" ufer
ring extension $A \sse B$. Notice that $A_{\mathfrak{w}} = A$,
supp$(\mathfrak{w}) = \cap_{\mathbf{a} \in M^\ast_+} \mathbf{a}$ 
is the conductor of $A$ in $B$, and the reduced ring extension 
$A/{\rm supp}(\mathfrak{w}) \sse B/{\rm supp}(\mathfrak{w})$
is Pr\" ufer. The relation between the $l$-valuations 
$w : B \lra \wh{\Lam}$ and $\mathfrak{w} : B \lra \wh{M^\ast}$, 
as well as that between the morphisms $\wh{w} : M \lra \wh{\Lam}$ and
$\wh{\mathfrak{w}} : M \lra \wh{M^\ast}$, are described by Lemma 1.20.
In particular, $M^\ast \cong \wh{w}(M^\ast) = \wh{w}(M)^\ast =
\wh{w}(M) \cap \Gam$, where $\Gam$ is the convex hull of
$\wh{w}(M)^\ast$ in $\Lam$. Using Lemma 1.20.(3), we introduce 
the following class of Pr\" ufer cr-qsrings.

\begin{de}
A Pr\" ufer cr-qsring $R$ is called {\em Pr\" ufer-Manis} if
$\Lam = \Ep(R)$ is the convex hull of
$\wh{w}(M)^\ast \cong M^\ast$, in particular, $\wh{w}(M) \cap \Lam =
\wh{w}(M)^\ast$ and $w$ is equivalent with the Pr\" ufer-Manis 
$l$-valuation $\mathfrak{w}$ associated to the Pr\" ufer extension 
$A \sse B$. It is called {\em Pr\" ufer-Manis in a strong sense}
if $M^\ast \cong \wh{w}(M)^\ast = \Lam$.
\end{de}

\begin{lem}
Let $R$ be a superrigid directed cr-qsring. If $R$ is complete
and $M^\ast \cong \wh{w}(M)^\ast = \Lam$ then $R$ is Pr\" ufer-Manis
in a strong sense and lff, in particular, median.
\end{lem}

\bp
By Lemma 6.9, $R$ is a Pr\" ufer cr-qsring, and hence 
Pr\" ufer-Manis in a strong sense cf. Definition 6.10.

To show that $R$ is lff, let $x, y \in R$ be such that
$x - y \in \Lam = \Ep(R)$. According to Lemma 5.4 and 
Definition 5.5, we have to show that $x \vee y \neq \infty$,
i.e. $x \leq z, y \leq z$ for some $z \in R$.
Put $\epl(x) = \alpha, \epl(y) = \beta$, so 
$x - y = v(x - y) = \epl(x - y) = \alpha+ \beta = \alpha \wedge \beta$
by assumption. Let $\gam := \alpha_+ \vee \beta_+$. 
By completeness, it follows that there exist 
$\varphi, \psi \in B := B(R)$ such that 
$\wt{\varphi}(\alpha) := \varphi(\gam) + \alpha = 
x, \wt{\psi}(\beta) := \psi(\gam) + \beta = y$,
i.e. $x \leq \varphi(\gam), y \leq \psi(\gam)$.
Consequently, $(\varphi - \psi)(\gam) = \varphi(\gam) -
\psi(\gam) \geq x - y = \alpha \wedge \beta$, therefore
$w(\varphi - \psi) \geq w(\varphi - \psi)(\gam) \geq
\alpha \wedge \beta$, i.e. $\varphi - \psi \in I_{\alpha \wedge \beta}$.
On the other hand, since $\wh{w}(M)^\ast = \Lam$ by assumption,
it follows that $I_\alpha, I_\beta \in M^\ast$, and
$I_{\alpha \wedge \beta} = I_\alpha + I_\beta \in M^\ast$.
Thus $\varphi - \psi = \rho + \mu$ for some $\rho \in I_\alpha,
\mu \in I_\beta$. Setting $\theta := \varphi - \rho = \psi +
\mu$, we deduce that $\theta(\gam) = \varphi(\gam) - \rho(\gam) \geq
x, \theta(\gam) = \psi(\gam) + \mu(\gam) \geq y$, and hence
$x \vee y \neq \infty$ as desired. 

As $R$ is lff, it is median by Remarks 5.6.(1).
\ep

Denote by $\mathfrak{R}$ the category of Pr\" ufer-Manis 
cr-qsrings, with morphisms $f:R \lra R'$ 
satisfying the condition from Remark 6.4 : $\Lam' := \Ep(R')$
is the convex hull of $f(\Lam)$. The correspondence above
$R \mapsto (A \sse B)$ extends as in Remark 6.4 to a covariant
functor $\cB : \mathfrak{R} \lra \mathfrak{P}$, where $\mathfrak{P}$ denotes
the category of Pr\" ufer extensions with morphisms
$g : (A \sse B) \lra (A' \sse B')$ satisfying the 
supplementary natural condition : the Abelian $l$-group
$M_{B'/A'}^\ast$ of invertible $A'$-submodules of $B'$ is
the convex hull of the image $\{A' g(\mathbf{a})\,|\,\mathbf{a} \in 
M_{B/A}^\ast\}$ of $M_{B/A}^\ast$ through the morphism $g$.

Moreover we obtain

\begin{lem}
The functor $\cB : \mathfrak{R} \lra \mathfrak{P}$ takes 
values in the full subcategory $\mathfrak{CP}$ of $\mathfrak{P}$ 
consisting of the {\em complete Pr\" ufer extensions}, i.e. those
Pr\" ufer extensions $A \sse B$ for which the canonical
morphisms $A \lra \displaystyle\lim_{\stackrel{\lla}
{\mathbf{a} \in (M_{B/A}^\ast)_+}} A/\mathbf{a}$ and 
$B \lra \displaystyle\lim_{\stackrel{\lla}{\mathbf{a} 
\in (M_{B/A}^\ast)_+}} B/\mathbf{a}$ are isomorphisms. 
\end{lem}

\bp
Let $R$ be a Pr\" ufer-Manis cr-qsring,
with $\cB(R) = (A \sse B)$, its associated 
Pr\" ufer ring extension, and $M := M_{B/A}$ the $sl$-monoid
of f.g. $A$-submodules of $B$. Since $\Lam := \Ep(R)$ is the convex
hull of $\wh{w}(M^\ast) = \wh{w}(M) \cap \Lam \cong M^\ast$, 
it follows that $A \cong \displaystyle\lim_{\stackrel{\lla}
{\mathbf{a} \in M^\ast_+}} A/I_{\wh{w}(\mathbf{a})}$,
$B \cong \displaystyle\lim_{\stackrel{\lla}
{\mathbf{a} \in M^\ast_+}} B/I_{\wh{w}(\mathbf{a})}$,
where $I_\alpha := \{\varphi \in A\,|\,w(\varphi) \geq \alpha\}$
for $\alpha \in \Lam_+$. It remains to recall that 
$I_{\wh{w}(\mathbf{a})} = \mathbf{a}$ for all 
$\mathbf{a} \in M^\ast_+$.
\ep


\section{Directed commutative regular quasi-semirings associated
to Pr\" ufer extensions}

$\q$ Extending \cite{Res}, we construct in this section a covariant
functor from the category of Pr\" ufer extensions to the category
of superrigid directed cr-qsrings, whose image turns out to be
equivalent with the category of complete Pr\" ufer extensions. 

Let $A \sse B$ be a Pr\" ufer ring extension, and $\Lam := M^\ast_{B/A}$
be the multiplicative Abelian $l$-group of the invertible 
$A$-submodules of $B$, with $A$
as neutral element, and $\Lam_+ := \{\alpha \in \Lam\,|\,\alpha \sse A\}$,
the $l$-monoid of nonnegative elements of $\Lam$. 
As shown in 1.4, the map $w : B \lra \wh{\Lam}$, defined by
$w(x)(\alpha) = \alpha + A x$ for $x \in B, \alpha \in \Lam_+$,
is the Pr\" ufer-Manis $l$-valuation associated to the 
Pr\" ufer ring extension $A \sse B$.

Let $R = R_{B/A} := \sqcup_{\alpha \in \Lam} B/\alpha$ be the 
disjoint union of the factor $A$-modules $B/\alpha$
for $\alpha$ ranging over $\Lam$. Notice that $R$ is 
a singleton if and only if $A = B$. 

The injective map $\Lam \lra R, \alpha \mapsto 0\,{\rm mod}\,\alpha$ 
identifies $\Lam$ with a subset of $R$, while the ring structure of 
$B$ induces on the residue set $R$ an algebra 
$(R, +, \bullet, -, ^{- 1}, \eps)$ of signature 
$(2, 2, 1, 1, 0)$ as follows.

For $\mathbf{x} := x\,{\rm mod}\,\alpha,\, \mathbf{y} := 
y\,{\rm mod}\,\beta \in R$, we define
$$\mathbf{x} + \mathbf{y} := (x + y)\,{\rm mod}\,(\alpha + \beta),$$
$$- \mathbf{x} := (-x)\,{\rm mod}\,\alpha,$$
and
$$\mathbf{x} \bullet \mathbf{y} := x y\,{\rm mod}\,\gam,$$
with 
$$\gam = \alpha \beta + x \beta + y \alpha = \alpha v(\mathbf{y})
+ \beta v(\mathbf{x}),$$
where $v(\mathbf{x}) := w(x)(\alpha) = \alpha + A x \in \Lam$.

We see at once that the operations above are well defined,
making $(R, +, -)$ a regular commutative semigroup with
$\Ep = (\Lam, +)$ as its semilattice of idempotents, and
$(R, \bullet)$ a commutative semigroup with 
$\Eb = \{1\,{\rm mod}\,\alpha\,|\,\alpha \in \Lam_+\}$
as its semilattice of idempotents, canonically isomorphic
through the map $1\,{\rm mod}\,\alpha \mapsto \alpha$ with
the semilattice $(\Lam_+, +)$. Put $\eps := 0\,{\rm mod}\,A =
1\,{\rm mod}\,A$, and notice that $\Ep \cap \Eb = \{\eps\}$.
Moreover the commutative semigroup $(R, \bullet)$ is regular
with the quasi-inverse $\mathbf{x}^{- 1}$ of any element
$\mathbf{x} = x\,{\rm mod}\,\alpha$ defined as follows.
Since $v(\mathbf{x}) = \alpha + A x \in \Lam$, we obtain
$A = v(\mathbf{x})^{- 1} (\alpha + A x) = v(\mathbf{x})^{- 1} \alpha
+ v(\mathbf{x})^{- 1} A x$, therefore there exists 
$z \in v(\mathbf{x})^{- 1}$ such that $x z \equiv 1\,{\rm mod}\,
v(\mathbf{x})^{-1} \alpha$. One checks easily that the
element 
$$\mathbf{x}^{- 1} := z\,{\rm mod}\,\alpha v(\mathbf{x})^{- 2}$$
is well defined and is the required quasi-inverse of $\mathbf{x}$.  

With the notation from Section 2, we get $\epl(\mathbf{x}) = \alpha,
\eb(\mathbf{x}) = 1\,{\rm mod}\,\alpha v(\mathbf{x})^{- 1}$ for
$\mathbf{x} = x\,{\rm mod}\,\alpha \in R$. The canonical
partial order $\leq$ on $R$ is given by
$$x\,{\rm mod}\,\alpha \leq y\,{\rm mod}\,\beta \Llra
\beta \sse \alpha\,{\rm and}\,x \equiv y\,{\rm mod}\,\alpha.$$
The axioms from Definition 2.3 are easily verified, so $R$,
with the operations as defined above, is a cr-qsring.
Moreover, as $\eps \leq \mathbf{x}$ for all $\mathbf{x} \in \Eb$,
the equivalent conditions of Proposition 3.1. are satisfied,
and hence the cr-qsring $R$ is directed, with the meet-semilattice 
operation with respect to the partial order $\leq$ given by
$$x\,{\rm mod}\,\alpha \wedge y\,{\rm mod}\,\beta =
x\,{\rm mod}\,(\alpha + \beta + A (x - y)).$$
In addition, $R$ is superrigid and, according to Section 4, 
a $\Lam$-metric space with the $\Lam$-valued distance map 
$d : R \times R \lra \Lam_+$ defined by
$$d(x\,{\rm mod}\,\alpha, y\,{\rm mod}\,\beta) = 
\alpha \beta (\alpha + \beta + A (x - y))^{-2}.$$

Proceeding as in 6.1, we associate to the superrigid
directed cr-qsring $R$ above the commutative ring 
$\wh{B} := B(R)$ with carrier $\displaystyle\lim_{\stackrel{\lla}
{\alpha \in \Lam_+}} B/\alpha$ and the $l$-valuation
$\mathbf{w} : \wh{B} \lra \wh{\Lam}$, whose $l$-valuation
subring is $\wh{A} := \displaystyle\lim_{\stackrel{\lla}
{\alpha \in \Lam_+}} A/\alpha$. Composing $\mathbf{w}$ with 
the canonical ring morphism $B \lra \wh{B}$, we obtain
the Pr\" ufer-Manis $l$-valuation $w : B \lra \wh{\Lam}$ 
associated to the Pr\" ufer ring extension $A \sse B$,
in particular, $\cM_{\mathbf{w}}^\ast = \cM_w^\ast = \Lam$ and
$\mathbf{w}$ is a Manis $l$-valuation.

\begin{lem}
 Let $A \sse B$ be a Pr\" ufer commutative ring extension,
$\Lam := M^\ast_{B/A}$ the Abelian $l$-group of invertible
$A$-submodules of $B$, and $R := R_{B/A}$ the associated 
superrigid directed cr-qsring. Then $R$ is complete 
and Pr\" ufer-Manis in a strong sense.

In particular, $R$ is lff, and hence median, with 
the median operation $m : R^3 \lra R$ defined by
$$m(\mathbf{x}_1, \mathbf{x}_2, \mathbf{x}_3) = 
\vee_{1 \leq i < j \leq 3} (\mathbf{x}_i \wedge \mathbf{x}_j) = 
x\,{\rm mod}\,\alpha\,{\rm for}\,\mathbf{x}_i = x_i\,
{\rm mod}\,\alpha_i, i = 1, 2, 3,$$
where
$$\alpha = \cap_{1 \leq i < j \leq 3} (\alpha_i + \alpha_j + A (x_i - x_j)),$$
and $x\,{\rm mod}\,\alpha$ is uniquely determined by the
conditions
$$x\,\equiv\,x_i\,\equiv\,x_j\,{\rm mod}\,(\alpha_i + \alpha_j +
A(x_i - x_j)),\,1 \leq i < j \leq 3.$$
\end{lem}

\bp
To check that the directed cr-qsring $R$ is complete, let
$\mathbf{x} = x\,{\rm mod}\,\gam \in R$. For any $y \in B$
satisfying $y \equiv x\,{\rm mod}\,\gam$, the map
$\varphi_y : \Lam \lra R$, defined by $\varphi_y(\alpha) :=
y\,{\rm mod}\,\alpha$, satisfies the condition (1) from 
Definition 6.8. as desired. In particular, by Corollary 6.9,
$R$ is a Pr\" ufer cr-qsring, i.e. the commutative ring
extension $\wh{A} \sse \wh{B}$ associated to $R$ is Pr\" ufer.
Moreover $R$ is Pr\" ufer-Manis in a strong sense (cf.
Definition 6.10.) since $\Lam = \wh{w}(M_{B/A})^\ast \sse
\wh{\mathbf{w}}(M_{\wh{B}/\wh{A}})^\ast \sse \Lam$,
and hence $M_{\wh{B}/\wh{A}}^\ast \cong \wh{\mathbf{w}}
(M_{\wh{B}/\wh{A}})^\ast = \Lam$ as required. By
Lemma 6.11, $R$ is lff, and hence median.
\ep

With the notation from 6.3, the correspondence above
$(A \sse B) \mapsto R_{B/A}$ extends to a covariant
functor $\cR : \mathfrak{P} \lra \mathfrak{R}$. By
Lemma 7.1, the functor $\cR$ takes values in the full
subcategory $\mathfrak{CR}$ of $\mathfrak{R}$ 
consisting of those superrigid directed cr-qsrings
which are complete and Pr\" ufer-Manis 
in a strong sense.

The relationship of adjunction between $\cR$ and the functor
$\cB$ as defined in 6.3. is described by the next statement.

\begin{te}
\ben
\item[\rm (1)] The covariant functor $\cR : \mathfrak{P} \lra \mathfrak{R}$ 
is the left adjoint of the covariant functor 
$\cB : \mathfrak{R} \lra \mathfrak{P}$.
\item[\rm (2)] The endofunctor $\cB \circ \cR : \mathfrak{P} \lra \mathfrak{P}$,
factorizing through $\mathfrak{CP}$, sends a Pr\" ufer extension 
$A \sse B$ to its completion $\wh{A} := \displaystyle\lim_{\stackrel{\lla}
{\alpha \in (M_{B/A}^\ast)_+}} A/\alpha \sse \wh{B} := 
\displaystyle\lim_{\stackrel{\lla}{\alpha \in (M_{B/A}^\ast)_+}} B/\alpha$. 
\item[\rm (3)] By restriction, $\cR$ and $\cB$ yield inverse
equivalences of the categories $\mathfrak{CP}$ and $\mathfrak{CR}$.
\een
\end{te}

\bp
We define two natural transformations $\mathfrak{r} : \cR \circ
\cB \lra 1_{\mathfrak{R}}, \mathfrak{p} : 1_{\mathfrak{P}} \lra
\cB \circ \cR$ as follows. For each $R$ in $\mathfrak{R}$,
$\cR(\cB(R))$ is an object of $\mathfrak{CR}$,
canonically identified with a subobject of $R$ in $\mathfrak{R}$. Let 
$\mathfrak{r}_R : \cR(\cB(R)) \lra R$ denote this embedding,
and notice that $\mathfrak{r}_R$ is an isomorphism if and
only if $R$ is in $\mathfrak{CR}$. On the other hand, for
each $(A \sse B)$ in $\mathfrak{P}$, $\cB(\cR(A \sse B)) = 
(\wh{A} \sse \wh{B})$ is an object of $\mathfrak{CP}$, the
completion of the Pr\" ufer extension $A \sse B$. Let
$\mathfrak{p}_{A \sse B} : (A \sse B) \lra (\wh{A} \sse \wh{B})$
be the canonical completion morphism in $\mathfrak{P}$, and
notice that $\mathfrak{p}_{A \sse B}$ is an isomorphism if
and only if $(A \sse B)$ is in $\mathfrak{CP}$. It follows
easily that the families of morphisms $\mathfrak{r}_R$
and $\mathfrak{p}_{A \sse B}$, for $R$, $(A \sse B)$ ranging over
the objects of the categories $\mathfrak{R}$ and $\mathfrak{P}$
respectively, are natural transformations. Moreover it
follows that $\cR \mathfrak{p} = (\mathfrak{r} \cR)^{- 1} :
\cR \lra \cR \circ \cB \circ \cR$ and 
$\mathfrak{p} \cB = (\cB \mathfrak{r})^{- 1} :
\cB \lra \cB \circ \cR \circ \cB$ are natural isomorphisms, 
so, in particular, the {\em counit-unit equations} 
$1_{\cR} = \mathfrak{r} \cR \circ \cR \mathfrak{p},
1_{\cB} = \cB \mathfrak{r} \circ \mathfrak{p} \cB$ are
satisfied. We conclude that $(\cR, \cB)$ form an {\em adjoint pair}
with {\em counit} $\mathfrak{r}$ and {\em unit} $\mathfrak{p}$,
inducing by restriction inverse equivalences of the full
subcategories $\mathfrak{CR}$ and $\mathfrak{CP}$. 

The associated {\em monad} in the category $\mathfrak{P}$ 
is the triple $(T, \mathfrak{p}, \mu)$,
where the endofunctor $T : \mathfrak{P} \lra \mathfrak{P}$
is given by $T = \cB \circ \cR$, the {\em unit of the monad} is
just the unit $\mathfrak{p} : 1_{\mathfrak{P}} \lra T$ of 
the adjunction, and the {\em multiplication}
$\mu : T^2 = T \circ T \lra T$
is the natural isomorphism $\cB \mathfrak{r} \cR$.
Dually, the endofunctor $\cR \circ \cB : \mathfrak{R} \lra
\mathfrak{R}$ together with the natural transformation
$\mathfrak{r} : \cR \circ \cB \lra 1_{\mathfrak{R}}$, as
{\em counit}, and the natural isomorphism 
$\cR \mathfrak{p} \cB : \cR \circ \cB \lra (\cR \circ \cB)^2$, 
as {\em comultiplication}, defines a {\em comonad} in 
the category $\mathfrak{R}$. 
\ep 

\begin{rems} \em 
\ben
\item[\rm (1)] Let us denote by $\mathfrak{P}_0$ the 
full subcategory of $\mathfrak{P}$ consisting of those 
Pr\" ufer extensions $A \sse B$ for which $M_{B/A}^\ast$ 
is nontrivial and totally ordered, i.e. $B \neq A$ 
and $B \sm A$ is multiplicatively closed 
(cf. Corollary 1.22). The objects of its full subcategory 
$\mathfrak{CP}_0 := \mathfrak{P}_0 \cap \mathfrak{CP}$ are
the complete valued fields. On the other hand, we denote by 
$\mathfrak{R}_0$ the full subcategory of $\mathfrak{R}$
consisting of those $R$ for which $\Lam = \Ep(R)$ is
nontrivial and totally ordered. Let $\mathfrak{CR}_0 :=
\mathfrak{R}_0 \cap \mathfrak{CR}$. By restriction, we
obtain an adjoint pair $(\cR_0 : \mathfrak{P}_0 \lra
\mathfrak{R}_0, \cB_0 : \mathfrak{R}_0 \lra \mathfrak{P}_0)$, 
inducing inverse equivalences of $\mathfrak{CP}_0$ 
and $\mathfrak{CR}_0$.
\item[\rm (2)] Given a family $(A_i \sse B_i)_{i \in I}$ 
of Pr\" ufer extensions, let $A := \prod_{i \in I} A_i,
B := \prod_{i \in I} B_i$, with the canonical projections
$p_i : B \lra B_i$, satisfying $p_i(A) = A_i, i \in I$. 
By \cite[I, Proposition 5.17.]{Kne}, the ring extension 
$A \sse B$ is Pr\" ufer, so the projections

$p_i : (A \sse B) \lra (A_i \sse B_i), i \in I$, are morphisms
in $\mathfrak{P}$. The commutative $sl$-monoid $M_{B/A}$ of
f.g. $A$-submodules of $B$ is identified through the 
canonical morphism $M_{B/A} \lra \prod_{i \in I} M_{B_i/A_i}, 
\alpha \mapsto (p_i(\alpha))_{i \in I}$ to the subdirect
product consisting of those families 
$\alpha := (\alpha_i)_{i \in I} \in \prod_{i \in I} M_{B_i/A_i}$
for which there is a bound $n_\alpha \in \N$ such that
the $A_i$-submodule $\alpha_i$ of $B_i$ admits a system of at most
$n_\alpha$ generators for all $i \in I$. By restriction,
the Abelian $l$-group $\Lam := M_{B/A}^\ast$ is identified
with a subdirect product of the family of Abelian $l$-groups
$\Lam_i := M_{B_i/A_i}^\ast, i \in I$, consisting of
those $\alpha := (\alpha_i)_{i \in I} \in \prod_{i \in I} \Lam_i$
satisfying the boundedness condition above.
In particular, $\Lam = \prod_{i \in I} \Lam_i$ whenever either
the index set $I$ is finite or there is a bound $n \in \N$ such
that each invertible $A_i$-submodule of $B_i$ admits a system of
at most $n$ generators for all $i \in I$.

By Lemma 1.3 and Remark 1.5, the commutative $l$-monoid
$\wh{\Lam}$ is identified with the subdirect product
of the family $(\wh{\Lam_i})_{i \in I}$ consisting of
those $\varphi \in \prod_{i \in I} \wh{\Lam_i} \cong
\wh{\prod_{i \in I} \Lam_i}$ for which $\varphi \wedge 
\alpha \in \Lam$ for all $\alpha \in \Lam$. 

Let $R_i := \cB(A_i \sse B_i), i \in I, R := \cB(A \sse B)$.
$R$ is identified with the subdirect product of the
family $(R_i)_{i \in I}$ consisting of those $\mathbf{x} \in
\prod_{i \in I} R_i$ for which $\epl(\mathbf{x}) \in \Lam$.  
Applying the functor $\cB : \mathfrak{R} \lra \mathfrak{P}$
to the objects $R_i, i \in I$, and $R$ of the category $\mathfrak{CR}$,
we obtain the completions $\cB(R_i) := (\wh{A_i} \sse \wh{B_i}),
i \in I, \cB(R) := (\wh{A} \sse \wh{B})$ of the Pr\" ufer
extensions $A_i \sse B_i, i \in I$, and $A \sse B$ respectively.
The complete Pr\" ufer extension $\wh{A} \sse \wh{B}$ is identified
with the subdirect product of the family $\wh{A_i} \sse \wh{B_i},
i \in I$, consisting of those elements $(x_i)_{i \in I} \in 
\prod_{i \in I} \wh{B_i}$ for which $(w_i(x_i))_{i \in I}
\in \wh{\Lam}$, where $w_i : \wh{B_i} \lra \wh{\Lam_i}, i \in I$,
are the Pr\" ufer-Manis $l$-valuations associated to the 
complete Pr\" ufer extensions $\wh{A_i} \sse \wh{B_i}, i \in I$.
In particular, $R \cong \prod_{i \in I} R_i$,
and $(\wh{A} \sse \wh{B}) \cong \prod_{i \in I} (\wh{A_i} \sse
\wh{B_i})$ whenever either $I$ is finite or there is $n \in \N$
such that for all $i \in I$, each $\alpha \in \Lam_i$ admits 
a system of at most $n$ generators.

Notice that, in general, the Pr\' ufer extension $A \sse B$ together
with the projections $p_i : (A \sse B) \lra (A_i \sse B_i), i \in I$,
is not a direct product in the category $\mathfrak{P}$ since,
given an arbitrary family $f_i : (C \sse D) \lra (A_i \sse B_i), i \in I$,
of morphisms in $\mathfrak{P}$, the canonical morphism of 
Pr\" ufer extensions $f : (C \sse D) \lra (A \sse B)$ is not 
necessarily a morphism in $\mathfrak{P}$. However 
$f$ is in $\mathfrak{P}$ whenever the index set $I$ is finite, and 
hence the category $\mathfrak{P}$ has finite products. Similarly, 
the category $\mathfrak{R}$ has finite products, and for each 
finite family $(R_i)_{i \in I}$ of objects of $\mathfrak{R}$, 
it follows by the continuity of the right adjoint functor 
$\cB : \mathfrak{R} \lra \mathfrak{P}$ that 
$\cB(\prod_{i \in I} R_i) \cong \prod_{i \in I} \cB(R_i)$.
In particular, the equivalent full subcategories 
$\mathfrak{CP} \sse \mathfrak{P}$ and $\mathfrak{CR} \sse \mathfrak{R}$
are both closed under finite products.
\item[\rm (3)] Assume that $A$ is a Dedekind (in particular,
Pr\" ufer) domain with its quotient field $B \neq A$. Let
us denote by $\cP := \cP(A)$ the nonempty set of the non-null
prime (maximal) ideals of $A$. For every $\mathfrak{p} \in \cP$, we denote 
by $w_{\mathfrak{p}} : B \lra \Z \cup \{\infty\}$ the 
corresponding discrete valuation with
valuation ring $A_{\mathfrak{p}}$ and maximal ideal 
$\mathfrak{p} A_{\mathfrak{p}}$. With the notation above,
we have $M := M_{B/A} = \Lam \cup \{ \{0\} \}$, where $\Lam :=
M^\ast \cong \Z^{(\cP)}$ is the free Abelian multiplicative group
generated by $\cP$, with the canonical lattice operations. The
corresponding commutative $l$-monoid extension $\wh{\Lam}$ of
$\Lam$ is identified, as in Example 1.1, with the collection of
those formal products 
$\prod_{\mathfrak{p} \in \cP} \mathfrak{p}^{n_{\mathfrak{p}}}$
with $n_{\mathfrak{p}} \in \Z \cup \{\infty\}$ for which the set
$\{\mathfrak{p}\,|\,n_{\mathfrak{p}} < 0\}$ is finite. The
Pr\" ufer-Manis $l$-valuation $w : B \lra \wh{\Lam}$ associated
to the Pr\" ufer extension $A \sse B$ sends any $0 \neq x \in B$
to $A x = \prod_{\mathfrak{p} \in \cP} \mathfrak{p}^{w_{\mathfrak{p}}(x)}
\in \Lam$, while $w(0) = \om$, so $M \cong \cM_w$. Let 
$R := \cR(A \sse B), R_{\mathfrak{p}} := \cR(A_{\mathfrak{p}}
\sse B)$ for $\mathfrak{p} \in \cP$, so $\Lam_{\mathfrak{p}} :=
\Ep(R_{\mathfrak{p}}) \cong (\Z, +), \Lam = \Ep(R) = 
\oplus_{\mathfrak{p} \in \cP} \Lam_{\mathfrak{p}}$. By the 
{\em approximation theorem in Dedekind domains} 
\cite[VII, 4, Proposition 2]{Bou} it follows
that the canonical subdirect representation 
$$R \lra \prod_{\mathfrak{p} \in \cP} R_{\mathfrak{p}},\, x \,{\rm mod}\,\alpha
\in R_\alpha := B/\alpha \mapsto 
(x\,{\rm mod}\,\mathfrak{p}^{w_{\mathfrak{p}}(\alpha)})_
{\mathfrak{p} \in \cP},$$
for $\alpha := \prod_{\mathfrak{p} \in \cP} \mathfrak{p}
^{w_{\mathfrak{p}}(\alpha)} \in \Lam$, identifies $R$ with the
substructure of the product $\prod_{\mathfrak{p} \in \cP} R_{\mathfrak{p}}$
consisting of those families
$\mathbf{x} := (\mathbf{x}_{\mathfrak{p}})_{\mathfrak{p} \in \cP}$ for which 
$d(\mathbf{x}, \eps := 0\,{\rm mod}\,A) := (d_{\mathfrak{p}}
(\mathbf{x}_{\mathfrak{p}}, \eps_{\mathfrak{p}} := 0\,{\rm mod}\,A_{\mathfrak{p}}))
_{\mathfrak{p} \in \cP} \in \Lam$, i.e. the set  
$\{\mathfrak{p} \in \cP\,|\,\mathbf{x}_{\mathfrak{p}} \neq \eps_{\mathfrak{p}}\}$ 
is finite. In particular, $R \cong \prod_{\mathfrak{p} \in \cP} R_{\mathfrak{p}}$
$\Llra A$ is semilocal, i.e. $\cP$ is finite. As median set, $R$ is
{\em simplicial}, i.e. each cell $[\mathbf{x}, \mathbf{y}] \sse R$ has finitely
many elements. Moreover for each finite subset $X \sse R$ with $|X| \geq 2$, the 
convex closure $[X]$ of $X$ in the median set $R$ is isomorphic to
a finite product of finite simplicial trees : $[X] = \prod_{\mathfrak{p} \in 
\cP} [X_{\mathfrak{p}}] \cong \prod_{|X_{\mathfrak{p}}| \geq 2} [X_{\mathfrak{p}}]$,
where $[X_{\mathfrak{p}}]$ is the subtree of the simplicial tree $R_{\mathfrak{p}}$
generated by the image $X_{\mathfrak{p}}$ of $X$ trough the projection
map $R \lra R_{\mathfrak{p}}$. Applying the functor $\cB : \mathfrak{R} \lra
\mathfrak{P}$ to the object $R = \cR(A \sse B)$ of $\mathfrak{CR}$, we obtain
the complete Pr\" ufer extension $\cB(R) = (\wh{A} \sse \wh{B})$ - the
completion of the Pr\" ufer extension $A \sse B$ - with
$\wh{A} = \prod_{\mathfrak{p} \in \cP} \wh{A_{\mathfrak{p}}}$,
$\wh{B} = \{(x_{\mathfrak{p}})_{\mathfrak{p} \in \cP} \in \prod_{\mathfrak{p}
\in \cP} \wh{B}_{\mathfrak{p}}\,|\,\{\mathfrak{p} \in \cP\,|\,
x_{\mathfrak{p}} \not \in \wh{A_{\mathfrak{p}}}\}\,{\rm is\,finite}\}$-the
ring of {\em restricted adeles} of the Dedekind domain $A$, where 
$\wh{B}_{\mathfrak{p}} (\wh{A_{\mathfrak{p}}})$ denotes the completion of
the field $B$ (the domain $A$) with respect to the valuation $w_{\mathfrak{p}}$.
The Pr\" ufer-Manis $l$-valuation $\mathfrak{w} : \wh{B} \lra \wh{\Lam}$
associated to the Pr\" ufer extension $\wh{A} \sse \wh{B}$, sending
any $x := (x_{\mathfrak{p}})_{\mathfrak{p} \in \cP} \in \wh{B}$ to the 
formal product $\prod_{\mathfrak{p} \in \cP} \mathfrak{p}^{w_{\mathfrak{p}}
(x_{\mathfrak{p}})} \in \wh{\Lam}$, is onto, inducing the isomorphisms
$$M_{\wh{B}/\wh{A}} = \{\wh{A} x\,|\,x \in \wh{B}\} \cong \cM_{\mathfrak{w}}
\cong \wh{\Lam},\, M_{\wh{B}/\wh{A}}^\ast = \{\wh{A} x\,|\,x \in \wh{B}^\ast\}
\cong \wh{B}^\ast/\wh{A}^\ast \cong \Lam.$$
\een
\end{rems}

\end{document}